\documentclass{article}

\usepackage{amscd} 
\usepackage{diagramm,curvesm}
\usepackage{amstex,amssymb}

\numberwithin{equation}{section}

\newfam\msafam
\newfam\msbfam

\font\bfzwei=cmbx12 scaled \magstep1

\font\bfzwei=cmbx12 scaled \magstep1
\font\fivemsa=msam5
\font\fivemsb=msbm5

\font\sc=cmcsc10

\font\sevenmsa=msam7
\font\sevenmsb=msbm7

\font\tenmsa=msam10
\font\tenmsb=msbm10
\textfont\msafam=\tenmsa
\textfont\msbfam=\tenmsb
\scriptfont\msafam=\sevenmsa
\scriptfont\msbfam=\sevenmsb
\scriptscriptfont\msafam=\fivemsa
\scriptscriptfont\msbfam=\fivemsb

\def\bfcm{\fontfamily{cmss}\fontseries{bx}\selectfont}

\newtheorem{theorem}{{\bf{T}}{\fontsize{7pt}{11}\fontshape{n}
\fontseries{bx}\selectfont{$\!\!\!$HEOREM}}}[section]
\newtheorem{proposition}[theorem]{{\bf{P}}{\fontsize{7pt}{11}\fontshape{n}
\fontseries{bx}\selectfont{$\!\!\!$ROPOSITION}}}
\newtheorem{corollary}[theorem]{{\bf{C}}{\fontsize{7pt}{11}\fontshape{n}
\fontseries{bx}\selectfont{$\!\!\!$OROLLARY}}}
\newtheorem{lemma}[theorem]{{\bf{L}}{\fontsize{7pt}{11}\fontshape{n}
\fontseries{bx}\selectfont{$\!\!\!$EMMA}}}
\newtheorem{definition}[theorem]{\sc Definition}
\newtheorem{remark}{\sc Remark}
\newtheorem{example}{\sc Example}
\newenvironment{proof}{\par\noindent{\sc Proof.}}

\newcommand{\ouraddress}{{\sc National Academy of Sciences, Bogolyubov
Institute for Theoretical Physics}
 \nl
 Metrologichna str.\ $14^b$, 252-143 Kiev-143, UKRAINE
 \nl
 e-mail: {\bfcm mmtpitp@@gluk.apc.org}
 \abs\abs
 {\sc Department of Mathematics, Katholieke Universiteit Leuven}
 \nl
 Celestijnenlaan 200 B, B-3001 Leuven-Heverlee, BELGIUM
 \nl
 e-mail: {\bfcm bernhard.drabant@@wis.kuleuven.ac.be}
 \abs
 B.\ D.'s address after January 1997:
 \nl
 {\sc Department of Theoretical Physics, Universidad de Valencia}
 \nl
 Avda.\ Dr.-Moliner 50, E-46100 Burjassot (Valencia), SPAIN
 \nl
 e-mail: {\bfcm drabant@@lie1.ific.uv.es}}

\newcommand{\mult}{\put(0,0){\arc(-5,0){180}}
                    \put(0,-5){\line(0,-1){5}}}
\newcommand{\halfmult}{\put(0,0){\arc(-2.5,0){180}}
                    \put(0,-2.5){\line(0,-1){7.5}}}

\newcommand{\comult}{\put(0,0){\line(0,-1){5}}
                     \put(0,-10){\arc(5,0){180}}}
\newcommand{\halfcomult}{\put(0,0){\line(0,-1){7.5}}
                     \put(0,-10){\arc(2.5,0){180}}}
\newcommand{\lact}{\put(0,0){\line(0,-1){10}}
                   \put(0,0){\arc(-5,0){90}}}
\newcommand{\ract}{\put(0,0){\line(0,-1){10}}
                   \put(0,0){\arc(0,-5){90}}}
\newcommand{\lcoact}{\put(0,0){\line(0,-1){10}}
                     \put(0,-10){\arc(0,5){90}}}
\newcommand{\rcoact}{\put(0,0){\line(0,-1){10}}
                     \put(0,-10){\arc(5,0){90}}}
\newcommand{\morphcirc}[1]{\put(0,0){\circle{8}}
                           \put(0,0){\makebox(0,0){$#1$}}}
\newcommand{\morphoval}[1]{\put(0,0){\oval(20,10)}
                           \put(0,0){\makebox(0,0){$#1$}}}
\newcommand{\braid}{\put(0,10){\curve(0,0, 0,-1, 1,-2, 5,-5, 9,-8, 10,-9,
                    10,-10)}
                    \put(0,0){\curve(0,0, 0,1, 1,2, 4,4)}
                    \put(6,6){\curve(0,0, 3,2, 4,3, 4,4)}}

\newcommand{\figeins}{\put(0,0){\line(0,1){10}}
\put(5,0){\line(0,1){10}}
\put(5,20){\lcoact}
\put(5,24){\morphcirc{f}}
\put(5,28){\line(0,1){2}}
\put(5,40){\lact}
\put(0,40){\line(0,1){10}}
\put(5,44){\morphcirc{\dif_X}}
\put(5,48){\line(0,1){2}}}

\newcommand{\figzwei}{\put(0,0){\line(0,1){10}}
\put(5,0){\line(0,1){10}}
\put(5,20){\lcoact}
\put(5,20){\line(0,1){10}}
\put(5,40){\lact}
\put(0,40){\line(0,1){10}}
\put(5,44){\morphcirc{\dif_Y}}
\put(5,48){\line(0,1){2}}}

\newcommand{\figdrei}{\put(0,0){\line(0,1){10}}
\put(20,0){\line(0,1){10}}
\put(0,20){\mult}
\put(20,20){\lact}
\put(5,20){\braid}
\put(-5,20){\line(0,1){10}}
\put(20,20){\line(0,1){10}}
\put(0,40){\comult}
\put(20,40){\lcoact}
\put(0,40){\line(0,1){10}}
\put(20,44){\morphcirc{\dif_Y}}
\put(20,48){\line(0,1){2}}}

\newcommand{\figvier}{\put(0,0){\line(0,1){10}}
\put(20,0){\line(0,1){2}}
\put(20,6){\morphcirc{f}}
\put(0,20){\mult}
\put(20,20){\lact}
\put(5,20){\braid}
\put(-5,20){\line(0,1){10}}
\put(20,20){\line(0,1){10}}
\put(0,40){\comult}
\put(20,40){\lcoact}
\put(0,40){\line(0,1){10}}
\put(20,44){\morphcirc{\dif_X}}
\put(20,48){\line(0,1){2}}}

\newcommand{\figfuenf}{\put(0,0){\line(0,1){10}}
\put(5,0){\line(0,1){2}}
\put(5,6){\morphcirc{f}}
\put(5,20){\lcoact}
\put(5,20){\line(0,1){10}}
\put(5,40){\lact}
\put(0,40){\line(0,1){10}}
\put(5,44){\morphcirc{\dif_X}}
\put(5,48){\line(0,1){2}}}

\newcommand{\hopfbieins}{\put(0,10){\ract}
\put(20,10){\mult}
\put(0,10){\line(0,1){10}}
\put(5,10){\braid}
\put(25,10){\line(0,1){10}} 
\put(0,30){\rcoact}
\put(20,30){\comult}}

\newcommand{\hopfbizwei}{\put(0,10){\rcoact}
\put(0,20){\ract}}

\newcommand{\crossedeins}{\put(0,0){\line(0,1){10}}
\put(12.5,10){\halfmult}
\put(0,10){\braid}
\put(15,10){\line(0,1){10}}
\put(0,20){\line(0,1){10}}
\put(10,30){\rcoact}
\put(0,50){\line(0,1){10}}
\put(12.5,60){\halfcomult}
\put(0,40){\braid}
\put(15,40){\line(0,1){10}}
\put(0,30){\line(0,1){10}}
\put(10,40){\ract}}

\def\abs{\par\vskip 0.3cm\goodbreak\noindent}
\def\Abs{\par\vskip 1.7cm\goodbreak\noindent}
\def\Aut{{\rm Aut}\,}

\def\BiCD{\operatorname{{\rm Bi}_{\cal C}\text{-\rm Der}}}
\def\C{{\cal C}}
\def\CC{{\mathbb C}}
\def\classification#1#2{\medskip \\ {\em #1:\ } #2}

\def\coim{\operatorname{coim}}
\def\Coim{\operatorname{Coim}}
\def\coker{\operatorname{coker}}
\def\Coker{\operatorname{Coker}}
\def\Der{\operatorname{Der}}
\def\dif{{\rm d}}
\def\Dif{{\rm D}}

\def\DY{\operatorname{\cal DY}}
\def\e{{\rm e}}
\def\E{{\bf 1 \!\! {\rm l}}}
\def\End{\operatorname{End}}
\def\endproof{\lfl$\scriptstyle\blacksquare$}
\def\hhchh{{{}^H_H{\C}^H_H}}
\def\Hom{\operatorname{Hom}}
\def\i{{\rm i}}
\def\I{{\rm I}}
\def\II{{\cal I}}
\def\im{\operatorname{im}}
\def\Im{\operatorname{Im}}
\def\id{{\rm id}}
\def\j{{\rm j}}
\def\JJ{{\cal J}}
\def\keywords#1{\classification{Keywords}#1}
\def\ker{\operatorname{ker}}
\def\Ker{\operatorname{Ker}}
\def\LCD{\operatorname{{\rm L}_{\C}\text{-\rm Der}}}
\def\lfl{\leaders\hbox to 1em{\hss \hss}\hfill}
\def\m{{\rm m}}
\def\nl{\par\noindent}
\def\NN{{\mathbb N}}
\def\Ob{\operatorname{Ob}}
\def\p{{\rm p}}

\def\RCD{\operatorname{{\rm R}_{\C}\text{-\rm Der}}}

\def\scripthhchh{{{}^{\scriptscriptstyle H}_{\scriptscriptstyle H}
 {\C}^{\scriptscriptstyle H}_{\scriptscriptstyle H}}}
\def\uNN{{\underline{\NN}}}

\newcommand\subdif[1]{#1^{\hbox{\scriptsize Diff}}}

\newcommand\uotimes{\mathop{\underline\otimes}}
\newcommand\op{{\rm op}}
\newcommand\Mat[1]{\operatorname{M}(#1)}

\newcommand\Mcat{{\rm Mcat}}

\author{Yuri Bespalov
\and
 Bernhard Drabant\thanks{For parts of this work the author had been 
 supported by a grant from the Research Council of the University of Leuven}}

\title{\bfzwei Differential Calculus in Braided Abelian Categories}

\date{December 1996}

\begin{document}

\baselineskip=14pt
\maketitle

\begin{abstract}
\nl

Braided non-commutative differential geometry is studied.
In particular we investigate the theory of (bicovariant) differential calculi
in braided abelian categories. Previous results on crossed modules and
Hopf bimodules in braided categories are used
to construct higher order bicovariant differential calculi over braided
Hopf algebras out of first order ones. These graded objects are shown to be
braided differential Hopf algebras with universal bialgebra properties.
The article especially extends
Woronowicz's results on (bicovariant) differential calculi to
the braided non-commutative case.
\nl
\classification{1991 Mathematics Subject Classification}
               {16W30, 18D10, 18E10, 17B37}
\end{abstract}
\abs
\section*{Introduction}\label{intro}  
Differential geometry and group theory interact very fruitfully within
the theory of Lie groups. Tangent Lie algebras, invariant differential
forms, infinitesimal representations, principal bundles, gauge theory,
etc.\ emerged out of this interplay. It was discovered by Woronowicz that
many of these geometrically related structures can be generalized to
non-commutative geometry \cite{Wor}. He introduced and studied the
differential calculus on
compact quantum groups or more generally on Hopf algebras over a field
$\Bbbk$ with ${\rm char}\,(\Bbbk)=0$. His theory is built on the base
category of $\Bbbk$-modules with the usual tensor product and the
involutive tensor transposition $\tau:a\otimes b\mapsto b\otimes a$.
In what follows we refer to the conditions in \cite{Wor} as the classical
conditions in contrast to our investigations in the braided case.
There have been a lot of publications along the classical lines of
\cite{Wor}. However, it would be beyond the scope of this introduction to
give an appropriate appreciation to all of them. Nevertheless we would like
to mention three articles \cite{Brz,Mal,Man} besides \cite{Wor} which
particularly influenced our work from the differential geometrical point
of view in a considerable way. The differential graded algebra approach
to quantum groups and the quantum de Rham complexes were studied in
\cite{Mal,Man}. In Manin's work \cite[Proposition 2.6.1]{Man} very
general conditions are found under which an operator ring (of differentials)
over an algebra has a bialgebra structure. The
non-commutative differential calculus discovered by Woronowicz is
indeed of this type. The differential Hopf
algebra structure of the higher order differential calculi of \cite{Wor}
has been unfolded explicitely in \cite{Brz}. Essentially all of these
articles proceed on the assumption that the classical
symmetric conditions are at the bottom of the theory of (bicovariant)
differential calculi.

Braided or quasisymmetric monoidal categories had been initially
investigated by Joyal and Street in \cite{JS}. Footing on their work,
Majid generalized the notion of (tensor) algebras, bi- and Hopf
algebras to categories with a tensor product which has a not necessarily
involutive commutation behaviour expressed by the so-called braiding.
The braiding takes over the r\^ole of the tensor transposition $\tau$.
As standard references to Majid's investigations on this subject
we refer to \cite{Ma1,Ma6} where many foundations of braided
mathematics and in particular of braided geometry can be found.

In our article we are interested in the amalgamation of differential
geometrical methods and braided category theory leading to
{\it braided non-commutative differential geometry}. It extends the
results of \cite{BD2}. Several successful
steps in bringing together both methods have been done for instance in
\cite{BM,IV,Ma4}. We show that the theory of differential calculus of
Lie groups and of Hopf algebras can be generalized to braided
non-commutative differential geometry. We examine (bicovariant)
differential calculi over Hopf algebras in a braided abelian category.
Nowhere in the main results of our article we suppose the objects under
consideration to be
sets. Our investigation is footed completely on categorical language
which allows us to lay bare fundamental structures of the theory of
differential calculus. Working with these key-stones we
are able to build in the theory into the framework of a very general class
of braided abelian categories. Our categorical set-up is general enough
to cover and extend all known examples of categories where braided
differential geometry has already been applied. On the other side it is
rich enough to obtain a theory of braided (bicovariant) differential
calculus including the notations of \cite{Mal,Man,Wor}.
We show that every first order bicovariant differential calculus over a
braided Hopf algebra induces in a natural way a higher order
bicovariant differential calculus. It is a
differential Hopf algebra
with certain universal properties -- this
generalizes results of \cite{Brz,Mal,Man} to braided categories.
Especially, our results can be applied to the study of (bicovariant)
differential calculi of quantum planes \cite{FRT,Ma4,Ma5} and of braided
groups \cite{Ma6} but also for the investigation of differential geometric
structures of, say, degenerate Sklyanin algebras \cite{Ma9} and other
interesting objects in braided categories. In connection with
\cite{BD,Bes} bicovariant differential calculi of cross product Hopf
algebras can be investigated in relation to the differential calculi of
their particular factors. Our results may also contribute to the
investigation of braided gauge theory and quantum bundle theory
\cite{BM}. 

The article is devided into \ref{sec-diff-calc} chapters.
Chapter \ref{prelim} mainly reviews and reformulates definitions and results 
of \cite{BD} which are necessary for the understanding of the subsequent
chapters. The proofs of the particular statements can be found in
\cite{BD}. The isomorphism theorem of Hopf bimodule bialgebras and
bialgebra projections in braided monoidal categories is of special
importance. It initiated our investigations on braided bicovariant
differential calculi. Hopf bimodules or bicovariant bimodules are the
fundamental objects for the definition of bicovariant differential
calculi over a Hopf algebra. The isomorphism theorem is an essential tool 
for the construction of higher order differential calculi.

In Chapter \ref{ideal-section} we fix our categorical basis and state
essential results on it. We describe modules generated by morphisms which
will be useful in particular for the definition and
investigation of differential calculi in Chapter \ref{sec-diff-calc}.
Furthermore we are concerned with factor (co-)algebras and bialgebras
arising from canonical epi-mono decompositions of (co-)algebra and
bialgebra morphisms respectively.
Later in Chapter \ref{grade-alg-com} this decomposition allows us to
derive very naturally the antisymmetric and exterior Hopf algebra objects
needed for the construction of higher order differential calculi.
A further restriction of the properties of the categories under
consideration allows us to use a definition of (co-)ideals and bi-ideals
in such categories and to establish well known results of ring theory
in our more general context. This will be applied in Chapters
\ref{grade-alg-com} and \ref{sec-diff-calc}.

In Chapter \ref{grade-alg-com} we investigate graded categories.
We are using the notion of a coalgebra category \cite{CF}
to study under which conditions graded categories, or more general,
functor categories are monoidal. We introduce a quasitriangular structure
on a coalgebra category and show that the corresponding functor categories
are braided in this case. As a special case we are interested in graded
categories and categories of complexes over $\{0,1\}$ and $\NN_0$.
We investigate (co-)algebras,
bi- and Hopf algebras in such categories. We introduce braided
combinatorics in the gereralized sense of \cite{Ma4,Wor}.
Using the results of Chapters \ref{prelim} and \ref{ideal-section} we
construct analogues of antisymmetric tensor Hopf algebras over an arbitrary
object in the braided category $\C$. For a Hopf bimodule $X$ over a
Hopf algebra $H$ in $\C$ we build the exterior Hopf algebra
of forms $X^{\wedge_H}$. This object is the braided version of
Woronowicz's so-called external algebra $\Gamma^\wedge$ over a bicovariant
bimodule $\Gamma$ \cite{Wor}. It is universal as a bialgebra generated
algebraically by its components of order 0 and 1.

Finally Chapter \ref{sec-diff-calc} is
concerned with differential calculi.
In the sense of \cite{Mal,Man,Wor} we define a differential calculus as a
differential graded algebra which is generated by its 0th component and
the image of it under the differential $\dif_0$.
This generalized point of view is due to \cite{Mal,Man}. From
\cite[Chapters 0 and 2]{Man} and \cite[Theorem 1.2.3]{Mal} one can see
that both first and higher order differential calculi of \cite{Wor} fit
into this scheme. 
Manin calls the differential calculi quantum de Rham complexes.
For any graded bialgebra $\hat B$ in the braided category we define
$\hat B$-left, $\hat B$-right and $\hat B$-bicovariant
differential calculi as the corresponding $\hat B$-comodule differential
calculi.
In particular we examine braided bicovariant (first order) differential
calculi $(X,\dif)$ over a Hopf algebra $H$.
We define the (unique) maximal differential calculus of a given
differential algebra and study its structure. Using comma
category technique we establish a braided Woronowicz construction of a
higher order differential Hopf algebra out of a given first order
bicovariant differential calculus. From this object
we extract the exterior Hopf algebra of differential
forms $(X^{\wedge_H},\dif^\wedge)$ over $H$ as its maximal differential
calculus. $(X^{\wedge_H},\dif^\wedge)$ coincides with the first order
calculus $(X,\dif)$ in the 0th and 1st component and is therefore its
unique higher order extension as a differential calculus.

The appendix provides facts from classical category theory connected with
the present article and \cite{BD}. Mostly we will work with these results
without reference.
\abs
\abs
\section{Preliminaries and Notations}\label{prelim}  
We review results and techniques from
\cite{BD,JS,Ma2,Ma6} essential for our purposes.
We assume that the reader is familiar with the notion of a braided
monoidal category
$({\C},\E,\otimes,\alpha,\rho,\lambda,\Psi)$ \cite{FY,JS,Ma1}.
If not otherwise stated we suppose henceforth that $\C$ is braided and
admits split idempotents.
I.e.\ for every idempotent $e=e^2:X\to X$ in $\C$ there exist morphisms
$\i_e:X_e\to X$, $\p_e:X\to X_e$ in $\C$, such that $e=\i_e\circ\p_e$
and $\id_e=\p_e\circ\i_e$. 
This is not a severe restriction, because every (braided monoidal) category 
can be canonically embedded into a (braided monoidal) category that admits 
split idempotents. Abelian categories automatically admit split idempotents.

The notations of (co-)algebras, bi- and Hopf algebras, (co-)modules
and bi-(co-)modules can be generalized obviously 
to the category $\C$ \cite{Ma1,Ma2,Ma6}.
Loosely speaking we replace in the defining relations of the objects
the tensor transposition $\tau$ by the braiding $\Psi$.
Since we are dealing with braided categories we occasionally use graphical
presentation of morphisms. The graphical calculus turns out to be a very 
convenient tool for dealing with complicated expressions of morphisms.
There are several detailed expositions on this subject
\cite{FY,Lyu,Ma1,Ma2,Yet}. We are using the notation of \cite{BD}.
In the article we will sometimes obtain statements simply by dualizing
previous results. This means that we consider the original outcome in
the opposite category and transpose it to its dual version in the category.
In diagrammatic language this may be interpreted as reversing the arrows.
The main results of \cite{BD} which we need in the sequel are
concerned with Hopf bimodules and crossed modules in braided categories.
Hopf bimodules or bicovariant bimodules in symmetric categories
are the essential objects in the definition of bicovariant differential
calculi. Crossed modules are related to the invariant vector fields of a
bicovariant differential calculus \cite{Wor}. They appear naturally in the
representation theory of quantum groups \cite{Yet}.
In the braided case Hopf bimodules and crossed modules will turn out to be of
similar importance for the investigation of braided bicovariant
differential calculi. We recall the definition of Hopf bimodules and
crossed modules in braided categories \cite{BD}
\abs
\begin{definition}\label{hopf-bimod}
Let $(B,\m,\eta,\Delta,\varepsilon)$ be a bialgebra in $\C$. An object
$(X,\mu_r,\mu_l,\nu_r,\nu_l)$ is called a $B$-Hopf bimodule
if $(X,\mu_r,\mu_l)$ is a $B$-bimodule, and if $(X,\nu_r,\nu_l)$ is a
$B$-bicomodule in the category of $B$-bimodules
where the regular action $\m$ on $B$ and the diagonal
action on tensor products of modules are used.
$B$-Hopf bimodules together with the $B$-bimodule-$B$-bicomodule morphisms
form a category which will be denoted by ${}_B^B{\C}^B_B$.
\end{definition}
\abs
The diagonal action of the tensor product of, say, two $B$-left modules
$X$ and $Y$ is given by
\begin{equation}\label{diag-act}
 \mu_l^{X\otimes Y}=(\mu_l^X\otimes\mu_l^Y)\circ(\id_B\otimes\Psi_{B,X}
 \otimes\id_Y)\circ(\Delta\otimes\id_X\otimes\id_Y)\,.
\end{equation}
Similarly the diagonal action of right modules is defined. 
\abs
\begin{definition}\label{crossed-mod}
A right crossed module $(X,\mu_r,\nu_r)$ over the Hopf algebra $H$ in
$\C$ is an $H$-right module and an $H$-right comodule obeying the
compatibility relations
\begin{equation}\label{crossed-mod1}
\begin{gathered}
(\id_X\otimes m)\circ(\Psi_{H,X}\otimes\id_H)\circ
(\id_H\otimes\nu_r\circ\mu_r)\circ
(\Psi_{X,H}\otimes H)\circ(\id_X\otimes\Delta)\\
=(\mu_r\otimes m)\circ(\id_X\otimes\Psi_{H,H}\otimes\id_H)\circ
(\nu_r\otimes\Delta).
\end{gathered}
\end{equation}
The category $\DY({\C})^H_H$ is the category of
crossed modules. The morphisms are right-module-right-comodule
morphisms over $H$. In a similar way all other combinations of
crossed modules will be defined.
\end{definition}
\abs
In the first part of Figure 1 we exemplarily
represent the right module morphism property of the right
coaction of a Hopf bimodule. The defining identity (\ref{crossed-mod1}) of
a crossed module is represented graphically in the second part of Figure 1.
\begin{figure}
\begin{center}
\unitlength=0.4ex   
\begin{picture}(140,66)(0,5)
\put(25,65){\makebox(0,0){Hopf bimodule identity}}
\put(0,15){\hopfbieins}
\put(35,30){\makebox(0,0){$=$}}
\put(45,20){\hopfbizwei}
\put(65,30){\makebox(0,0){;}}
\put(110,65){\makebox(0,0){Crossed module identity}}
\put(80,15){\hopfbieins}
\put(115,30){\makebox(0,0){$=$}}
\put(125,0){\crossedeins}
\end{picture}
\end{center}
\caption{Some defining relations for Hopf bimodules and crossed modules}
\end{figure}
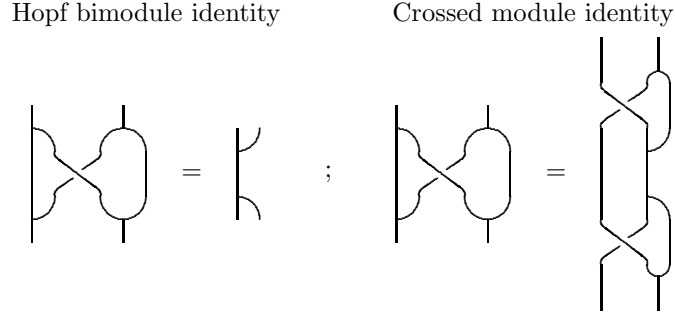
\abs
\begin{example}\label{r-adj-cross}
 A Hopf algebra $H$ is a crossed module $H_{\rm ad}$ over itself through
 the right adjoint action
 \begin{equation}\label{r-adj-act}
 \mu_r^{\rm ad}:=
 \mu_l\circ(\id_H\otimes\mu_r)\circ(\Psi_{X,H}\otimes \id_H)\circ
    (\id_X\otimes (S\otimes\id_H)\circ\Delta)
 \end{equation}
 and the comultiplication $\Delta$ as regular coaction. Dually
 $H$ is a crossed module $H^{\rm ad}$ through the regular action and the
 right coadjoint coaction.
\end{example}
\abs
The following description of the functor ${}_H(\_\,)$ holds for Hopf
modules in general. However we restrict to Hopf bimodules here for
conceptional reasons. We state
the corresponding results of \cite{BD} in invariant form. 
If $H$ is a Hopf algebra in the category $\C$ and
$(X,\mu_r,\mu_l,\nu_r,\nu_l)$ is an $H$-Hopf bimodule then 
there exists an object ${}_HX$ such that
${}_HX\cong \E\underset{H}{\otimes} X$, 
${}_HX\cong \E\underset{H}{\square} X$ and ${}_X\p\circ{}_X\i=\id_{{}_HX}$
where ${}_{\scriptscriptstyle X}\p:\E\otimes X\cong X\to
{}_HX\cong \E\underset{H}{\otimes}X$
and ${}_X\i:\E\underset H\square X \cong {}_HX\to X \cong \E\otimes X$ are the 
corresponding universal morphisms. This means in particular that ${}_HX$ is at
the same time tensor product and cotensor product over $H$ of $\E$ and $X$. 
It is unique up to isomorphism, ${}_{\scriptscriptstyle X}\p'=\xi\circ
{}_{\scriptscriptstyle X}\p$ and ${}_{\scriptscriptstyle X}\i'=
{}_{\scriptscriptstyle X}\i\circ\xi^{-1}$. 
The object of (co-)invariants ${}_HX$ is a crossed module. The assignment
${}_H(\_\,):\hhchh\longrightarrow \DY({\C})^H_H$, which is given through
${}_H(X):= \left({}_HX,\,{}_{\scriptscriptstyle X}\p\circ\mu_r\circ
({}_{\scriptscriptstyle X}\i\otimes\id_H),\,
({}_{\scriptscriptstyle X}\p\otimes\id_H)\circ\nu_r\circ
{}_{\scriptscriptstyle X}\i\right)$
for an object $X$, and through ${}_H(f)= {}_{\scriptscriptstyle Y}\p\circ
f\circ{}_{\scriptscriptstyle X}\i$ for a Hopf bimodule morphism $f:X\to Y$, 
defines a functor. It is a braided monoidal equivalence as will be 
described below in Theorem \ref{yd-hopfbi}. 

For the category of Hopf bimodules over a Hopf algebra $H$ a tensor bifunctor 
can be defined such that $\hhchh$ is a (braided) monoidal category \cite{BD}.
We will give an invariant formulation of this fact in the following 
theorem. As pre-requisites we define two bifunctors
$\odot$ and $\boxdot$ on the category $\hhchh$ over a Hopf algebra $H$.
Two objects $X$ and $Y$ of $\hhchh$ yield the $H$-Hopf bimodule
$X\odot Y$ with diagonal left and right actions
$\mu_{d,l}^{X\odot Y}$ and
$\mu_{d,r}^{X\odot Y}$ according
to (\ref{diag-act}), and with induced left and right coactions
$\nu_{i,l}^{X\odot Y}=\nu_l^X\otimes\id_Y$ and
$\nu_{i,r}^{X\odot Y}=\id_X\otimes\nu_r^Y$.
The Hopf bimodule $X\boxdot Y$ is obtained in the dual
symmetric manner. For Hopf bimodule morphisms $f$ and $g$
we define $f\odot g=f\boxdot g=
f\otimes g$. Then the categories $(\hhchh, \odot,
\alpha)$ and $(\hhchh, \boxdot,\alpha)$
are semi-monoidal, i.\ e.\ they are categories which are almost
monoidal, except that the unit object and the relations involving it are
not required. One verifies easily that the identity functor
$\id_\hhchh:\displaystyle \bigl(\hhchh,\boxdot\bigr)
\longrightarrow \bigl(\hhchh,\odot{}^\op\bigr)$
is semi-monoidal with respect to the natural transformation
$\Theta:\displaystyle \boxdot\overset\bullet\to
\odot{}^\op$ defined through
\begin{equation}\label{Hopf-pseudobraid}
 \Theta_{X,Y}:=(\mu_l^Y\otimes\mu_r^X)\circ
 (\id_H\otimes\Psi_{X,Y}\otimes\id_H)\circ
 (\nu_l^X\otimes\nu_r^Y)\,:\, X\otimes Y \to Y\otimes X\,.
\end{equation}
\abs
\begin{theorem}
\label{Hopf-br}
Let $H$ be a Hopf algebra in $\cal C$.
The category $\hhchh$ of Hopf bimodules over $H$ is monoidal with
the regular Hopf bimodule $H$ as unit object and with the
tensor product $\otimes_H$ uniquely defined (up to monoidal equivalence)
by one of the following equivalent conditions for any pair of $H$-Hopf
bimodules $X$ and $Y$.
\begin{itemize}
\item
The $H$-Hopf bimodule $X\otimes_H Y$ is the tensor product over $H$ of the
underlying modules, and the canonical morphism
$\lambda^H_{X,Y}:X\odot Y\to X\otimes_H Y\cong X\underset H\otimes Y$
is functorial in $\hhchh$, i.\ e.\ 
$\lambda^H:\odot\overset\bullet\to \underset H\otimes$.
\item
The $H$-Hopf bimodule $X\otimes_H Y$ is the cotensor product over $H$ of the
underlying comodules, and the canonical morphism
$\rho^H_{X,Y}:X\underset H\square Y\cong X\otimes_H Y\to
X\boxdot Y$ is functorial in $\hhchh$, i.\ e.\ 
$\rho^H:\underset H\square\overset\bullet\to \boxdot $.
\end{itemize}
Then for given natural morphisms $\lambda^H$ and $\rho^H$ it holds
\begin{equation}
\label{rho-lambda}
\rho^H_{X,Y}\circ\lambda^H_{X,Y}=
(\mu^{X}_{r}\otimes\mu^{Y}_{l})\circ
(\id_X\otimes\Psi_{H\,H}\otimes\id_Y)\circ
(\nu^{X}_{r}\otimes\nu^{Y}_{l})\,.
\end{equation}
The category $\hhchh$ is pre-braided through the pre-braiding
${}^{{}^\scripthhchh}\!\Psi_{X,Y}$ uniquely defined by the condition
\begin{equation}\label{Hopf-braid0}
\rho^H_{Y,X}\circ{}^\scripthhchh\Psi_{X,Y}\circ\lambda^H_{X,Y}=
\Theta_{X,Y}\,.
\end{equation}
${}^H_H{\C}^H_H$ is braided if the antipode of $H$ is an isomorphism.
\end{theorem}

\begin{proof} The existence of a braided monoidal structure
on $\hhchh$ obeying the conditions of Theorem \ref{Hopf-br} has been
deduced explicitely in \cite{BD}with the help of the functor ${}_H(\_)$. 
The universal morphisms for the tensor product and the cotensor
product are
\begin{equation}\label{tsr-prod}
\begin{split}
 &\lambda^H_{{}_0\,N,M} : N\otimes M
 \to N{\underset H\otimes}M
\cong N\otimes{}_H(M)\\
&\lambda^H_{{}_0\,N,M} = (\mu_r^N\otimes\id_{{}_HM})\circ(\id_N\otimes
 (\id_H\otimes{}_{\scriptscriptstyle M}\p)\circ\nu_l^M)
\end{split}
\end{equation}
and
\begin{equation} \label{cotsr-prod}
\begin{split}
&\rho^H_{{}_0\,N,M} : N\otimes{}_H(M)\cong N{\underset H\square} M
  \to N\otimes M\\
&\rho^H_{{}_0\,N,M} = (\id_N\otimes \mu^M_l\circ(\id_H\otimes
 {}_{\scriptscriptstyle M}\i))\circ
 (\nu_r^N\otimes\id_{{}_HM})\,.
\end{split}
\end{equation}
Assume now that there exists another natural transformation $\lambda^H:
\odot\overset\bullet\to\underset H\otimes$ obeying the first condition of
the theorem. Once such a morphism $\lambda^H$ is given, the induced
Hopf bimodule structure on the objects $X\otimes_H Y$
is uniquely determined, since $\lambda^H_{X,Y}$ and therefore
$\lambda^H_{H\boxdot X,Y}\cong\id_H\otimes\lambda^H_{X,Y}$
are Hopf bimodule epimorphisms. Then there are functorial 
$H$-Hopf bimodule isomorphisms $\xi$ between different tensor
products, $\xi_{X,Y}: X\otimes {}_H(Y)\to X\otimes_H Y$, such that
$\xi_{X,Y}\circ\lambda^H_{{}_0\,X,Y}=\lambda^H_{X,Y}$. Then the second
condition of the theorem is immediately verified by defining
$\rho^H_{X,Y}:= \rho^H_{{}_0\,X,Y}\circ \xi_{X,Y}$. In a canonical way
the pre-braided monoidal structure induced by 
$\lambda^H$ can be established. Starting with a natural morphism 
$\rho^H:\underset H\square \overset\bullet\to
\boxdot$ obeying the second condition of the theorem,
analogous results are obtained in the dual manner. 
\end{proof}
\abs
\begin{remark} {\rm If the antipode $S$ of $H$ is an isomorphism, the
inverse ${}^\scripthhchh\Psi^{-1}$ of the braiding ${}^\scripthhchh\Psi$
is given by 
\begin{equation}\label{inv-braid}
{}^\scripthhchh{\Psi^{-1}}_{X,Y} = \lambda^H_{Y,X}\circ
(\mu^Y_r\circ{\Psi^{-1}}_{H,Y}\otimes\id_X)\circ
(S^{-1}\otimes{\Psi^{-1}}_{X,Y})\circ
({\Psi^{-1}}_{X,H}\circ\nu^X_r\otimes\id_Y)\circ
\rho^H_{X,Y}\,.
\end{equation}
In a symmetric way \eqref{inv-braid} can be expressed analogously through
the left (co-)actions $\mu_l^X$ and $\nu^Y_l$ because of the universal
property of the (co-)tensor product. The proof of \eqref{inv-braid} is
based on the following observation. For $f:X\underset H\otimes Y\to U
\underset H\square V$ we put $\tilde f:=\rho^H_{U,V}\circ f\circ
\lambda^H_{X,Y}$. Then $\widetilde{fg}=\tilde f \circ
(\mu_r\otimes\id)(\id\otimes S\otimes\id)(\nu_r\otimes\id)\circ\tilde g$.}
\end{remark}
\abs
We outline the correspondence of Hopf bimodules and crossed modules in
braided categories (with split idempotents).
A full inclusion functor of the category of $H$-right crossed modules into
the category of $H$-Hopf bimodules, $H\!\ltimes\!(\_\,):
{\DY({\C})}^H_H\to {}^H_H{\C}^H_H$, is defined by
$H\ltimes (X) ={(H\otimes X,\mu_{i,l}^{H\otimes X}, \nu_{i,l}^{H\otimes X},
\mu_{d,r}^{H\otimes X}, \nu_{d,r}^{H\otimes X})}$
for any right crossed module $X$ and by $H\ltimes(f)= \id_H\otimes f$ for
any crossed module morphism $f$.
The action $\mu_{i,l}^{H\otimes X}$ is the left action
induced by $H$ and $\mu_{d,r}^{H\otimes X}$
is the diagonal action according to (\ref{diag-act}).
Analogously, in the dual manner the coactions
$\nu_{i,l}^{H\otimes X}$ and $\nu_{d,r}^{H\otimes X}$ are defined.
We formulate the equivalence theorem of crossed modules
and Hopf bimodules in braided categories \cite{BD}.
\abs
\begin{theorem}\label{yd-hopfbi}
Let $H$ be a Hopf algebra with isomorphic antipode in $\C$. Then the
categories $\DY({\C})^H_H$ and ${}^H_H{\C}^H_H$ are
braided monoidal equivalent. The equivalence is given by
${H\ltimes (\_\,): \DY({\C})^H_H\longrightarrow \hhchh}$ and
${{}_H(\_\,): \hhchh\longrightarrow \DY({\C})^H_H}$.
\nl\ \endproof
\end{theorem}
\abs
A tuple $(H,B,\underline\eta,\underline\varepsilon)$ of bialgebras $H$ and
$B$, and of bialgebra morphisms $H\buildrel {\underline\eta}\over
\rightarrow B\buildrel{\underline\varepsilon}\over\rightarrow H$ in
$\C$ is called a bialgebra projection on $H$ if
${\underline\varepsilon}\circ{\underline\eta}=\id_H$ \cite{Rad}.
Analogously Hopf algebra projections are defined.
If $(H,B,\underline\eta,\underline\varepsilon)$ is a bialgebra
projection on $H$ in $\C$, then $\underline B=(B,\mu_r^B,\mu_l^B,
\nu_r^B,\nu_l^B)$ is an $H$-Hopf bimodule through the actions and coactions
\begin{equation}\label{proj-hhopf}
\begin{alignedat}{2}
&\mu_l^B=\m_B\circ(\underline\eta\otimes\id_B)\,, &\quad
&\mu_r^B=\m_B\circ(\id_B\otimes\underline\eta)\,,\\
&\nu_l^B=(\underline\varepsilon\otimes\id_B)\circ\Delta_B\,,&\quad
&\nu_r^B=(\id_B\otimes\underline\varepsilon)\circ\Delta_B\,.
\end{alignedat}
\end{equation}
For the formulation of the main theorem of this chapter the notation
of a relative antipode of a Hopf bimodule $(X,\mu_r,\mu_l,\nu_r,\nu_l)$
is needed. It is the Hopf bimodule morphism $S_{X/H}: X\to X$
defined by 
\begin{equation}\label{rel-ant}
S_{X/H}:=M_X\circ(S\otimes \id_X\otimes S)\circ N_X
\end{equation}
where $M_X:=\mu_l\circ(\id_H\otimes\mu_r)$ and 
$N_X:=(\id_H\otimes\nu_r)\circ\nu_l$.
Every bialgebra projection
$(H,B,\underline\eta_B,\underline\varepsilon_B)$
yields an object $\underline B$ according to (\ref{proj-hhopf}) which
is a bialgebra
$(\underline B,\underline\m_B,\underline\eta_B,\underline\Delta_B,
\underline\varepsilon_B)$ in ${}^H_H{\C}^H_H$. The multiplication
$\underline\m_B$ and the comultiplication $\underline \Delta_B$ are defined
uniquely through
\begin{equation}\label{bialg-trans1}
\underline\m_B\circ(\mu_r^B\otimes\id_B)
=\m_B\circ(\id_B\otimes{}^{{}^\scripthhchh}\lambda_B)
 \quad\hbox{and}\quad
(\nu_r^B\otimes\id_B)\circ\underline\Delta_B
=(\id_B\otimes({}^{{}^\scripthhchh}\lambda_B)^{-1})\circ\Delta_B\,.
\end{equation}
Conversely every bialgebra $\underline B=(B,\underline\m_B,\underline\eta_B,
\underline\Delta_B,\underline\varepsilon_B)$ in ${}^H_H{\C}^H_H$
yields the bialgebra $B=(B,\m_B,\eta_B,\Delta_B,\varepsilon_B)$ in $\C$
whose structure morphisms are given by
\begin{equation}\label{bialg-trans2}
\m_B=\underline\m_B\circ\lambda^H_{B,B}\,,\quad
  \eta_B=\underline\eta_B\circ\eta_H\,,\quad
  \Delta_B=\rho^H_{B,B}\circ\underline\Delta_B\,,\quad
  \varepsilon_B=\varepsilon_H\circ\underline\varepsilon_B\,,
\end{equation}
and furthermore $(H,B,\underline\eta_B,\underline\varepsilon_B)$ is a
bialgebra projection on $H$. There also exists
a correspondence of Hopf algebra structures. The antipode of
$\underline B$ is given by 
${\underline S_B=M_B\circ(\id_H\otimes S_B\otimes\id_H)\circ N_B}$
and for any Hopf algebra $\underline B$ in ${}^H_H{\C}^H_H$ the
antipode of $B$ is given by
${S_B=\underline S_B\circ S_{B/H}=S_{B/H}\circ \underline S_B}$.

Let us denote by $H$-Bialg-$\C$ the category
of bialgebra projections on the Hopf algebra $H$. Its objects are the 
projections $B=\{ B,\m_B ,\eta_B ,\Delta_B ,\varepsilon_B,
\underline\eta_B,\underline\varepsilon_B\}$ and its morphisms 
are bialgebra morphisms
$f: B\to D$ such that $f\circ\underline\eta_B=\underline\eta_D$ and 
$\underline\varepsilon_D\circ f=\underline\varepsilon_B$. 
Then the following central theorem can be formulated \cite{BD}.
\abs
\begin{theorem}\label{bialg-proj}
Let $H$ be a Hopf algebra with isomorphic antipode in $\C$. Then
the category $H$-{\rm Bialg}-\,$\C$ of bialgebra projections on $H$
and the category of $H$-Hopf bimodule bialgebras
{\rm Bialg}-\,$\hhchh$ are isomorphic. The isomorphism
$\displaystyle {\rm Bialg\hbox{-}}\,\hhchh
\mathop{\rightleftarrows}^G_F H\hbox{-}{\rm Bialg}\hbox{-}\,{\C}$
is given by (\ref{bialg-trans1}) and (\ref{bialg-trans2}) on the objects
and through the identity mapping on the morphisms.\endproof
\end{theorem}
\abs
\begin{remark} {\rm Equation \eqref{Hopf-braid0} for the braiding
${}^\hhchh\Psi$ in $\hhchh$ admits the following reformulation.
\begin{equation}
\begin{split}
\label{Hopf-braid2}
&\rho_{X^\prime\otimes_HY,X\otimes_HY^\prime}\circ
(\id_{X^\prime}\otimes_H{}^\hhchh\Psi_{X,Y}\otimes_H\id_{Y^\prime})
\circ\lambda_{X^\prime\otimes_HX,Y\otimes_HY^\prime}\\
&=(\lambda_{X^\prime,Y}\otimes\lambda_{X,Y^\prime})\circ
(\id_{X^\prime}\otimes\Psi_{X,Y}\otimes\id_{Y^\prime})
\circ(\rho_{X^\prime,X}\otimes\rho_{Y,Y^\prime})
\end{split}
\end{equation}
for Hopf bimodules $X,Y,X^\prime,Y^\prime$. Then we can give another proof
of the algebra morphism property of the comultiplication $\Delta_B$ in 
\eqref{bialg-trans2} (compare \cite{BD}). By definition
\eqref{bialg-trans2} it holds
$\Delta_B\circ\m_B=
 \rho_{B,B}\circ\underline\Delta_B\circ\underline\m_B\circ\lambda_{B,B}$.
Using \eqref{Hopf-braid2} we derive
\begin{equation}
\label{cat-equiv}
\begin{split}
&(\m_B\otimes\m_B)\circ(\id_B\otimes\Psi_{B,B}\otimes\id_B)\circ
 (\Delta_B\otimes\Delta_B)\\
&=(\underline\m_B\otimes\underline\m_B)\circ
   (\rho_{B,B}\otimes\rho_{B,B})\circ
   (\id_B\otimes\Psi_{B,B}\otimes\id_B)\circ
   (\lambda_{B,B}\otimes\lambda_{B,B})\circ
   (\underline\Delta_B\otimes\underline\Delta_B)\\
&=(\underline\m_B\otimes\underline\m_B)\circ
   \rho_{B\otimes_HB,B\otimes_HB}\circ
   (\id_B\otimes_H{}^\hhchh\Psi_{B,B}\otimes_H\id_B)\circ
   \lambda_{B\otimes_HB,B\otimes_HB}\circ
   (\underline\Delta_B\otimes\underline\Delta_B)\\
&=\rho_{B,B}\circ(\underline\m_B\otimes_H\underline\m_B)\circ
   (\id_B\otimes_H{}^\hhchh\Psi_{B,B}\otimes_H\id_B)\circ
   (\underline\Delta_B\otimes_H\underline\Delta_B)\circ\lambda_{B,B}\\
&=\rho_{B,B}\circ\underline\Delta_B\circ\underline\m_B\circ\lambda_{B,B}
 =\Delta_B\circ\m_B\,.
\end{split}
\end{equation}
Conversely one deduces from \eqref{cat-equiv} the algebra morphism property
of $\underline\Delta_B$ if it holds for $\Delta_B$,
since $\rho_{B,B}$ is monomorphism and $\lambda_{B,B}$ is epimorphism.}
\end{remark}
\abs
\abs
\section{Factor Algebras and Ideals}
\label{ideal-section}
In Chapter \ref{ideal-section} we investigate factor (co-)algebras,
bi- and Hopf algebras in abelian categories which arise from the canonical
decomposition of (co-)algebra, bi- and Hopf algebra morphisms
into monomorphism and epimorphism. Under more restrictive conditions to
the category we will also consider (co-)ideals,
bi- and Hopf ideals. The obtained results provide appropriate
tools for the study of the braided exterior tensor algebras which will be
constructed in Chapter \ref{grade-alg-com}. On the other hand these
objects also play an important r\^ole in the examination of higher order
bicovariant differential calculi in Chapter \ref{sec-diff-calc}. 

We suppose henceforth that the categories are (braided) monoidal abelian
with bi-additive tensor product $\otimes$.
Before we will fix additional conditions of the categories in the 
Definitions \ref{tensor-factor} and \ref{tensor-exact} we provide an
appropriate definition of submodules generated by a given morphism. 
\abs
\begin{definition}\label{29-D30}
Let $A$ be an algebra and $(Y,\mu_l)$ be an $A$-left module in the monoidal
abelian category $\C$. For
any morphism $f:X\to Y$ in $\C$ we define
$A\langle f\rangle := \Im\big(\mu_l\circ(\id_A\otimes f)\big)=
\Ker\coker\big(\mu_l\circ(\id_A\otimes f)\big)$.
In an analogous manner
$\langle f\rangle A$ is defined for an $A$-right module. 
If $(Y,\mu_l,\mu_r)$ is an $A$-bimodule then 
$A\langle f\rangle A :=\Im\big(\mu_r\circ
 (\mu_l\circ(\id_A\otimes f)\otimes\id_A)\big)$.
\end{definition}
\abs
\begin{lemma}\label{1-1ZD30}
Suppose that for any pair of epimorphisms $g$ and $h$ in $\C$
the tensor product $g\otimes h$ is again epimorphic. 
Then $B\langle f\rangle$, $\langle f\rangle B$ and
$B\langle f\rangle B$ are
$B$-left submodules, $B$-right submodules and $B$-sub-bimodules
of $Y$ respectively.
\end{lemma}

\begin{proof}
We consider the $B$-left module $Y$ and use the abbreviation
$\phi:=\mu_l\circ(id_B\otimes f)$. Then we obtain
$\coker\phi\circ\mu_l\circ(\id_B\otimes\phi)=0$.
By assumption on the category $\C$ there exists a unique morphism
$\mu_l':B\otimes\Im\phi\to\Im\phi$ such that
$\im\phi\circ\mu_l'=\mu_l\circ(\id_B\otimes\im\phi)$. It follows
\begin{equation}
\begin{gathered}
\im\phi\circ\mu_l'\circ(\id_B\otimes\mu_l')
=\mu_l\circ(\m_B\otimes\id_Y)\circ(\id_{B\otimes B}\otimes\im\phi)
=\im\phi\circ\mu_l'\circ(\m_B\otimes\id_{\Im\phi})\,,\\
\im\phi\circ\mu_l'\circ(\eta_B\otimes\id_{\Im\phi})=\im\phi\,.
\end{gathered}
\end{equation}
This proves that $\big(B\langle f\rangle=\Im\phi,\mu_l'\big)$ is a $B$-left
submodule of $Y$. Analogously the other cases of the proposition
can be verified.\end{proof}
\abs
The conditions on the category $\C$ used in Lemma \ref{1-1ZD30}
will be needed frequently in the following. We will therefore collect
these properties in the next definition.
\abs
\begin{definition}\label{tensor-factor}
The (braided) monoidal abelian category $\C$ is said to fulfill
the $\otimes$-epimorphism property if the tensor product
$f\otimes g$ of any two epimorphisms $f$ and $g$ in $\C$ is
again epimorphic. The $\otimes$-monomorphism property is defined in
the dual way. The category $\C$ is called $\otimes$-factor
(braided) abelian if both the $\otimes$-epimorphism and the
$\otimes$-monomorphism property hold in $\C$.
\end{definition}
\abs
For $\otimes$-factor (braided) abelian categories 
we provide some facts on factor (co-)algebras, factor bi- and Hopf
algebras emerging out of the corresponding morphisms. This lemma will
be used in the sequel, especially in Chapters
\ref{grade-alg-com} and \ref{sec-diff-calc}.

\begin{lemma}\label{8-4Z31}
Let $f:A\to C$ be a (co-)algebra morphism in the $\otimes$-factor (braided)
abelian category $\C$. Then the decomposition of $f$ into its
image and coimage $A@>{\coim f}>>B@>{\im f}>>C$ admits
a unique (co-)algebra structure on $B$ turning $\im f$ and $\coim f$ into
(co-)algebra morphisms.

The analogous results hold for bi- and Hopf algebras
in $\otimes$-factor braided abelian categories.
\end{lemma}

\begin{proof}
The proof works similar as in classical ring
theory. We will therefore only sketch it for the algebra case.
Since $\C$ is a $\otimes$-factor abelian category there exists a
unique morphism $\m_B$ making the following diagram commutative.
\begin{equation}
\bfig
\putsquare<1`1`1`1;1200`500>(0,-250)[A\otimes A`B\otimes B`A`B;
 \coim f\otimes\coim f`\m_A`\m_B`\coim f]
\putsquare<0`0`1`1;1150`500>(1250,-250)[`C\otimes C``C;
 ``\m_C`\im f]
\putmorphism(1335,250)(1,0)[``\im f\otimes\im f]{930}1a 
\efig
\end{equation}
Then the definition $\eta_B:=(\coim f)\circ\eta_A$ renders
$(B,\m_B,\eta_B)$ a unital associative algebra.
\end{proof}
\abs
If not otherwise stated we will work in the subsequent chapters with
$\otimes$-factor (braided) abelian categories.
However, when we consider (co-)ideals, bi- and Hopf ideals, we will use
categories which have a more restricted structure. This will be discussed
in the next section.

\subsection*{Ideals}

\begin{definition}\label{tensor-exact}
We call the (braided) abelian category $\C$
$\otimes$-right-exact$/$$\otimes$-left-exact (braided) abelian if
the following conditions are fulfilled.
For any two epimorphisms$/$monomorphisms $f_i:X_i\to Y_i$,
$i\in\{1,2\}$ in $\C$
we suppose that the push-out$/$pull-back of the pair of morphisms
$\{f_1\otimes \id\,,\,\id\otimes f_2\}$ is given by
$\{\id\otimes f_2\,,\,f_1\otimes \id\}$.
We say that the category $\C$ is $\otimes$-exact (braided) abelian if
it is $\otimes$-left-exact abelian and $\otimes$-right-exact abelian.
\end{definition}
\abs
Both the dual categories of a $\otimes$-factor abelian and a $\otimes$-exact
abelian category are again $\otimes$-factor abelian and
$\otimes$-exact abelian respectively.
In some cases we formulate statements which are not invariant under
categorical duality. However, in this case the reader may immediately
verify the dual assertion.

The $\otimes$-right-exact property means that for any pair of morphisms
$l:X_1\otimes Y_2\to Z$ and $h:Y_1\otimes X_2\to Z$, which obey the
identity $l\circ(\id_{X_1}\otimes f_2)=h\circ (f_1\otimes\id_{X_2})$,
there exists a unique morphism
$k:Y_1\otimes Y_2\to Z$ such that $l=k\circ(f_1\otimes\id_{Y_2})$
and $h=k\circ(\id_{Y_1}\otimes f_2)$. As a consequence 
$f\otimes g$ is epimorphic whenever $f$ and $g$ are epimorphic.
Similarly in $\otimes$-left-exact abelian categories
$k\otimes l$ is monomorphic when $k$ and $l$ are monomorphic. Hence
$\otimes$-exact abelian categories are $\otimes$-factor abelian.
\abs
Suppose that $f_k:Y_k\to Z,\, k\in\{1,2\}$ are morphisms in an abelian
category and the push-out of 
$\{\coker f_1,\coker f_2\}$ is given by $\{g_1,g_2\}$.
Then 
\begin{equation}\label{factor-latt-eq}
g_1\circ f_1=g_2\circ f_2=
\coker\big( Y_1\oplus Y_2 @>{(f_1,f_2)}>>Z\big)\,.
\end{equation}
Of course, an analogous dual symmetric result holds.
An immediate implication of this fact is the next corollary.
\begin{corollary}\label{15-31}
Let $\C$ be a $\otimes$-right-exact abelian category
and $X_1$, $X_2$ be objects in $\C$ such that
$X_1\otimes\id_{\C}$ and $\id_{\C}\otimes X_2$ are right exact.
Then for any pair of morphisms $h_k:Y_k\to X_k$, $k\in\{1,2\}$ it holds
\begin{equation}
\label{15-31-rexact}
\coker(h_1)\otimes\coker(h_2)=\coker
\big((h_1\otimes\id_{X_2}),(\id_{X_1}\otimes h_2)\big)
\end{equation}
\end{corollary}

\begin{proof}
Applying equation (\ref{factor-latt-eq}) to 
$f_1:=h_1\otimes\id_{X_2}$ and $f_2:=\id_{X_1}\otimes h_2$, and
keeping in mind that $\C$ is $\otimes$-right-exact abelian and
that $X_1\otimes\id_{\C}$, $\id_{\C}\otimes X_2$ are right exact
yields the result.
\end{proof}
\abs
\begin{remark}
{\rm Suppose that $X\otimes\id_{\C}$ and $\id_{\C}\otimes X$ are
right exact for every object $X$ in the abelian monoidal category $\C$ --
this holds if $\C$ is closed, for instance.
Then the following statements are equivalent.
\begin{enumerate}
\item
$\C$ is $\otimes$-right-exact abelian.
\item
For any pair of epimorphisms $f_k:X_k\to Y_k$, $k\in\{1,2\}$
the tensor product $f_1\otimes f_2$ is epimorphic. If
$h:X_1\otimes X_2\to Z$ factorizes over $f_1\otimes\id_{X_2}$ and
$\id_{X_1}\otimes f_2$, then $h$ factorizes over $f_1\otimes f_2$.
\item
For any pair of morphisms $h_1$, $h_2$ the equation
(\ref{15-31-rexact}) holds.
\end{enumerate}

\begin{proof}
Obviously (2) implies (1) which on the other hand yields (3) because
of Corollary \ref{15-31}. It remains to show that (3) implies (1). 
For the epimorphisms $f_k,\,k\in\{1,2\}$ the equation
(\ref{15-31-rexact}) yields the relation
$f_1\otimes f_2=
 \coker\left((\ker f_1\otimes\id_{X_2},\id_{X_1}\otimes\ker f_2)\right)$.
If the morphism $h$ factorizes over $f_1\otimes\id_{X_2}$ and
$\id_{X_1}\otimes f_2$ then $h\circ(\ker f_1\otimes\id_{X_2})=0$
and $h\circ(\id_{X_1}\otimes\ker f_2)=0$. Hence $h$ factorizes over
$f_1\otimes f_2$.\end{proof}}
\end{remark}
\abs
We suppose in the remainder of this section that
the categories under consideration are $\otimes$-exact (braided) abelian,
although many results can be derived under weaker conditions.

In the following we imitate the well known definition of ideals, co-, bi-
and Hopf ideals of classical ring theory.
\abs
\begin{definition}\label{16-17-32}
Let $(A,\m,\eta)$ be an algebra in $\C$.
If there exists a subobject of $A$ (represented by)
$I\overset{\i}\hookrightarrow A$ such that
\begin{equation}\label{ideal}
\im{\left(\m\circ
        \big((\i\otimes\id_A)\,,\,(\id_A\otimes\i)\big) \right)}
\subset \i
\end{equation}
then $(I,\i)$ is called ideal in $A$.
If $(C,\Delta,\varepsilon)$ is a coalgebra in $\C$ and
$J\overset{\j}\hookrightarrow C$ is a subobject of $C$ with
\begin{equation}\label{coideal}
 \im{(\Delta\circ\j)}\subset
 \im{\big((\j\otimes\id_A)\,,\,(\id_A\otimes\j)\big)},\quad\text{and}\quad
 \varepsilon\circ\j=0
\end{equation}
then $(J,\j)$ is called coideal in $C$. 
For a bialgebra $B$ in $\C$ a subobject $(I,\i)$ of $B$ is called
bi-ideal if both $(I,\i)$ is ideal and coideal in $B$;
if $B$ is a Hopf algebra then $(I,\i)$ is called Hopf ideal if in addition
$\coker(\i)\circ S\circ\i=0$.
\end{definition}
\abs
\begin{remark}{\rm 
The ideals of an algebra $A$ in $\C$ are in one-to-one correspondence
with the (equivalence classes of) subobjects of the regular bimodule $A$ in
${}_A{\C}_A$.

Of course the definition of an ideal and of a coideal are not mutually
dual in the categorical sense. For instance in the category of vector
spaces the dual of a coideal is the quotient space of an algebra and
a subalgebra together with the canonical surjective mapping.}
\end{remark}
\abs
The following propositions generalize classical results
of ring theory to $\otimes$-exact braided abelian categories. Proposition
\ref{17-32} describes the correspondence of ideals and factor algebras.
In Proposition \ref{2-3ZD30} ideals generated by morphisms are described.
\abs
\begin{proposition}\label{17-32}
Let $\C$ be a $\otimes$-right-exact braided abelian category. Suppose
that $A$ is an algebra in $\C$ and $A\otimes\id_{\C}$ is right exact.
Then every ideal of $A$ represented by $I\overset\i\hookrightarrow A$
induces a unique algebra structure on the morphism $\coker\i:A\to
\bar A:=\Coker\i$. Conversely every epimorphic algebra morphism
$p:A\to\bar A$ is isomorphic to the cokernel of a morphism $\i$
representing an ideal of $A$.

Analogous relations hold in the case of coalgebras, bi- and Hopf algebras
and their corresponding co-, bi- and Hopf ideals respectively.
\end{proposition}

\begin{proof}
If $(A,\m,\eta)$ is an algebra and $\i:I\hookrightarrow A$ is an ideal of
$A$, then (\ref{ideal}) holds, and we can apply Corollary \ref{15-31} to
obtain a unique morphism $\bar\m:\Coker\i\otimes \Coker\i\to \Coker\i$ such
that $\coker\i\circ\m=\bar\m\circ(\coker\i\otimes\coker\i)$.
Obviously $(\Coker\i,\bar\m,\bar\eta:=\coker\i\circ\eta)$ is the unique
algebra making $\coker\i$ an algebra morphism.
Conversely, let $\p:A\to \bar A$ be an epimorphic algebra morphism.
Then for $\i:=\ker\p$ the ideal property (\ref{ideal}) is verified directly,
$\p\circ\m_A\circ\big((\i\otimes\id_A)\,,\,(\id_A\otimes\i)\big)=\m_P\circ
(\p\otimes\p)\circ\big((\i\otimes\id_A)\,,\,(\id_A\otimes\i)\big)=0$.

Let $(C,\Delta,\varepsilon)$ be a coalgebra and $\j:J\to C$ be a
coideal of $C$. Then because of (\ref{coideal}) there exists a
unique morphism $\bar\Delta:\Coker\j\to \Coker\j\otimes\Coker\j$
and a unique morphism $\bar\varepsilon:\Coker\j\to\E$
such that $\bar\Delta\circ \coker\j= \Delta\circ(\coker\j\otimes\coker\j)$
and $\bar\varepsilon\circ\coker\j=\varepsilon$. Similarly as in the algebra
case it follows that $(\Coker\j,\bar\Delta,\bar\varepsilon)$ is the unique
coalgebra so that $\coker\i$ is a coalgebra morphism.
Conversely an epimorphic coalgebra morphism $\p:C\to \bar C$
leads to $(\p\otimes\p)\circ\Delta\circ\ker\p=0$ and
$\varepsilon\circ\ker\p=\bar\varepsilon\circ\p\circ\ker\p=0$.
Therefore (\ref{coideal}) holds for $\j=\ker\p$.

Bi-ideals of a bialgebra $B$ induce factor objects which are at the same time
algebras and coalgebras. The bialgebra properties of $B$ then yield
the compatibility of both structures leading to the factor bialgebra.
The existence of the antipode and the Hopf algebra relations of the factor
bialgebra arising from a Hopf ideal of a Hopf algebra $H$ can be proven
similarly.\end{proof}
\abs
\begin{remark}
{\rm 1. In the coalgebra case of the previous proposition
the $\otimes$-right-exactness of $\C$ and the right
exactness of $C\otimes\id_{\C}$ need not be required.
\nl
2. Under the assumptions of Proposition \ref{17-32}
it follows from Proposition \ref{8-4Z31} that the kernel of any
algebra morphism is an ideal whoose structure is compatible with the
given algebra structure. Analogous statements are obtained for co-, bi-
and Hopf algebras.}
\end{remark}
\abs
\begin{proposition}
\label{2-3ZD30} Let $A$ be an algebra and $g:X\to A$ be any morphism in the
category $\C$. Then $(g):=A\langle g\rangle A$ is an ideal in $A$.
Suppose that $B$ is a bialgebra (Hopf algebra) with right exact functor
$B\otimes\id_\C$, and $f:X\rightarrow B$ is a morphism in $\C$
satisfying the condition
\begin{equation}
\label{Delta-f}
(\coker f\otimes\coker f)\circ\Delta\circ f=0\,,\quad 
\varepsilon_B\circ f=0\quad (and\quad \coker f\circ S\circ f=0)\,. 
\end{equation}
Then the ideal $(f)$ generated by $f$ is a bi-ideal (Hopf ideal) in $B$.
\end{proposition}

\begin{proof}
Given an algebra $A$, we consider the morphism
$\tilde\phi_A:=\m_A\circ(\id_A\otimes\m_A)\circ
(\id_A\otimes g\otimes\id_A)$. Analogously the morphism $\tilde\phi_B$ is
defined for a bialgebra $B$. Because of Lemma \ref{1-1ZD30} we know that
\begin{equation}
\coker\tilde\phi\circ\m\circ\big((\im\tilde\phi\otimes\id),
 (\id\otimes\im\tilde\phi)\big)
=\coker\tilde\phi\circ(\im\tilde\phi\circ\mu_r',\im\tilde\phi\circ\mu_l')
=0\,.
\end{equation}
Thus $(g)$ and $(f)$ are ideals according to Definition \ref{16-17-32}.
Since $B\otimes\id_\C$ is supposed to be right exact we obtain the following
identities involving a unique morphism $\xi$.
\begin{equation}\label{coideal-proof1}
 \coker\tilde\phi_B\circ\m_B\circ(\id_B\otimes\m_B)=
 \xi\circ(\id_B\otimes\coker f\otimes\id_B)
\end{equation}
where use has been made of the fact that $\tilde\phi_B$ factorizes over
$(\id_B\otimes f\otimes\id_B)$. Keeping in mind that $\Delta_B$ is an
algebra morphism we arrive at
\begin{equation}\label{coideal-proof2}
 (\coker\tilde\phi_B\otimes \coker\tilde\phi_B)\circ\Delta_B\circ
 \im\tilde\phi_B=0\,.
\end{equation}
For the derivation of (\ref{coideal-proof2}) we used (\ref{coideal-proof1})
and the assumption (\ref{Delta-f}) for the morphism $f$. The explicit
calculation yielding (\ref{coideal-proof2}) is presented graphically in
Figure 2.
\begin{figure}
\begin{center}
\unitlength=0.4ex   
\begin{picture}(240,100)

\put(15,0){\line(0,1){10}}
\put(15,15){\morphoval{\coker\tilde\phi}}
\put(41,0){\line(0,1){10}}
\put(41,15){\morphoval{\coker\tilde\phi}}
\put(15,20){\curve(0,0, 0,5, 0,10, 4,15, 8,20, 8,25, 8,30)}
\put(28,60){\comult}
\put(41,20){\curve(0,0, 0,5, 0,10, -4,15, -8,20, -8,25, -8,30)}
\put(28,70){\mult}
\put(23,80){\mult}
\put(28,80){\line(0,1){6}}
\put(28,90){\morphcirc{f}}
\put(28,94){\line(0,1){6}}
\put(18,80){\line(0,1){20}}
\put(33,70){\line(0,1){30}}

\put(61,50){\makebox(0,0){$=$}}

\put(81,0){\line(0,1){10}}
\put(81,15){\morphoval{\coker\tilde\phi}}
\put(107,0){\line(0,1){10}}
\put(107,15){\morphoval{\coker\tilde\phi}}
\put(81,20){\line(0,1){10}}
\put(107,20){\line(0,1){10}}
\put(81,40){\mult}
\put(76,40){\line(0,1){5}}
\put(107,40){\mult}
\put(112,40){\line(0,1){5}}
\put(76,55){\mult}
\put(71,55){\line(0,1){20}}
\put(112,55){\mult}
\put(117,55){\line(0,1){20}}
\put(81,75){\curve(0,0, 0,-3, 0,-5, 4,-10, 13,-20, 21,-30, 21,-32, 21,-35)}
\put(99,75){\curve(0,0, 0,-3, 0,-5, 4,-10, 8,-15, 8,-17, 8,-20)}
\put(89,75){\curve(0,0, 0,-1, 0,-3, -3,-9)}
\put(107,75){\curve(0,0, 0,-1, 0,-3, -3,-9)}
\put(102,64){\curve(0,0, -7,-8)}
\put(81,55){\curve(0,0, 0,1, 0,2, 1,5, 2,7, 3,9)}
\put(86,40){\curve(0,0, 0,1, 0,3, 6,12, 7,13)}
\put(76,85){\comult}
\put(76,85){\line(0,1){15}}
\put(112,85){\comult}
\put(112,85){\line(0,1){15}}
\put(94,85){\comult}
\put(94,89){\morphcirc{f}}
\put(94,93){\line(0,1){7}}

\put(127,50){\makebox(0,0){$=$}}

\put(147,0){\line(0,1){10}}
\put(147,15){\morphoval{\xi}}
\put(225,0){\line(0,1){10}}
\put(225,15){\morphoval{\xi}}
\put(142,20){\line(0,1){70}}
\put(230,20){\line(0,1){70}}
\put(147,20){\curve(0,0, 0,3, 6,13, 16,22, 18,24)}
\put(152,20){\curve(0,0, 0,1, 5,5, 20,11, 33,15)}
\put(220,20){\curve(0,0, 0,1, -20,10, -34,16, -54,25, -68,35, -68,45, -68,70)}
\put(225,20){\curve(0,0, 0,3, 0,5, -6,13, -13,18, -19,23, -26,30, -26,32, -26,35)}
\put(167,46){\curve(0,0, 1,1, 5,5, 6,8, 6,9)}
\put(188,36){\curve(0,0, 16,7)}
\put(220,90){\curve(0,0, 0,-30, -4,-35, -13,-45)}
\put(147,100){\comult}
\put(186,100){\line(0,-1){7}}
\put(186,89){\morphcirc{f}}
\put(186,85){\comult}
\put(225,100){\comult}
\put(181,75){\curve(0,0, 0,-2, 0,-3, -4,-5, -8,-7, -8,-8, -8,-10)}
\put(191,75){\curve(0,0, 0,-2, 0,-3, 4,-5, 8,-7, 8,-8, 8,-10)}
\put(173,60){\morphoval{\coker f}}
\put(199,60){\morphoval{\coker f}}
\end{picture}
\end{center}
\caption{Proof of the coideal property of $\im\tilde\phi$.}
\end{figure}
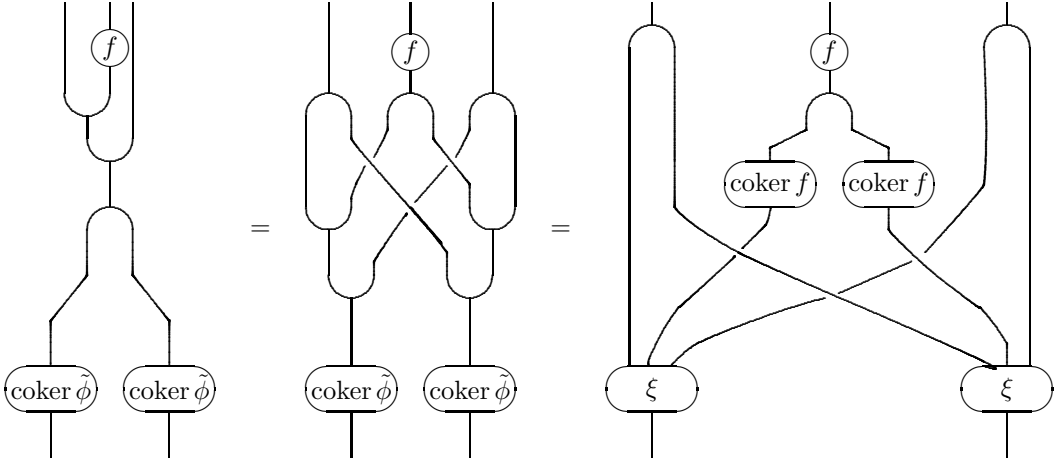  
To verify the counit relation one observes that equivalent statements
${\varepsilon_B\circ\im\tilde\phi_B=0}\Leftrightarrow
{\varepsilon_B\circ\tilde\phi_B=0}\Leftrightarrow
{\varepsilon_B\otimes\varepsilon_B\circ f\otimes\varepsilon_B=0}$
are fulfilled under the assumption (\ref{Delta-f}). Using similar
techniques as in the proof of (\ref{coideal-proof2}) the antipode relation
${\coker\tilde\phi_B\circ S_B\circ\im\tilde\phi_B=0}$ can be checked
if $B$ is a Hopf algebra. Therefore all the conditions of
Definition \ref{16-17-32} hold and the proposition is proved.\end{proof}
\abs
\abs
\section{Graded Categories, Complexes and Exterior Hopf Algebras}
\label{grade-alg-com}
In this chapter we are considering graded categories and categories of
complexes over braided abelian base categories $\C$. We are investigating
(co-)algebras, bi- and Hopf algebras in these categories. In particular
we study graded tensor algebras over an arbitrary object $X$ in $\C$ and
derive Hopf algebra structures on them using braided combinatorics. The
corresponding braided antisymmetric tensor Hopf algebras and the braided
exterior Hopf algebras of forms are deduced, and their universal properties
are shown. In Chapter \ref{sec-diff-calc} these objects are shown to
have a certain unique differential structure if $X$ is a bicovariant
differential calculus.
\nl
For the investigation of graded and complex categories we go some steps
beyond and study functor categories ${\C}^{\II}$ over coalgebra
categories $\II$. The notion of a coalgebra, bi- or Hopf algebra
category is first given in a very general sense in \cite{CF}.
The properties of (quasitriangular) coalgebra categories allow us
to define a (braided) monoidal structure on ${\C}^{\II}$
if $\C$ is (braided) monoidal. 
The reader who wants to avoid reading this rather formal investigation
may skip the following section and should continue at the definition of
graded categories and complexes over $\NN_0$ or $\{0,1\}$ which are special
examples of the general situation.
\abs
\subsection*{Matrix Categories, Coalgebra Categories and Functor Categories}

Let $\Bbbk$ be a commutative, unital ring.
In this section we suppose that the categories under consideration
are $\Bbbk$-linear categories and the functors are $\Bbbk$-linear.
In particular the ring $\Bbbk$ itself is a $\Bbbk$-linear category with a
single object, $\Ob({\Bbbk})=\{\Bbbk\}$ and $\End_{\Bbbk}({\Bbbk})
=\Bbbk$. For a given
$\Bbbk$-linear category $\C$ we consider the matrix category
$\Mat{\C}$. It consists of objects $\vec X=(X_1,\ldots,X_n)$, 
$n\in\NN_0$ which are finite ordered sets of objects $X_i\in\Ob(\C)$,
$\i\in\{1,\ldots,n\}$. For $n=0$ we formally denote by $\emptyset$
the ``empty" zero object in $\Mat{\C}$. The morphisms
$\hat g:\vec X\to \vec Y$ are given by $\hat g=\{g_{i\,j}\}_{(i,j)}$, where
$g_{i\,j}:X_j\to Y_i$ is a morphism in $\Hom_{\,\C}(X_j,Y_i)$. The
composition of morphisms is matrix-like for every single
component.
Then obviously $\Mat{\C}$ is also a $\Bbbk$-linear category.

A $\Bbbk$-linear functor $F:{\C}\to \Mat{\cal D}$ extends
to a $\Bbbk$-linear functor $\hat F:\Mat{\C}\to\Mat{\cal D}$
through $\hat F(\vec X)=(F(X_1),\ldots,F(X_n))$ and $\hat F(\hat g)
=\{F(g_{i\,j})\}$. Then the canonical embedding
${\C}\hookrightarrow\Mat{\C}$ extends to the identity functor
on $\Mat{\C}$ and $(\hat F\circ G)^{\wedge}=\hat F\circ\hat G$.
These facts permit the construction of the (large)
category ${\Bbbk}\text{-}\Mcat$. Its objects are (small)
$\Bbbk$-linear categories, and morphisms
$F\in \Hom_{{\Bbbk}\text{-}\Mcat}({\C},{\cal D})$ are $\Bbbk$-linear
functors $F:{\C}\to \Mat{\cal D}$.

For $\Bbbk\text{-}\Mcat$ a bi-functor $\uotimes:\Bbbk\text{-}\Mcat
\times \Bbbk\text{-}\Mcat \longrightarrow \Bbbk\text{-}\Mcat$ can be defined
as follows. The object ${\C}\uotimes{\cal D}$ is the
$\Bbbk$-linear category with objects $(X,Y)\in \Ob(\C)\times
\Ob({\cal D})$, and with morphisms
$\Hom_{\,\C\uotimes{\cal D}}((X_1,Y_1),(X_2,Y_2))
=\Hom{\,\C}(X_1,X_2)\otimes_\Bbbk \Hom{\,{\cal D}}(Y_1,Y_2)$.
For $\Bbbk$-linear functors
$F:{\C}_1\longrightarrow{\C}_2$ and
$G:{\cal D}_1\longrightarrow{\cal D}_2$ the $\Bbbk$-linear functor
$F\uotimes G:{\C}_1\uotimes{\cal D}_1\longrightarrow
{\C}_2\uotimes{\cal D}_2$ is given by
$F\uotimes G(X,Y)=(F(X),G(Y))$ in lexicographic ordering on the
objects and by $F\uotimes G(f\otimes_{\Bbbk} g)= \{F(f)_{i\,j}
\otimes_{\Bbbk} G(g)_{k\,l}\}$ on the morphisms.
In addition functorial morphisms $\varphi_1:F_1\overset{\bullet}{\to}G_1$
and $\varphi_2:F_2\overset{\bullet}{\to}G_2$ compose to a functorial
morphism $\varphi_1\uotimes\varphi_2: F_1\uotimes G_1
\overset{\bullet}{\to}F_2\uotimes G_2$ which is defined through
$(\varphi_1\uotimes\varphi_2)_{(X,Y)}=\{(\varphi_1)_{X,\,k\,l}
\otimes_{\Bbbk}(\varphi_2)_{Y,\,m\,n}\}$.

Because of the symmetry of the tensor product in $\Bbbk$-mod
it is possible to define an involutive functorial isomorphism
$\underline{\tau}:\uotimes\overset{\bullet}{\cong}
\uotimes^{\rm op}$ which is the analogue of the tensor
transposition in $\Bbbk$-mod. It is given by
$\underline{\tau}_{\,{\C}\,{\cal D}}(X,Y)=(Y,X)$ on the objects and by
$\underline{\tau}_{\,{\C}\,{\cal D}}(f\otimes_{\Bbbk}g)=
g\otimes_{\Bbbk}f$ on the (generating) morphisms. We collect these
results in the following lemma.
\abs
\begin{lemma}\label{mcat-monoidal}
The category $(\Bbbk\text{-}\Mcat,\uotimes,{\Bbbk},
\underline{\tau})$ is symmetric monoidal. The functorial
isomorphisms which are ruling the associativity and the unit object
property are naturally induced by the monoidal structure of
$\Bbbk$-{\rm mod}.\endproof
\end{lemma}
\abs
Lemma \ref{mcat-monoidal} allows us to use the notion of a {\it coalgebra
category} \cite{CF} in $\Bbbk\text{-}\Mcat$. It is a category ${\II}$ in
$\Bbbk\text{-}\Mcat$ with $\Bbbk$-linear functors
$\underline{\varepsilon}_{\II}:{\II}\longrightarrow \Mat{\Bbbk}$
and $\underline{\Delta}_{\II}:{\II}\longrightarrow
\Mat{{\II}\uotimes{\II}}$ such that
$({\II},\underline{\Delta}_{\II},\underline{\varepsilon}_{\II})$
is a coalgebra in $\Bbbk\text{-}\Mcat$.
  
Since $\Bbbk\text{-}\Mcat$ is a 2-category where the functorial morphisms
are the 2-morphisms, it is in particular possible to define the
composition of functors and functorial morphisms \cite{MP,Mac}. Let
$F,G:{\C}\longrightarrow\Mat{\cal D}$,
$K:{\cal D}\longrightarrow \Mat{\cal E}$ and
$L:{\cal F}\longrightarrow\Mat{\C}$ be functors in $\Bbbk\text{-}\Mcat$.
Given any functorial morphism $\varphi:F\overset{\bullet}{\to}G$, then the
compositions $K\diamond\varphi:K\circ F\overset{\bullet}{\to} K\circ G$,
with $(K\diamond\varphi)_X=K(\varphi_X)$ and $\varphi\diamond L:
F\circ L\overset{\bullet}{\to} G\circ L$, with $(\varphi\diamond L)_Y=
\varphi_{L(Y)}$ are functorial morphisms. This will be used in the following
definition.
\abs
\begin{definition}\label{qtcc}
A coalgebra category $({\II},\underline{\Delta}_{\II},
\underline{\varepsilon}_{\II})$ is called quasitriangular if
there exists a functorial isomorphism 
$\underline{R}_{\II}:\underline{\Delta}_{\II}\overset{\bullet}{\to}
\underline{\Delta}^\op_{\II}:=\underline{\tau}_{{\II},{\II}}\circ
\underline{\Delta}_{\II}$ satisfying the following dual hexagon
identities.
\begin{equation}
\label{coHex}
\begin{split}
(\id_{\II}\uotimes\underline{\Delta}_{\II})\diamond
\underline{R}_{\II} &=
\big(\big((\underline{\tau}_{{\II},{\II}}\uotimes\id_{\II})
\diamond (\id_{\id_{\II}}\uotimes\underline{R}_{\II})\big)\circ
(\underline{R}_{\II}\uotimes\id_{\id_{\II}})\big)\diamond
\underline{\Delta}_{\II}\,,\\
(\underline{\Delta}_{\II}\uotimes\id_{\II})\diamond
\underline{R}_{\II} &=
\big(\big((\id_{\II}\uotimes\underline{\tau}_{{\II},{\II}})
\diamond (\underline{R}_{\II}\uotimes\id_{\id_{\II}})\big)\circ
(\id_{\id_{\II}}\uotimes\underline{R}_{\II})\big)\diamond
\underline{\Delta}_{\II}\,.
\end{split}
\end{equation}
In (\ref{coHex}) the composition ``$\,\circ$" denotes the usual
(pointwise) 2-composition of functorial morphisms.
\end{definition}
\abs
\begin{remark}{\rm
One can write the identities (\ref{coHex}) in a more
common form.
\begin{equation*}
(\id\uotimes\underline{\Delta})\diamond\underline{R}=
(\underline{R}_{13}*\underline{R}_{12})
\qquad\text{and}\qquad
(\underline{\Delta}\uotimes\id)\diamond\underline{R}=
(\underline{R}_{13}*\underline{R}_{23})
\end{equation*}
where obvious notations have been used.}
\end{remark}
\abs
Before we are going to formulate the proposition on braided monoidal
structures of functor categories, we have to provide another two functors.
For any $\Bbbk$-linear abelian monoidal category $\C$ (with
$\Bbbk$-bilinear functor $\otimes$ henceforth) we define the
$\Bbbk$-linear functor $\bullet_{\C}:\Mat{{\C}\uotimes{\C}}
\longrightarrow{\C}$
\begin{equation}\label{funct-CC-C}
\begin{split}
 \bullet_{\C}:\ \left\{
 \begin{matrix}
 (X_i,Y_i)_{i=1}^n &\longmapsto &\bigoplus_{i=1}^n X_i\otimes Y_i\\
 \Big\{\sum_{k(j,i)} f_{j\,i}^k\otimes_{\Bbbk} g_{j\,i}^k\Big\}_{j,i}
 &\longmapsto
 &\Big(\sum_{k(j,i)} f_{j\,i}^k\otimes g_{j\,i}^k\Big)_{j,i}\,.
 \end{matrix}\right.
\end{split}
\end{equation}
Similarly (using direct sums only) the functor
$\times_{\C}:\Mat{\C}\longrightarrow {\C}$ may be defined.
\abs
\begin{proposition}\label{functor-cat}
Let $({\II},\underline{\Delta}_{\II},
\underline{\varepsilon}_{\II},\underline{R}_{\II})$
be a quasitriangular coalgebra category and
$({\C},\otimes,\E_{\C},\Psi^{\C})$ be
a $\Bbbk$-linear braided abelian category.
Then the category
$(\Bbbk\text{-}\,{\C}^{\II},\otimes_{{\C}^{\II}},
\E_{{\C}^{\II}}, \Psi^{{\C}^{\II}})$
of $\Bbbk$-linear functors is again a $\Bbbk$-linear braided
abelian category. The tensor product is given by
$F\otimes_{\C^\II}G=\bullet_\C\circ(F\uotimes G)\circ\underline{\Delta}_\I$
on the functors and by $\varphi\otimes_{\C^\II}\chi=\bullet_\C\circ
(\varphi\uotimes \chi)\circ\underline{\Delta}_\I$ on the corresponding
functorial morphisms. The unit
object is $\E_{{\C}^{\II}}=\eta_{\C}\circ
\underline{\varepsilon}_{\II}$ with the functor $\eta_{\C}:
\Mat{\Bbbk}\longrightarrow{\C}$, $\eta_{\C}(\Bbbk)=\E_{\C}$,
$\eta_{\C}(\lambda)=\lambda\cdot\id_{\E_{\C}}$, and the braiding
is defined by
\begin{equation}\label{braid-grade}
\Psi^{{\C}^{\II}}_{F,G}=\big(\times_{\C}\diamond
\Psi^{\C}\diamond(F\uotimes G)\circ\underline{\Delta}_{\II}^\op\big)
\circ\big(\bullet_{\C}\circ(F\uotimes G)\diamond\underline{R}_{\II}
\big): F\otimes_{{\C}^{\II}}G\longrightarrow
       G\otimes_{{\C}^{\II}}F\,.
\end{equation}       
The isomorphisms of associativity and of left and right unit multiplication
are defined componentwise through the ones of $\C$.
If $\C$ is in addition $\otimes$-factor$/$$\otimes$-exact abelian then
$\Bbbk\text{-}\,{\C}^{\II}$ is also $\otimes$-factor$/$$\otimes$-exact
abelian.
\end{proposition}

\begin{proof} Since functor categories are abelian if the codomaine
category $\C$ is abelian (see e.\ g.\ \cite{Mac}) one immediately
verifies that $\Bbbk\text{-}\,{\C}^{\II}$ is abelian.
Straightforward computations prove the monoidal structure of
$(\Bbbk\text{-}\,{\C}^{\II},\otimes_{{\C}^{\II}},
\E_{{\C}^{\II}})$, because $\C$ is monoidal.
Without problems the functorial property
$\Psi^{{\C}^{\II}}:\otimes_{{\C}^{\II}}
\overset{\bullet}{\cong}\otimes_{{\C}^{\II}}^\op$ is proven.
It remains to show the hexagon identities for the braiding
$\Psi^{{\C}^{\II}}$. The use of the explicit form of the tensor
product of $\Bbbk\text{-}\,{\C}^{\II}$ leads (up to associativity) to
\begin{equation}\label{braid-grade1}
\begin{split}
\Psi^{{\C}^{\II}}_{F,G\otimes_{{\C}^{\II}}K}
&=\big(\times_{\C}\diamond\Psi^{\C}\diamond
(\id_\C\uotimes\bullet_\C)
\circ(F\uotimes G\uotimes K)\circ (\id_\II\uotimes\underline{\Delta}_\II)
\circ\underline{\Delta}_\II^\op\big)\circ\\
&\quad\circ
\big(\bullet_\C\circ (\id_\C\uotimes\bullet_\C)\circ(F\uotimes G\uotimes K)
\circ(\id_\II\uotimes\underline{\Delta}_\II)\diamond\underline{R}_\II\big)\,.
\end{split}
\end{equation}
On the other hand a similar calculation yields
\begin{equation}\label{braid-grade2}
\begin{split}
\Psi^{\C^\II}_{F,G}\otimes_{\C^\II}\id_K
&=\big(\bullet_\C\diamond(\times_\C\diamond\Psi^\C\uotimes\id_{\id_\C})
\diamond(F\uotimes G\uotimes K)\circ(\underline{\Delta}_\II^\op\uotimes
\id_\II)\circ\underline{\Delta}_\II\big)\circ\\
&\quad\circ
\big(\bullet_\C\circ(\bullet_\C\uotimes\id_\C)\circ(F\uotimes G\uotimes K)
\diamond(\underline{R}_\II\uotimes\id_{\id_\II})\circ\underline{\Delta}_\II
\big)
\end{split}
\end{equation}
and
\begin{equation}\label{braid-grade3}
\begin{split}
&\id_F\otimes_{\C^\II}\Psi^{\C^\II}_{G,K}\\
&=\big(\bullet_\C\diamond(\id_{\id_\C}\uotimes\times_\C\diamond\Psi^\C)
\diamond(\underline{\tau}_{\C,\C}\uotimes\id_\C)\circ
(F\uotimes G\uotimes K)\circ(\underline{\tau}_{\II,\II}\uotimes\id_\II)
\circ(\id_\II\uotimes\underline{\Delta}_\II^\op)\circ
\underline{\Delta}_\II\big)\circ\\
&\quad\circ
\big(\bullet_\C\circ(\id_\C\uotimes\bullet_\C)\circ
(\underline{\tau}_{\C,\C}\uotimes\id_\C)
\circ(F\uotimes G\uotimes K)\circ
(\underline{\tau}_{\II,\II}\uotimes\id_\II)
\diamond(\id_{\id_\II}\uotimes\underline{R}_\II)\circ\underline{\Delta}_\II
\big)
\end{split}
\end{equation} 
where in equation (\ref{braid-grade3}) additionally the involutivity of
$\underline{\tau}$ has been used. We compose (\ref{braid-grade2})
and (\ref{braid-grade3}) and obtain
\begin{equation}\label{braid-grade4}
\begin{split}
&(\id_G\otimes_{\C^\II}\Psi^{\C^\II}_{F,K})\circ
 (\Psi^{\C^\II}_{F,G}\otimes_{\C^\II}\id_K)\\
&=\bullet_\C\diamond
 \Big(\big((\id_{\id_\C}\uotimes\times_\C\diamond\Psi^\C)
\diamond(\underline{\tau}_{\C,\C}\uotimes\id_\C)\big)\diamond
(F\uotimes G\uotimes K)\circ
(\id_\II\uotimes\underline{\Delta}_\II)\circ
\underline{\Delta}_\II^\op\Big)\circ\\
&\quad\circ
\big(\bullet_\C\circ(\id_\C\uotimes\bullet_\C)
\circ(F\uotimes G\uotimes K)\circ(\id_\II\uotimes\underline{\Delta}_\II)
\diamond\underline{R}_\II\big)\\
&=\big(\times_{\C}\diamond\Psi^{\C}\diamond
(\id_\C\uotimes\bullet_\C)
\circ(F\uotimes G\uotimes K)\circ (\id_\II\uotimes\underline{\Delta}_\II)
\circ\underline{\Delta}_\II^\op\big)\circ\\
&\quad\circ
\big(\bullet_\C\circ (\id_\C\uotimes\bullet_\C)\circ(F\uotimes G\uotimes K)
\circ(\id_\II\uotimes\underline{\Delta}_\II)\diamond\underline{R}_\II\big)\\
&=\Psi^{{\C}^{\II}}_{F,G\otimes_{{\C}^{\II}}K}\,.
\end{split}
\end{equation}
In the first equality of (\ref{braid-grade4}) we used the properties
of $\underline{R}_\II$ according to Definition \ref{qtcc}, the
functoriality of $\Psi^\C:\otimes\overset{\bullet}{\cong}\otimes^\op$
and the involutivity of $\underline{\tau}$.
The second equation follows because $\Psi^\C$ is a braiding and
fulfills the hexagon identities. In a similar manner the second hexagon
identity for $\Psi^{\C^\II}$ is verified.
\begin{equation}
\Psi^{\C^\II}_{F\otimes_{\C^\II}G,K}=
(\Psi^{\C^\II}_{F,K}\otimes_{\C^\II}\id_G)\circ
(\id_F\otimes_{\C^\II}\Psi^{\C^\II}_{G,K})\,.
\end{equation}
If $\C$ is $\otimes$-factor abelian and
$\phi_j:F_j\overset{\bullet}{\to} G_j$, $j\in\{1,2\}$ are two epimorphisms
in $\Bbbk\text{-}\,{\C}^{\II}$, then one verifies easily componentwise
that $(\phi_1\otimes_{\C^\II}\phi_2)$ is epimorphic in
$\Bbbk\text{-}\,{\C}^{\II}$. Similarly the $\otimes$-monomorphy
property will be checked. Hence $\Bbbk\text{-}\,\C^\II$ is $\otimes$-factor
abelian.

Suppose that $\C$ is $\otimes$-exact abelian. Let again
$\phi_j:F_j\overset{\bullet}{\to} G_j$, $j\in\{1,2\}$ be two epimorphisms
in $\Bbbk\text{-}\,{\C}^{\II}$. For an object $\i\in\II$ denote
$\underline{\Delta}_\II(\i)=\big((\i_1^{(1)},\i_1^{(2)}),\ldots,
(\i_n^{(1)},\i_n^{(2)})\big)$. If there are morphisms
$\rho_1:G_1\otimes_{\C^\II}F_2\overset{\bullet}{\to}
G_1\otimes_{\C^\II}G_2$ and $\rho_2:F_1\otimes_{\C^\II}G_2
\overset{\bullet}{\to}G_1\otimes_{\C^\II}G_2$ such that
$\rho_1\circ(\phi_1\otimes_{\C^\II}\id_{F_2})=
\rho_2\circ(\id_{F_1}\otimes_{\C^\II}\phi_2)$, then it follows
$\rho_1(\i)_j\circ(\phi_1(\i^{(1)}_j)\otimes_{\C}\id_{F_2(\i^{(2)}_j)})=
 \rho_2(\i)_j\circ(\id_{F_1(\i^{(1)}_j)}\otimes_{\C}\phi_2(\i^{(2)}_j))$.
There are unique morphisms $\rho(\i)_j$ such that
$\rho(\i)_j\circ(\id_{G_1(\i^{(1)}_j)}\otimes_{\C}\phi_2(\i^{(2)}_j))=
\rho_1(\i)_j$ and $\rho(\i)_j\circ(\phi_1(\i^{(1)}_j)\otimes_{\C}
\id_{G_2(\i^{(2)}_j)})=\rho_2(\i)_j$ because $\C$ is $\otimes$-exact
abelian by assumption. From this fact it follows immediately that $\rho$
is a natural morphism and that the $\otimes$-right-exact condition of
Definition \ref{tensor-exact} is fulfilled for $\Bbbk\text{-}\,\C^\II$.
Analogously the $\otimes$-left-exact property can be shown. Hence
$\Bbbk\text{-}\,\C^\II$ is $\otimes$-exact abelian.\end{proof}
\abs
Every coalgebra morphism in $\Bbbk\text{-}\Mcat$ induces an exact
$\Bbbk$-linear monoidal functor on the corresponding functor categories.
More precisely it holds
\abs
\begin{corollary}\label{ind-grade-funkt}
Let $\II$ and $\JJ$ be (quasitriangular) coalgebra
categories and $\Theta:\II\longrightarrow\JJ$
be a coalgebra functor in $\Bbbk\text{-}\Mcat$ (such that
$(\Theta\uotimes\Theta)\diamond \underline{R}_\II=\underline{R}_\JJ\diamond
\Theta$). If $(\C,\otimes,\E_\C,(\Psi^\C))$ is a $\Bbbk$-linear
(braided) monoidal abelian category then an exact $\Bbbk$-linear
(braided) monoidal functor $\C^\Theta:\Bbbk\text{-}\,\C^\JJ\longrightarrow
\Bbbk\text{-}\,\C^\II$ can be defined by
\begin{equation}\label{ind-grade-funkt1}
\C^\Theta\ :\ \left\{ \begin{matrix}
                    F & \mapsto & \times_\C\circ F\circ\Theta\\
              \varphi & \mapsto & \times_\C\diamond\varphi\diamond\Theta
                    \end{matrix}\right.
\end{equation}
where $F$ is a functor in $\Ob(\Bbbk\text{-}\,\C^\JJ)$ and $\varphi$ is a
functorial morphism in $\Bbbk\text{-}\,\C^\JJ$. If $\Theta$ is an
equivalence of categories then $\C^\Theta$ is an equivalence.
\end{corollary}

\begin{proof}
One proves without difficulties that $\C^\Theta$ is a
$\Bbbk$-linear functor. Since $\C^\II$ and $\C^\JJ$ are abelian
the exactness of $\C^\Theta$ is shown on the several components of the
natural morphisms. The coalgebra morphism properties
$\underline{\Delta}_\JJ\circ\Theta= (\Theta\uotimes\Theta)\circ
\underline{\Delta}_\II$ and $\underline{\varepsilon}_\JJ\circ\Theta=
\underline{\varepsilon}_\II$ show immediately that $\C^\Theta$ is monoidal.

If $(\Theta\uotimes\Theta)\diamond \underline{R}_\II=\underline{R}_\JJ
\diamond\Theta$ then the following equations hold.
\begin{equation}
\begin{split}
 \C^\Theta(\Psi^{\C^\JJ}_{F,G})&=\times_\C\diamond\Big(\big(\times_{\C}
 \diamond\Psi^\C\diamond(F\uotimes G)\circ\underline{\Delta}_\JJ^\op\big)
 \circ\big(\bullet_\C\circ(F\uotimes G)\diamond\underline{R}_\JJ\big)\Big)\\
 &=\big(\times_{\C}
 \diamond\Psi^\C\diamond(F\circ\Theta\uotimes G\circ\Theta)\circ
 \underline{\Delta}_\II^\op\big)\circ\big(\bullet_\C\circ
 (F\circ\Theta\uotimes G\circ\Theta)\diamond\underline{R}_\II\big)\\
 &=\big(\times_{\C}
 \diamond\Psi^\C\diamond(\C^\Theta(F)\uotimes \C^\Theta(G))\circ
 \underline{\Delta}_\II^\op\big)\circ\big(\bullet_\C\circ
 (\C^\Theta(F)\uotimes \C^\Theta(G))\diamond\underline{R}_\II\big)\\
 &=\Psi^{\C^\II}_{\C^\Theta(F),\C^\Theta(G)}\,.
\end{split}
\end{equation} 
Thus the functor $\C^\Theta$ respects the braiding.\end{proof}
\abs
The investigations of braided graded categories and of categories of
complexes require special $\Bbbk$-linear
categories $\uNN\subset\uNN_{c}$. Their objects are the natural
numbers $\NN_0=\{0,1,2,\ldots\}$, and the morphisms are given by
\begin{equation*}
\Hom_{\uNN}(m,n)=\cases\Bbbk&\text{if $m=n$,}\\
                       \{0\}&\text{if $m\ne n$,}\endcases\quad
                       \text{and}\quad
\Hom_{\uNN_{c}}(m,n)=\cases\Bbbk&\text{if $m=n$,}\\
                         \Bbbk\cdot\partial_m&\text{if $n=m+1$},\\
                         \{0\}&\text{else.}\endcases
\end{equation*}
The compositions are defined by $\Bbbk$-linearity and the
generating relations $\partial_{n+1}\cdot\partial_n=0$ for all $n\in\NN_0$
in $\uNN_c$. The category $\uNN$ is a quasitriangular coalgebra category
through the $\Bbbk$-linear functors $\underline\Delta$ and
$\underline\varepsilon$ defined on every object $n\in\Ob(\uNN)$ by
\begin{equation}
\label{Delta-Ob}
\underline{\Delta}(n)=\left( (0,n),(1,n-1),\dots,(n,0)\right)\quad
 \text{and}\quad
\underline{\varepsilon}(n)=\cases \Bbbk\in\Ob(\Mat{\Bbbk})&\text{if $n=0$}
            \\
            \emptyset\in\Ob(\Mat{\Bbbk})&\text{if $n\ne0$}\endcases
\end{equation}
and by the obvious definition of $\underline{\Delta}$ and
$\underline{\varepsilon}$ on the morphisms. Compatible
quasitriangular structures $\underline{R}_n:\Delta(n)\overset{\bullet}{\to}
\Delta^\op(n)$ on $\uNN$ are characterized by invertible elements
$\lambda\in\Bbbk^*$. The nonzero components of the corresponding morphisms
$\underline{R}^{(\lambda)}_n$ are given by
$(\underline{R}^{(\lambda)}_n)_{(\ell,k),(k,\ell)}=\lambda^{k\,\ell}$
for $k+\ell=n$. Hence $(\uNN,\underline{\Delta},\underline{\varepsilon},
\underline{R}^{(\lambda)})$ is a quasitriangular coalgebra category
for any $\lambda\in \Bbbk^*$ . The next lemma states that there are also
quasitriangular structures on $\uNN_{c}$.
\abs
\begin{lemma}
\label{Delta-Mor} The category
$(\uNN_{c},\underline{\Delta}_c,\underline{\varepsilon}_c,
\underline{R}_c)$ is a quasitriangular coalgebra category where the
coalgebra structure on $\uNN_{c}$ is given by
${\underline{\Delta}_c}_{\vert \uNN}=\underline{\Delta}\,$,
${\underline{\varepsilon}_c}_{\vert \uNN}=\underline{\varepsilon}$ and
$\underline{\Delta}_c(\partial_n)_{(r,s),(k,\ell)}=
\delta_{(r,s),(k+1,\ell)}\,\partial_k\otimes_\Bbbk\id_\ell+
\delta_{(r,s),(k,\ell+1)}\,(-1)^k\,\id_k\otimes_\Bbbk\partial_\ell$ and
$\underline{\varepsilon}_c(\partial_n)=0$. The quasitriangular isomorphism
$\underline{R}_c$ coincides with $\underline{R}^{(-1)}$.
\end{lemma}

\begin{proof}
Obviously $\underline{\Delta}_c$ is a functor because of
the definition of $\underline{\Delta}_c(\partial_n)$ which is adapted to
the identity $\partial_{n+1}\circ\partial_n=0$. For the prove of the
functoriality of $\underline{R}^{(-1)}$ we note that the nonzero elements
of both matrices $\underline{R}^{(-1)}_{n+1}\circ\Delta(\partial_n)$ and
$\Delta^\op(\partial_n)\circ\underline{R}^{(-1)}_n$ coincide and are
given by $(-1)^{(k+1)\ell}\cdot\partial_k\otimes\id_\ell$ and
$(-1)^{k\ell}\cdot\id_k\otimes\partial_\ell$.
\end{proof}
\abs
\begin{remark}{\rm
A more general form for $\underline{\Delta}_c(\partial_n)$ which does not
necessarily induce a quasitriangular structure on $\uNN_{c}$, is the
morphism $\underline{\Delta}_c(\partial_n)_{(r,s),(k,\ell)}=
\delta_{(r,s),(k+1,\ell)}\,\alpha_{k,\ell}\cdot
\partial_k\otimes_\Bbbk\id_\ell+
\delta_{(r,s),(k,\ell+1)}\,(-1)^k\,\beta_{k,\ell}\cdot
\id_k\otimes_\Bbbk\partial_\ell$ where
$\alpha_{k,\ell},\beta_{k,\ell}\in\Bbbk$. The coassociativity of
$\underline{\Delta}_c$ and the identity $\partial_{n+1}\circ\partial_n=0$
are equivalent to $\alpha_{k,\ell}=\alpha_{k+1}\cdot\alpha_{k+2}\cdot
\ldots\cdot\alpha_{k+\ell}$ and $\beta_{k,\ell}=\beta_{\ell+1}\cdot
\beta_{\ell+2}\cdot\ldots\cdot\beta_{\ell+k}$
and $\alpha_{k,\ell}\cdot\beta_{k+1,\ell}+
\beta_{k,\ell}\cdot\alpha_{k,\ell+1}=\alpha_{k+1}\cdot\ldots\cdot
\alpha_{k+\ell}\cdot\beta_{\ell+1}\cdot\ldots\cdot\beta_{\ell+k}\cdot
(\alpha_{k+\ell+1}+\beta_{k+\ell+1})=0$, for any $\alpha_i,\,\beta_i\in
\Bbbk$. In the case where all $\alpha_i,\,\beta_i$ are invertible, these
conditions imply that this coalgebra structure is equivalent to the one
described in Lemma \ref{Delta-Mor}. There also exist other solutions for 
$\underline{\Delta}_c$ which are deviating from the structure of
$\underline{\Delta}_c$ in Lemma \ref{Delta-Mor}. For example the elements
$\alpha_1=\beta_1=0$, $\alpha_2=1$, $\beta_2$ arbitrary, and
$\alpha_n=-\beta_n=1$ for $n\ge 3$ fulfill the conditions.}\end{remark}
\abs
\subsection*{Graded Categories and Complexes}

We return to the initial conditions fixed in Chapter \ref{ideal-section}
and direct our attention therefore to $\otimes$-factor
braided abelian categories $\C$ and to the commutative ring
$\Bbbk:=\End_{\C}(\E_\C)$. The full subcategories ${\mathbf 2}$ and
${\mathbf 2}_c$ of $\uNN$ and
$\uNN_{c}$ respectively, with objects $\{0,1\}$ are quasitriangular
colagebra subcategories. Denote by\/ $\I$ either the discrete category\/
$\I={\mathbf 2}$ or\/ $\I=\uNN$ and by $\I_c$ the categories
$\I_c={\mathbf2}_c$ or $\I_c=\uNN_{c}$. The $\I$-graded category
${\C}^\I$ over $\C$ and the category ${\C}^{\I_c}$ of
$\I$-graded complexes are the functor categories which we are interested in
in the sequel (see e.\ g.\ \cite{Hus,Mac}). The categories
${\C}^\I$, ${\C}^{\I_c}$ are $\otimes$-factor braided
abelian, because of Proposition \ref{functor-cat}. For completeness we will
give below the explicit braided monoidal structure of these categories.
\abs
\begin{enumerate}
\item
The objects of the graded category ${\C}^\I$ over $\I={\mathbf 2}$ or
$\I=\uNN$ are given by $\hat X=(X_0, X_1,\ldots)$ where
$X_j\in \Ob({\C})$ for all $j\in \I$. The morphisms of ${\C}^\I$
are of the form $\hat f=(f_0,f_1,\ldots):\hat X\to \hat Y$ where
$f_j:X_j\to Y_j$ is a morphism in $\C$ for all $j\in \I$.
\item
The category $\C^{\I_c}$ of\/ $\I$-graded complexes consists of
objects $(\hat X, \hat \dif)$ where $\hat X\in\Ob({\C}^\I)$,
and $\hat \dif=(\dif_0,\dif_1,\ldots)$ is a differential such that
for the morphisms $\dif_j:X_j\to X_{j+1}$ in $\C$ it holds
${\dif_{j+1}\circ \dif_{j}=0}$ for all $j\in\I$ if\/ $\I=\uNN$.
The morphisms $\hat f:(\hat X,\hat \dif)\to (\hat Y,\hat \dif')$
in $\C^{\I_c}$ are morphisms in ${\C}^\I$ which obey the relations
$f_{j+1}\circ \dif_j= \dif'_j\circ f_j$ $\forall\ j, j+1\in \I$.
\item
The unit objects of $\C^\I$ and $\C^{\I_c}$ are both given by
$\Hat\E=(\E_\C,0,0,\ldots)$. The tensor product is defined through
$(\hat X\otimes\hat Y)_n=\bigoplus_{k+l=n}X_k\otimes Y_l$ where
$n,k,l\in\I$, and the tensor product of\/ $\I$-graded and complex
morphisms is built analogously. The differential of a tensor product in
$\C^{\I_c}$ is of the form $(\hat\dif_{\hat X\otimes\hat Y})_n=
\sum_{k+l=n}[\dif_{\hat X,k}\otimes\id_{Y_l}+ (-1)^k\,\id_{X_k}\otimes
\dif_{\hat Y,l}]$. The category ${\C}^\I$ admits a family of
braidings given by $(\hat \Psi_{\hat X,\hat Y}^{(\lambda)})_n=
\bigoplus_{k+l=n}\lambda^{k\,l}\Psi_{X_k,Y_l}$ for any
$\lambda\in \Aut_\C(\E_\C)$\footnote{The
braidings in ${\C}^\I$ and $\C^{\I_c}$ are bi-additive,
$\hat \Psi^{(\lambda)}_{\hat X\oplus\hat Y, \hat Z}
 =\hat \Psi^{(\lambda)}_{\hat X,\hat Z}\oplus
  \hat \Psi^{(\lambda)}_{\hat Y,\hat Z}$, because the braiding in $\C$
is bi-additive.}. The braiding in $\C^{\I_c}$ is given by
$(\hat \Psi_{\hat X,\hat Y}^{(-1)})$ forced by the nilpotency of the
differentials $\hat\dif$.
\end{enumerate}
\abs
Consider the category $\uNN$ with the quasitriangular structure
$\underline{R}^{(-1)}$. Then there are coalgebra category morphisms 
$\Bbbk\leftrightarrows\uNN\leftrightarrows\uNN_c$ which are compatible
with $\underline{R}^{(-1)}$. The right arrows are
natural inclusions where $\Bbbk$ is identified with $\End(0)$. The
functor $\uNN\to\Bbbk$ is the counit, and
$\uNN_c\to\uNN$ is identical on objects and sends $\partial_n$ to 0.
As an implication of Corollary \ref{ind-grade-funkt}
the following exact braided monoidal functors between $\C$,
${\C}^\I$ and ${\C}^{\I_c}$ can be deduced.
\begin{equation}
\label{1-1b-1ZD11}
{\C}\leftrightarrows{\C}^\I\leftrightarrows{\C}^{\I_c}
\end{equation}
The right arrows are the well known canonical inclusions
given by the assignments $X\mapsto \hat X:=(X,0,0,\ldots)$ for all
$X\in \Ob({\C})$, and $\hat Y\mapsto (\hat Y, \hat\dif :=0)$ for all
$\hat Y\in\Ob({\C}^\I)$. Analogously the assignments for the morphisms
will be defined. In the sequel we will use this identification for
embedding the objects and morphisms of the several categories into the
graded category or the category of complexes respectively.
\abs
\begin{remark}{\rm We have to give a comment on the usual notation of
matrix elements of graded morphisms which we use subsequently.
If $\hat f:\hat X\to\hat Y$ is a morphism in ${\C}^\I$
and for all $n\in\I$ the components
$X_n=\bigoplus_{k+l=n}U_k^{(n)}\otimes V_l^{(n)}$ and
$Y_n=\bigoplus_{k+l=n}W_k^{(n)}\otimes Z_l^{(n)}$ are direct sums of
tensor products we denote by
\begin{equation}\label{matrix-elem1}
f_{r,s;k,l}:U_k^{(n)}\otimes V_l^{(n)}\to W_r^{(n)}\otimes Z_s^{(n)}
\end{equation}
the corresponding matrix element of $\hat f$. If it is clear from the
context we denote by $f_{m,n}$ either the morphism
\begin{equation}\label{matrix-elem2}
f_{m,n}:=\left(\begin{matrix}&f_{0,m+n;m,n}\\ &\vdots\\ &f_{m+n,0;m,n}
         \end{matrix}\right): U_m^{m+n}\otimes V_n^{m+n}\to Y_{m+n}
\end{equation}
or the morphism          
\begin{equation}\label{matrix-elem3}
f_{m,n}:=\big(f_{m,n;0,m+n},\ldots , f_{m,n;m+n,0}\big)
         : X_{m+n}\to W_m^{m+n}\otimes Z_n^{m+n}
\end{equation}
Of course, the same notation is used if only $\hat X$ or $\hat Y$
consists of components which are direct sums of tensor products.
One proceeds analogously if there are more than two tensor factors.}
\end{remark}
\abs
In the following we are investigating algebraic structures on objects
in graded or complex categories. Similar considerations as in
\cite{Mal,Man,Wor} immediately lead to 
\abs
\begin{proposition}\label{12-8IIIZD34}
Let $\hat X$ be an object in the $\otimes$-factor braided abelian category
${\C}^\I$ where\/ $\I={\mathbf 2}$ or\/ $\I=\uNN$. Then it holds in
particular
\begin{enumerate}
\item
If $\hat X$ is an algebra in ${\C}^\I$ then $X_0$ is an algebra
in ${\C}$ and $X_n$ is an $X_0$-bimodule for any $n\in\I$.
\item
If $\hat X$ is a coalgebra in ${\C}^\I$ then $X_0$ is a coalgebra
in ${\C}$ and $X_n$ is an $X_0$-bicomodule for any $n\in\I$.
\item
If $\hat X$ is a bi- or Hopf algebra in ${\C}^\I$ then $X_0$ is a bi-
or Hopf algebra in ${\C}$ respectively and $X_n$ is an
$X_0$-Hopf bimodule for any $n\in\I$.
\end{enumerate}
\end{proposition}

\begin{proof}
Suppose that $(\hat X,\hat\m,\hat\eta)$ is an algebra in ${\C}^\I$
then $\hat\m\circ(\id_{\hat X}\otimes\hat\m)=
 \hat\m\circ(\hat\m\otimes\id_{\hat X})$. Looking in particular at the
matrix components $(0,0,0)$ and $(0,0,n)$, $(n,0,0)$, $(0,n,0)$ of this
equation, it
follows that $(X_0,\m_{0,0})$ is an algebra and $(X_n,\m_{n,0},\m_{0,n})$
is an $X_0$-bimodule respectively. The unital property of $X_0$ and $X_n$
is given by the $0$th and $n$th component of the equations
$\hat\m\circ(\id_{\hat X}\otimes\hat\eta)=\id_{\hat X}=
\hat\m\circ(\hat\eta\otimes\id_{\hat X})$ respectively.
The dual statement for coalgebras is obtained analogously.
In the case that $\hat X$ is a bialgebra the additional relation
$\hat\Delta\circ\hat\m=(\hat\m\otimes\hat\m)\circ(\hat\id\otimes
\hat\Psi^{(\lambda)}\otimes\hat\id)\circ(\hat\Delta\otimes\hat\Delta)$
leads to the third statement of the proposition.\end{proof}
\abs
\begin{remark}\label{mod-alg-equiv}
{\rm If $\I={\mathbf 2}$ the particular statements (1), (2) and (3) in
Proposition \ref{12-8IIIZD34} are equivalences.}
\end{remark}
\abs
\begin{remark}
{\rm For complexes we obtain analogous statements as in
Propositon \ref{12-8IIIZD34}. The differential $\hat\dif$ of a coalgebra
$(\hat X,\hat\dif)\in\Ob(\C^{\I_c})$ is an $X_0$-bicomodule morphisms.
If $(\hat X,\hat\dif)$ is an algebra in $\C^{\I_c}$ then the
differential $\hat\dif$ obeys a kind of Leibniz rule since the
multiplication of $(\hat X,\hat\dif)$ is a complex morphism by assumption.
In particular it holds $\dif_0\circ\m_{X_0}=\mu^{X_1}_r\circ
(\dif_0\otimes\id_{X_0})+\mu_l^{X_1}\circ(\id_{X_0}\otimes\dif_0)$, where
$\mu_r^{X_1}$ and $\mu_l^{X_1}$ are the right and left action of $X_0$ on
$X_1$ respectively.}
\end{remark}
\abs
\begin{remark}\label{bialg-proj-canon}
{\rm If $\Hat X$ is a Hopf algebra in ${\C}^\I$
and $X_0$ is flat with isomorphic antipode in $\C$ then the canonical
morphisms $X_0\to \Hat X\to X_0$ form a Hopf algebra projection
on $X_0$ in ${\C}^\I$. Hence $(\underline{\Hat X},
\underline\m_{\Hat X},\underline\eta_{\Hat X},\underline\Delta_{\Hat X},
\underline\varepsilon_{\Hat X},\underline S_{\Hat X})$ according to
Theorem \ref{bialg-proj} is a Hopf
algebra in ${}^{X_0}_{X_0}({\C}^\I)^{X_0}_{X_0}$.}
\end{remark}
\abs
Using Proposition \ref{12-8IIIZD34} we derive necessary and sufficient
criteria for a graded bialgebra to be a Hopf algebra. 
This is a generalization of the ideas of \cite{Ma5}. The antipode of the
form (\ref{antipode}) has been used implicitely in \cite{Brz} in the
special case of Woronowicz's exterior algebra $\Gamma^\wedge$.
\abs
\begin{proposition}\label{7-8-6ZD11}
Let $(\hat B,\hat\m,\hat\Delta)$ be a bialgebra in ${\C}^\I$ or in
$\C^{\I_c}$. Then $\hat B$ is a Hopf algebra if and only if $B_0$
(with the bialgebra structure morphisms induced from $\hat B$) is
a Hopf algebra in $\C$. In this case the antipode of $\hat B$ is
successively given by
\begin{equation}\label{antipode}
S_n=-\sum_{k=1}^n \left(\m^{(2)}_{0,k,n-k}\circ
 (S_0\otimes\id_{B_k}\otimes S_{n-k})\circ
 \Delta^{(2)}_{0,k,n-k}\right)
\end{equation}
for all $n\in\I$, where
${\hat\m}^{(2)}=\hat\m\circ(\hat\id\otimes\hat\m)$ and
${\hat\Delta}^{(2)}=(\hat\id\otimes\hat\Delta)\circ\hat\Delta$.
\end{proposition}

\begin{proof}
The sufficiency has been proved in Proposition \ref{12-8IIIZD34}.3.
Suppose that $B_0$ is a Hopf algebra in $\C$ and $\hat B$ is a
bialgebra in $\C^\I$. We investigate
the morphism $\hat S$ in (\ref{antipode}) and show that it is an antipode
of $\hat B$. By assumption it holds $\m_{0,0}\circ(\id_{B_0}\otimes S_0)
\circ\Delta_{0,0}=\eta_0\circ\varepsilon_0=\m_{0,0}\circ
(S_0\otimes\id_{B_0})\circ\Delta_{0,0}$. For $n>0$ we obtain the following
equations with $\hat\m^{(3)}:=\hat\m\circ(\hat\id_{\hat B}\otimes\hat\m)
\circ(\hat\id_{\hat B}\otimes\hat\id_{\hat B}\otimes\hat m)$, and
$\hat\Delta^{(3)}$ defined in the dual analogous way.
\begin{equation}
\begin{split}\label{antipode2}
\big(\hat m\circ(\hat\id_{\hat B}\otimes\hat S)\circ\hat\Delta\big)_n
&= \m_{0,n}\circ(\id_{B_0}\otimes S_n)\circ\Delta_{0,n}+
\sum_{k=1}^n \m_{k,n-k}\circ(\id_{B_k}\otimes S_{n-k})\circ\Delta_{k,n-k}\\
&=\sum_{k=1}^n \Big(-\m^{(3)}_{0,0,k,n-k}\circ(\id_{B_0}\otimes S_0
 \otimes\id_{B_k}\otimes S_{n-k})\circ\Delta^{(3)}_{0,0,k,n-k}+\\
&\phantom{=}\,\,+\m_{k,n-k}\circ(\id_{B_k}\otimes S_{n-k})\circ
 \Delta_{k,n-k}\Big)\\
&=\sum_{k=1}^n (-S_{n-k}+S_{n-k})=0
\end{split}
\end{equation}
where we used (\ref{antipode}) in the second equation for $S_n$ and the
(co-)associativity and the antipode property of $S_0$ in the third
equation. Thus we obtain $\hat\m\circ(\hat\id_{\hat B}\otimes\hat S)
\circ\hat\Delta=\hat\eta\circ\hat\varepsilon$. Similarly
$\hat\m\circ(\hat S\otimes\hat\id_{\hat B})\circ\hat\Delta
=\hat\eta\circ\hat\varepsilon$ can be proven.

If $(\hat B,\hat \dif)$ is a bialgebra in $\C^{\I_c}$ and $B_0$ is a Hopf
algebra in $\C$ it follows that $\hat B$ is a Hopf algebra in $\C^\I$
according to the previous considerations. Then
\begin{equation}\label{antipode4}
\begin{split}
\hat S\circ\hat\dif &= \hat\m\circ (\hat S\otimes \hat\eta\circ
                     \hat\varepsilon)\circ\hat\Delta\circ\hat\dif\\
 &= \hat\m\circ (\hat S\otimes\hat\id)\circ \hat\dif^\otimes\circ
    (\hat\id\otimes\hat\eta\circ\hat\varepsilon)\circ\hat\Delta\\
 &= \hat\m\circ (\hat S\otimes\hat\id)\circ (\id\otimes\hat\m)\circ
    (\hat\Delta\otimes\hat\id)\circ\hat\dif^\otimes\circ
    (\hat\id\otimes\hat S)\circ\hat\Delta\\
 &=(\hat\varepsilon\otimes\hat\id)\circ\hat\dif^\otimes\circ
   (\hat\id\otimes\hat S)\circ\hat\Delta\\
 &=\hat\dif\circ \hat S\,.
\end{split}   
\end{equation}
In the second, the third and the fifth equation of (\ref{antipode4}) we
used that $(\hat B,\hat\dif)$ is a bialgebra in $\C^{\I_c}$.\end{proof}
\abs
\begin{remark}{\rm In (\ref{antipode4}) we only have to require that
$(\hat B,\hat\dif)$ is an algebra and a coalgebra in $\C^{\I_c}$. The
techniques we applied in (\ref{antipode4}) can also be used in particular
to show that any Hopf algebra which is a (co-)module bialgebra
is automatically a (co-)module Hopf algebra.}
\end{remark}
\abs
\subsection*{Exterior Hopf Algebras}

Next we consider the tensor algebra generated by an object
$X\in\Ob({\C})$. It is commonly defined as the graded object
$T_{\C}(X)$ in ${\C}^{\uNN}$ with
$(T_{\C}(X))_0:=\E_{\C}$ and $(T_{\C}(X))_j:=X^{\otimes\,j}$
the $j$-fold tensor product of $X$ for $j\ge 1$.
Before we construct Hopf algebra structures on $T_{\C}(X)$
we outline braided combinatorics which generalizes the usual
symmetric combinatorics. We extend the ideas of \cite{Ma4,Wor}.
Since $({\C},\E,\otimes,\Psi)$ is a braided category it is possible
to define a mapping of the symmetric group $S_j$ into representations
of the braid group $B_j$ \cite{Art} on the object $X^{\otimes\,j}$
for any $j\in \NN_0$. For that we look at the so-called reduced expression.
A reduced expression of a permutation $\sigma\in S_j$ is a product
decomposition of $\sigma$ into a minimal number of next neighbour
transpositions $t_a$, $a\in\{1,\ldots ,j-1\}$ which permute $a$ with
$a+1$. Two such decompositions of $\sigma$ can be transformed into each
other by applying braid group relations only \cite{Wor}. The involutivity
$t_a^2=\e$ will not be used. The minimal number is the length
$\ell(\sigma)$ of the permutation $\sigma$. Hence there is a well defined
mapping $S_j\to B_j$ given by
$\sigma=t_{a_1}\,t_{a_2}\cdots t_{a_{\ell(\sigma)}} \mapsto
 \psi_{a_1}\,\psi_{a_2}\cdots \psi_{a_{\ell(\sigma)}}$, where the set
$\{\psi_a\}_{a\in\{1,\ldots ,j-1\}}$ consists of the elementary generators
of the braid group $B_j$ \cite{Art}. The canonical composition
$S_j\to B_j\to S_j$ is obviously
the identity on $S_j$ and therefore the mapping $S_j\to B_j$ is a section.
Since the category $\C$ is braided, one can define
the mapping $S_j\to \End(X^{\otimes\,j})$ by
$\sigma=t_{a_1}\,t_{a_2}\cdots t_{a_{\ell(\sigma)}} \mapsto
 \sigma_{\C}(X):={(\Psi_{X,X})}_{a_1}\,{(\Psi_{X,X})}_{a_2}\cdots
 {(\Psi_{X,X})}_{a_{\ell(\sigma)}}$ where
${(\Psi_{X,X})}_{a}=\id_{X^{\otimes\,a-1}}\otimes\Psi_{X,X}\otimes
\id_{X^{\otimes\,n-a-1}}$ for all $a\in\{1,\ldots,n-1\}$ \cite{Wor}.
Now let $\pi=(j_1,\ldots,j_r)$ be any $\NN_0$-partition of $j$, i.\ e.\ 
$j=j_1+\cdots +j_r$ and $j_1,\ldots j_r\in\NN_0$.
We consider the shuffle permutations
$S^j_{\pi}\subset S_j$. For every $k\in\{1,\ldots,r\}$ they are mapping
$j_k$ elements of $\{1,\ldots ,j\}$ to
$\{(\sum_{l=1}^{k-1} j_l)+1,\ldots,(\sum_{l=1}^{k} j_l)\}$
without changing their
order. The set of inverse permutations of $S^j_{\pi}$ is denoted by
$S_j^{\pi}$. Based on the ideas of \cite{Ma4,Wor} we define braided
multinomials for every partition $\pi=(j_1,\ldots ,j_r)$ of $j$ and any
object $X\in\Ob({\C})$ according to
\begin{equation}\label{multinom1}
\Big[{\pi\atop j}\Big\vert X;\lambda\Big]=
\Big[{{j_1\dots j_r}\atop j}\Big\vert X;\lambda\Big] :=
\sum_{\sigma\in S_j^\pi}\lambda^{\ell(\sigma)}\,\sigma_{\C}(X)\,.
\end{equation}
where $\lambda\in\Aut(\E)$ is an automorphism of the unit object in
$\C$. The dual counterparts of (\ref{multinom1}) are given by
\begin{equation}\label{multinom2}
\Big[{j\atop \pi}\Big\vert X;\lambda\Big]=
\Big[{j\atop {j_1\dots j_r}}\Big\vert X;\lambda\Big] :=
\sum_{\sigma\in S^j_\pi}\lambda^{\ell(\sigma)}\,\sigma_{\C}(X)\,.
\end{equation}
Both types are endomorphisms of $X^{\otimes\,j}$ in $\C$. It is not
difficult to prove that
$\big[{j\atop j}\big\vert X;\lambda\big]=\id_{X^{\otimes\,j}}$, and
$\big[j\big\vert X;\lambda\big]!:=\big[{j\atop 1\dots 1}
 \big\vert X;\lambda\big]=\big[{1\dots 1\atop j}\big\vert X;\lambda\big]$
which corresponds to the antisymmetrizer $A_j$ in \cite{Wor}.
We obtain the braided form of the fundamental result on multinomials
of subpartitions.
\abs
\begin{proposition}\label{combi}
Let $\pi=(j_1,\dots, j_r)$ be a partition of $j$, and let
$\pi_k=(j_1^k,\dots ,j_{r_k}^k)$ be a partition of $j_k$ for any
$k\in\{1,\dots ,r\}$. With the notation from above it holds
\begin{equation}
\label{combi-eq1}
\Big[{(\pi_1,\pi_2,\dots ,\pi_r)\atop j}\Big\vert X;\lambda\Big]
=\Big[{\pi\atop j}\Big\vert X;\lambda\Big]\circ
   \Big(\Big[{\pi_1\atop j_1}\Big\vert X;\lambda\Big]\otimes
       \Big[{\pi_2\atop j_2}\Big\vert X;\lambda\Big]\otimes\cdots
      \otimes
         \Big[{\pi_r\atop j_r}\Big\vert X;\lambda\Big]\Big)
\end{equation}
and
\begin{equation}         
  \label{combi-eq2}
\Big[{j\atop (\pi_1,\pi_2,\dots ,\pi_r)}\Big\vert X;\lambda\Big]
=\Big(\Big[{j_1\atop \pi_1}\Big\vert X;\lambda\Big]\otimes
      \Big[{j_2\atop \pi_2}\Big\vert X;\lambda\Big]\otimes\cdots\otimes
         \Big[{j_r\atop \pi_r}\Big\vert X;\lambda\Big]\Big)\circ
         \Big[{j\atop\pi}\Big\vert X;\lambda\Big]\,.
\end{equation}
Equation (\ref{combi-eq2}) is the dual version of (\ref{combi-eq1}).
\end{proposition}

\begin{proof}
Let $j\in\NN_0$ and $\pi=(j_1,\ldots,j_r)$ be a partition of $j$.
Similarly as in \cite{Wor} one proves that every permutation
$\sigma\in S_j$ decomposes uniquely according to
\begin{equation}\label{decompose-perm}
\sigma=\sigma_{(j_1,\ldots,j_r)}\circ p_1\circ\cdots \circ p_r
\end{equation}
where $\sigma_{(j_1,\ldots,j_r)}\in S^j_\pi$ and for any
$k\in\{1,\ldots,r\}$ the $p_k$ permutes the elements
${\{(\sum_{l=1}^{k-1} j_l)+1,\ldots,(\sum_{l=1}^{k} j_l)\}} 
\subset\{1,\ldots,j\}$ only, i.\ e.\ 
$p_k$ is the tensor product (in the category of symmetric groups) of the
identity in $S_{j_1+\cdots +j_{k-1}}$ with a permutation in
$S_{j_k}$ and the identity in $S_{j-(j_1+\cdots +j_{k})}$.
Moreover the length of $\sigma $ is given by
$\ell(\sigma)=\ell(\sigma_{(j_1,\ldots,j_r)})+\ell(p_1)+\cdots+\ell(p_r)$.
Suppose that $\pi_k=(j_1^k,\dots ,j_{r_k}^k)$ is a partition of $j_k$
for any $k\in\{1,\dots ,r\}$. Then every $\sigma\in
S^j_{(\pi_1,\ldots,\pi_r)}$ decomposes uniquely like in
(\ref{decompose-perm}) where the nontrivial action of the permutation $p_k$
is given by a shuffle permutation in $S^{j_k}_{\pi_k}$.
Inserting these results in the Definition \ref{multinom1} 
for $\big[{j\atop (\pi_1,\ldots,\pi_r)}\big\vert X;\lambda\big]$
yields (\ref{combi-eq2}). The dual statement (\ref{combi-eq1}) is obtained
similarly.\end{proof}
\abs
For any $0\le k\le j$ we define the morphisms
$\big[{j\atop k}\big\vert X;\lambda\big]
:=\big[{j\atop {k , j-k}}\big\vert X;\lambda\big]$ and
$\big[{k\atop j}\big\vert X;\lambda\big]:=
\big[{{k , j-k}\atop j}\big\vert X;\lambda\big]$.
Then we derive by successive application of eq.\ (\ref{combi-eq2})
the identities (in abbreviated form)
\begin{equation}\label{combi-eq3}
\begin{split}
[j]! &= \Big(\id_{X^{\otimes\,j-2}}\otimes\Big[{2\atop 1}\Big]\Big)
  \circ\Big(\id_{X^{\otimes\,j-3}}\otimes\Big[{3\atop 1}\Big]\Big)
  \circ\dots\circ\Big[{j\atop 1}\Big]\,,\\
[j]! &= \Big(\Big[{2\atop 1}\Big]\otimes\id_{X^{\otimes\,j-2}}\Big)
   \circ\Big(\Big[{3\atop 2}\Big]\otimes\id_{X^{\otimes\,j-3}}\Big)
  \circ\dots\circ\Big[{j\atop j-1}\Big]
\end{split} 
\end{equation}
and
\begin{equation}\label{combi-eq4}
[j+k]! = ([j]!\otimes [k]!)\circ \Big[{j+k\atop j}\Big]\,.
\end{equation}
Similarly the equations which are dual to (\ref{combi-eq3}) and
(\ref{combi-eq4}) can be deduced. For the first equation in
(\ref{combi-eq3}) see also \cite{Ma4}.
We have provided all the definitions
and results which are necessary for the formulation and the proof of the
proposition on tensor Hopf algebras in braided categories extending the
corresponding theorem in \cite{Ma5}.
\abs
\begin{proposition}\label{3-4f-D11II}
The tensor algebra $T_{\C}(X)$ in ${\C}^{\uNN}$
of an object $X$ in a $\otimes$-factor braided abelian category $\C$
has a Hopf algebra structure which is given by
\begin{eqnarray}
\m_{n,m} &\cong& \id_{X^{\otimes\,n+m}}: X^{\otimes\,n}\otimes
X^{\otimes\,m}\to X^{\otimes\,n+m}\,,\nonumber\\
\eta_0 &=& \id_k,\ \eta_j=0\ {\it for}\ j\ne 0\,,\nonumber\\
\label{tens-alg1}
\Delta_{n,m} &\cong& \Big[{n+m\atop n}\Big\vert X;\lambda\Big]:
 X^{\otimes\,n+m}\to X^{\otimes\,n}\otimes X^{\otimes\,m}\,,\\
\varepsilon_0 &=& \id_k,\ \varepsilon_j=0\ {\it for}\ j\ne 0\,,\nonumber\\
S_n &=& (-1)^n\lambda^{{}^\binom{n}{2}}\,(\sigma_n^0)_{\C}(X)
: X^{\otimes\,n}\to X^{\otimes\,n}
\nonumber
\end{eqnarray}
where $\sigma_n^0=\left(\begin{smallmatrix}1 & \ldots & n\\ n & \ldots &1
\end{smallmatrix}\right)$.
Because of duality reasons another Hopf structure can be established on
$T_{\C}(X)$ according to
\begin{eqnarray}
{\overset\circ\m}_{n,m} &\cong&
\Big[{n\atop n+m}\Big\vert X;\lambda\Big]:
X^{\otimes\,n}\otimes X^{\otimes\,m}\to X^{\otimes\,n+m}\,,\nonumber\\
{\overset\circ\eta}_0 &=& \id_k,\
{\overset\circ\eta}_j=0\ {\it for}\ j\ne 0\,,\nonumber\\
\label{tens-alg2}
{\overset\circ\Delta}_{n,m} &\cong& \id_{X^{\otimes\,n+m}}:
X^{\otimes\,n+m}\to X^{\otimes\,n}\otimes X^{\otimes\,m}\,,\\
{\overset\circ\varepsilon}_0 &=& \id_k,\
{\overset\circ\varepsilon}_j=0\ {\it for}\ j\ne 0\,,\nonumber\\
{\overset\circ S}_n &=&
(-1)^n\lambda^{{}^\binom{n}{2}}\,(\sigma_n^0)_{\C}(X):
 X^{\otimes\,n}\to X^{\otimes\,n}\,.\nonumber
\end{eqnarray}
For distinction we denote by $T_{\C}(X)$ the Hopf algebra
defined through (\ref{tens-alg1}). The dual Hopf algebra corresponding
to (\ref{tens-alg2}) will be denoted by ${\overset\circ T}_{\C}(X)$.
\end{proposition}

\begin{proof}
Using the results on braided multinomials the proof of the proposition
is not difficult. To give an insight into the techniques we sketch 
the proof of the Hopf algebra structure (\ref{tens-alg1}).
Without problems one verifies the algebra structure of (\ref{tens-alg1}).
For the proof of the coalgebra properties one observes that Proposition
\ref{combi} implies
\begin{equation}
\Big(\Big[{m+n\atop m}\Big]\otimes\id_{X^{\otimes\,k}}\Big)\circ
\Big[{m+n+k\atop m+n}\Big]
=\Big(\id_{X^{\otimes\,m}}\otimes\Big[{n+k\atop n}\Big]\Big)\circ
\Big[{m+n+k\atop m}\Big]
\end{equation}
which yields the coassociativity. A similar investigation as in the proof
of Proposition \ref{combi} leads to the relation
\begin{equation}\label{comult-mult}
\Big[{n+m+p+q\atop n+p}\Big]=(\id_{X^{\otimes\,n}}\otimes
\Psi^{(\lambda)}_{X^{\otimes\,m},X^{\otimes\,p}}\otimes
\id_{X^{\otimes\,q}})\circ\Big(\Big[{n+m\atop n}\Big]\otimes
\Big[{p+q\atop p}\Big]\Big)\,.
\end{equation}
Then it can be shown directly that $\hat \Delta$ in (\ref{tens-alg1})
is an algebra morphism, that is
\begin{equation}
\hat\Delta\circ\hat\m = (\hat\m\otimes\hat\m)\circ
(\hat\id\otimes\hat\Psi^{(\lambda)}_{T,T}\otimes\hat\id)\circ
(\hat\Delta\otimes\hat\Delta)\,.
\end{equation}
Using (\ref{comult-mult}) and the relation
$S_n=-\Psi^{(\lambda)}_{X,X^{\otimes,n-1}}\circ(\id_X\otimes S_{n-1})$
we verify the Hopf equation
$\hat\m\circ(\hat\id\otimes \hat S)\circ\hat\Delta=
\hat\eta\circ\hat\varepsilon$ by induction on $n$, and analogously we
derive $\hat\m\circ(\hat S\otimes\hat\id)\circ\hat\Delta=
\hat\eta\circ\hat\varepsilon$.\end{proof}
\abs
Up to now braided multinomials are involved which
generally depend on an automorphism $\lambda$ of the unit object
in $\C$. However, the study of the differential
structure in Chapter \ref{sec-diff-calc} forces us to consider
$\lambda= -1$. Therefore we focus our consideration to $\lambda= -1$
without serious loss of generality. In what follows we will use the
notation
$\big[{{j_1\dots j_r}\atop j}\big\vert X\big] :=
 \big[{{j_1\dots j_r}\atop j}\big\vert X;(-1)\big]$ and
$\big[{j\atop {j_1\dots j_r}}\big\vert X\big] :=
 \big[{j\atop {j_1\dots j_r}}\big\vert X;(-1)\big]$ for the braided
multinomials, and $\hat\Psi :=\hat\Psi^{(-1)}$ for the braiding in
${\C}^{\uNN}$.
\abs
\begin{definition}\label{anti-sym}
Let $X$ be an arbitrary object in the category $\C$. Then the
antisymmetrizer $\hat A_X : T_{\C}(X)\to T_{\C}(X)$ is given by
$(\hat A_X)_n=[n\vert X]!: X^{\otimes\,n}\to X^{\otimes\,n}$.
\end{definition}
\abs
The antisymmetrizer $\hat A_X$ is compatible with the Hopf structures
of $T_{\C}(X)$ and ${\overset\circ T}_{\C}(X)$.
\abs
\begin{proposition}\label{3-1ZD11XII}
Let $X$ be an object in $\C$. Then $\hat A_X: T_{\C}(X)\to
{\overset\circ T}_{\C}(X)$ is a Hopf algebra morphism.
\end{proposition}

\begin{proof}
The proof is straightforward. We will only outline that $\hat A_X$ is
an algebra morphism. We use the relation dual to (\ref{combi-eq4}).
\begin{equation}
{\overset\circ\m}_{n,m}\circ(A_{X,n}\otimes A_{X,m})
=\Big[{n\atop n+m}\Big]\circ\big([n]!\otimes[m]!\big)=[n+m]!
\cong A_{X,n+m}\circ\m_{n,m}\,.
\end{equation}
\end{proof}
\abs
Proposition \ref{3-1ZD11XII} and Lemma \ref{8-4Z31} lead to the
definiton of the braided antisymmetric tensor algebra.
\abs
\begin{definition}\label{anti-tens-alg}
Let $X$ be an object in the category
$\cal C$. Then the antisymmetric tensor algebra of
$X$ is the Hopf algebra $T^\wedge_{\cal C}(X)$ in ${\cal C}^{\uNN}$
canonically induced by the decomposition
$\hat A_X=\im\hat A_x\circ\coim\hat A_x$.
\begin{equation}
\label{Anti-epi-mono}
T_{\cal C}(X)@>{\coim{\hat A_X}}>>T^\wedge_{\cal C}(X)
           @>{\im{\hat A_X}}>>{\overset\circ T}_{\cal C}(X)\,.
\end{equation}
\end{definition}
\abs
Definition \ref{anti-tens-alg} suggests to define the functors
$T_{\cal C}(\_\,),\ {\overset\circ T}_{\cal C}(\_\,),\
T^\wedge_{\cal C}(\_\,): \C\longrightarrow \text{Hopf-alg-}\,\C^{\uNN}$
on the objects according to (\ref{Anti-epi-mono}).
On morphisms $T_{\cal C}(\_\,)$ and
${\overset\circ T}_{\cal C}(\_\,)$ are given canonically, and
for any morphism $g:X\to Y$ in $\C$ the morphism
$T^\wedge_{\cal C}(g):T^\wedge_{\cal C}(X)\to T^\wedge_{\cal C}(Y)$ is
defined to be the unique morphism such that
$\im\hat A_Y\circ T^\wedge_{\cal C}(g) ={\overset\circ T}_{\cal C}(g)
\circ \im\hat A_X$. Then the following result can be derived which
will be needed in Chapter \ref{sec-diff-calc} for the construction of a
graded differential Hopf algebra out of
a first order bicovariant differential calculus.
\abs
\begin{lemma}\label{5-D11XVII}
Let $\cal C$ be a $\otimes$-factor braided abelian category.
Then the functors $T_{\cal C}(\_\,)$, $T^\wedge_{\cal C}(\_\,)$  and
${\overset\circ T}_{\cal C}(\_\,)$ map epi-/monomorphisms to
epi-/monomorphisms. Every morphism $g:X\to Y$ in $\C$ yields
$\big(T_{\cal C}(g),T^\wedge_{\cal C}(g),
{\overset\circ T}_{\cal C}(g)\big)$ which is a morphism of sequences of
the form (\ref{Anti-epi-mono}).
\end{lemma}

\begin{proof}
Obviously for any $g:X\to Y$ in $\C$ the triple
$\big(T_{\cal C}(g),T^\wedge_{\cal C}(g),
{\overset\circ T}_{\cal C}(g)\big)$ is a morphism of the corresponding
diagrams (\ref{Anti-epi-mono}). The assertion on the epi-/monomorphism
property is easily verified for $T_{\cal C}(\_\,)$ and
${\overset\circ T}_{\cal C}(\_\,)$. Then the commutativity of the diagrams
involving the morphisms $\big(T_{\cal C}(g),T^\wedge_{\cal C}(g),
{\overset\circ T}_{\cal C}(g)\big)$ yields the statement also for
$T^\wedge_{\cal C}(\_\,)$.\end{proof}
\abs
For a flat Hopf algebra $H$  with isomorphic antipode
in the $\otimes$-factor braided abelian category $\C$ we consider the
category $(\hhchh)^{\uNN}$ which is again
$\otimes$-factor braided abelian due to Proposition \ref{5-6-4ZD11}.
In analogy to Definition
\ref{anti-tens-alg} the functors $T_{\hhchh}(\_\,),\ 
{\overset\circ T}_{\hhchh}(\_\,),\ T_{\hhchh}^\wedge(\_\,):\hhchh\longrightarrow
(\hhchh)^{\uNN}$ can be defined. Then we apply Theorem \ref{bialg-proj} to
produce Hopf algebras in ${\C}^{\uNN}$.
\abs
\begin{definition}\label{ext-alg}
Given a flat Hopf algebra $H$ with bijective antipode in $\C$ and an
$H$-Hopf bimodule $X$. Then the tensor Hopf algebras
${\cal T}{}^H_{\C}(X)$, ${\overset\circ{\cal T}}{}^H_{\C}(X)$ and the
braided exterior tensor Hopf algebra of forms $X^{\wedge_H}$ over $X$ in
${\C}^{\uNN}$ are defined as
${\cal T}{}^H_{\C}(X):=G(T_{\hhchh}(X))$,
${\overset\circ{\cal T}}{}^H_{\C}(X):= G({\overset\circ T}_{\hhchh}(X))$
and $X^{\wedge_H}:=G(T_{\hhchh}^\wedge(X))$
respectively, where the corresponding functor $G$ of Theorem
\ref{bialg-proj} is used.
\end{definition}  
\abs

\begin{remark}
{\rm The Theorems \ref{yd-hopfbi} and
\ref{bialg-proj} and the results of the appendix may be used to 
prove that
\begin{equation}\label{smash}
X^{\wedge_H}\cong G\big(H\ltimes T^\wedge_{\DY({\C})^H_H}({}_HX)\big)
\end{equation}
are isomorphic Hopf algebra objects. In anticipation of Chapter
\ref{sec-diff-calc} suppose that $(X,\dif)$ is a braided first
order bicovariant differential calculus over $H$. Then
the object $(X^{\wedge_H},\dif^\wedge)$
is the braided exterior Hopf algebra of differential forms over $X$, and
(\ref{smash}) extends to an isomorphism of differential Hopf algebras
\begin{equation}\label{smash1}
(X^{\wedge_H},\dif^\wedge)\cong
\big(G\big(H\ltimes T^\wedge_{\DY({\C})^H_H}({}_HX)\big),
H\ltimes {}_H(\hat\dif)\big)\,.
\end{equation}
This particularly implies the results
(in the classical symmetric situation) of \cite{SZ} and dually of
\cite{Dr2}. There it was found that the higher order differential calculi
of \cite{Wor} and the quantum standard complexes are cross
products of the quantum group and the quantum enveloping algebra
respectively with a certain antisymmetric tensor algebra of invariant
vector fields.}
\end{remark}
\abs
Under the conditions of Definition \ref{ext-alg} we know from Proposition
\ref{12-8IIIZD34} that $(H,X)$ is a bialgebra in
${\C}^{\mathbf 2}$. By construction of the tensor algebras
${\cal T}{}^H_{\C}(X)$ and ${\overset\circ{\cal T}}{}^H_{\C}(X)$
we derive from Theorem \ref{bialg-proj}
and the corresponding definition in \eqref{bialg-trans2} that
the multiplication of ${\cal T}{}^H_{\C}(X)$ is given by
$\m^{{\cal T}{}^H_{\C}(X)}=
\lambda^H_{{\cal T}{}^H_{\C}(X),{\cal T}{}^H_{\C}(X)}$ and the
comultiplication of ${\overset\circ{\cal T}}{}^H_{\C}(X)$ equals
$\Delta^{{\overset\circ{\cal T}}{}^H_{\C}(X)}=
{\rho^H}{}_{{\overset\circ{\cal T}}{}^H_{\C}(X),
{\overset\circ{\cal T}}{}^H_{\C}(X)}$. This already indicates the
universality of the tensor algebras which will be demonstrated in
the following theorems and in Theorem \ref{8-10-D40}.
\abs
\begin{proposition}\label{21-8IXD34}
Suppose that $H$ is a flat Hopf algebra with isomorphic antipode in
$\C$ and $X$ is an $H$-Hopf bimodule.

If $\hat Y$ be an algebra resp.\ a bi-/Hopf algebra in ${\C}^{\uNN}$ and
$\hat f=(f_0,f_1):(A,X)\to (Y_0,Y_1)$ is an algebra resp.\ bialgebra
morphism in ${\cal C}^{\{0,1\}}$ then there exitsts a unique algebra
resp.\ bialgebra morphism
$\Hat{\Hat f}: {\cal T}{}^A_{\cal C}(X)\to \hat Y$ such that
$(\Hat{\Hat f})_0=f_0$ and $(\Hat{\Hat f})_1=f_1$.

A dually symmetric result holds for a coalgebra resp.\ bialgebra morphism
$\overset\circ g=(g_0,g_1): (Y_0,Y_1)\to (H,X)$ which can be extended
uniquely to a coalgebra resp.\ bialgebra morphism
$\hat{\overset\circ g}: \hat Y\to {\overset\circ{\cal T}}{}^H_{\cal C}(X)$
with $(\hat{\overset\circ g})_0=g_0$ and $(\hat{\overset\circ g})_1=g_1$.
\end{proposition}

\begin{proof}
Since $(f_0,f_1):(H,X)\to(Y_0,Y_1)$ is an algebra morphism by assumption 
and the multiplication of ${\cal T}{}^H_{\C}(X)$ is given by
the universal tensor product morphism $\m^{{\cal T}{}^H_{\C}(X)}=
\lambda^H_{{\cal T}{}^H_{\C}(X),{\cal T}{}^H_{\C}(X)}$
one immediately verifies successively that there exists a unique morphism
$f_n:X^{\otimes_H\,n}\to Y_n$ for every $n\in\NN_0$ such that
\begin{equation*}
\m^Y_{\underbrace{1,\ldots,1}_{n}}\circ
(\underbrace{f_1\otimes\ldots\otimes f_1}_{n})=f_{n}\circ
\m^{{\cal T}{}^H_{\cal C}(X)}_{\underbrace{1,\ldots,1}_{n}}\,.
\end{equation*}
Then the identity
$f_{m+n}\circ\m^{{\cal T}{}^H_{\cal C}(X)}_{m,n}\circ
(\m^{{\cal T}{}^H_{\cal C}(X)}_{1,\dots , 1}
\otimes
\m^{{\cal T}{}^H_{\cal C}(X)}_{1,\dots , 1})=
\m^Y_{m,n}\circ(f_m\otimes f_n)\circ
(\m^{{\cal T}{}^H_{\cal C}(X)}_{1,\dots , 1}
\otimes
\m^{{\cal T}{}^H_{\cal C}(X)}_{1,\dots , 1})$
can be derived which implies
$f_{m+n}\circ\m^{{\cal T}{}^H_{\cal C}(X)}_{m,n}
 =\m^Y_{m,n}\circ(f_m\otimes f_n)$
because the morphism
$\m^{{\cal T}{}^H_{\cal C}(X)}_{1,\dots , 1}$ is epimorphic.
Thus $\Hat{\Hat f}:=(f_n):{\cal T}{}^H_{\cal C}(X)\to\hat Y$ is an algebra
morphism.

If $\hat Y$ is a bi-/Hopf algebra and $(f_0,f_1):(H,X)\to (Y_0,Y_1)$ is a
bialgebra morphism, it follows by induction on $n>1$ and by use of
the fact that $\Hat{\Hat f}$ is an algebra morphism, that
\begin{equation*}
\Delta^{\hat Y}_{t,n-t}\circ f_n\circ\m^{{\cal T}{}^H_{\cal C}(X)}_{n-1,1}
= (f_t\otimes f_{n-t})\circ\Delta^{{\cal T}{}^H_{\cal C}(X)}_{t,n-t}\circ 
\m^{{\cal T}{}^H_{\cal C}(X)}_{n-1,1}\,.
\end{equation*}
Hence $\Hat{\Hat f}$ is a bialgebra morphism since
$\m^{{\cal T}{}^H_{\cal C}(X)}_{n-1,1}$ is epimorphic.
\end{proof}
\abs
\begin{remark} {\rm Note that under the assumption of right exactness of
$X\otimes(\_\,)$ for
every $X\in\Ob(\C)$ for any algebra $A$ in $\C$ the category ${}_A\C_A$ of
$A$-bimodules with the tensor product $\otimes_A$ over $A$ is monoidal.
And one can define a free tensor algebra ${\cal T}^A_{\C}(X)$ for every
$A$-bimodule $X$ with multiplication given by the canonical morphisms
into the tensor product over $A$.
This tensor algebra is characterized by the same universal property
as described in Proposition \ref{21-8IXD34}.}
\end{remark}
\abs
Consider as a special example of Proposition \ref{21-8IXD34} the Hopf
algebra $\hat Y := {\overset\circ{\cal T}}{}^H_{\C}(X)$.
Since $(\id,\id):(H,X)\to (H,X)$ is a bialgebra morphism we obtain with the
help of Proposition \ref{21-8IXD34} the following corollary.
\abs
\begin{corollary}\label{22-8XIIID34}
The unique higher order bialgebra
extension of the identity morphism ${\hat\id:(H,X)\to (H,X)}$ is given by
$\Hat{\Hat\id}=\hat A^H_X:{\cal T}{}^H_{\C}(X)\to
{\overset\circ{\cal T}}{}^H_{\C}(X)$ where $\hat A^H_X$ is the
antisymmetrizer of $X$ in the category $(\hhchh)^{\uNN}$.\endproof
\end{corollary}
\abs
The main theorem in Chapter \ref{grade-alg-com} states the universal
property of the exterior Hopf algebra $X^{\wedge_H}$
generated algebraically by its 0th and 1st component.
\abs
\begin{theorem}\label{23-8XIIID34}
Suppose that $H$ is a flat Hopf algebra with isomorphic antipode
and $X$ is an $H$-Hopf bimodule in $\C$. Then the exterior Hopf algebra
$X^{\wedge_H}$ is generated by $H$ and $X$ as an algebra, that is
$\Im(\m_{1,n}^{X^{\wedge_H}})=X^{\wedge_H}_{n+1}$ for any
$n\in\NN_0$. Let $\hat Y$ be a bi-/Hopf algebra in ${\C}^{\uNN}$.
Suppose that $\hat Y$ is generated by $Y_0$ and $Y_1$ as an algebra.
If $(g_0,g_1):(Y_0,Y_1)\to (H,X)$ is a bialgebra morphism in
${\C}^{\mathbf 2}$, then there exists a unique bialgebra morphism
$g^{\wedge_H}:\hat Y\to X^{\wedge_H}$ in ${\C}^{\uNN}$ such
that $(g^{\wedge_H})_0=g_0$ and $(g^{\wedge_H})_1=g_1$. In this case
it holds $\hat{\overset\circ g}=\im \hat A^H_X\circ g^{\wedge_H}$.
\end{theorem}

\begin{proof}
On $X^{\wedge_H}$ the multiplication is given by $\m^{X^{\wedge_H}}=
\m_{\Coim \hat A}\circ\lambda^H_{\Coim\hat A,\Coim\hat A}$ according to
the construction in Theorem \ref{bialg-proj}. Then it follows
\begin{equation}\label{mult-ext-hopf}
 \m_{X^{\wedge_H}}\circ(\coim\hat A\otimes\coim\hat A)
 =\coim\hat A\circ\lambda^H_{T_\hhchh(X),T_\hhchh(X)}
\end{equation}
and since the right hand side of (\ref{mult-ext-hopf}) is epimorphic,
one deduces that $X^{\wedge_H}$ is generated algebraically by $H$ and $X$.
Now let us suppose that there exists a bialgebra morphism
$(g_0,g_1):(Y_0,Y_1)\to(H,X)$. Proposition \ref{21-8IXD34} and
Corollary \ref{22-8XIIID34} then tell us that there is a bialgebra morphism
$\hat{\overset\circ g}: \hat Y\to{\overset\circ{\cal T}}{}^H_{\C}(X)$.
It holds $\coker(\hat A)\circ\hat{\overset\circ g}=0$ which can be seen by
induction as follows. For $n=0,1$ one verifies immediately that
$(\coker\hat A)_n\circ(\hat{\overset\circ g})_n=0$ since $(\hat A)_0=\id$
and $(\hat A)_1=\id$. Therefore $g_0=(\hat{\overset\circ g})_0=
(\im\hat A)_{0}\circ g^\wedge_0=:g^\wedge_0$ and $g_1=
(\hat{\overset\circ g})_1=(\im\hat A)_{1}\circ g^\wedge_1=:g^\wedge_1$.
If $n>1$ one derives the relation $(\coker\hat A)_n\circ
(\hat{\overset\circ g})_n\circ\m^{\hat Y}_{n-1,1}=0$ because 
$(\hat{\overset\circ g})$ is a bialgebra morphism,
$(\hat{\overset\circ g})_m=(\im\hat A)_{m}\circ g^\wedge_m$ for every $m<n$,
and $(\im\hat A)$ is a bialgebra morphism. Since $\m^{\hat Y}_{n-1,1}$ is
epimorphic by assumption, we obtain $(\hat{\overset\circ g})_n=
(\im\hat A)_{n}\circ g^\wedge_n$ for some unique $g^\wedge_n$.\end{proof}
\abs
\begin{remark}{\rm In the case of the antisymmetric and exterior Hopf
algebras it is not sufficient in general to restrict one's considerations
to the grades less than 3 and to argue by multiplicative continuation as
will be demonstrated in the following.
Suppose that the category $\C$ is $\otimes$-exact braided
abelian and $X\otimes\id_\C$ is right exact. Consider the ideal
$\big({\ker}[2\vert X;\lambda]!\big)$ in $T_{\hhchh}(X)$
generated by ${\ker}[2\vert X;\lambda]!$ according to Proposition
\ref{2-3ZD30}. Again $\lambda$ is an automorphism of the unit object.
It is not difficult to verify that $\big({\ker}[2\vert X;\lambda]!\big)$
is a Hopf ideal and a subobject of $\ker\Hat A_X$ because
$[n+m\vert X;\lambda]!=\big([n-2\vert X;\lambda]!\otimes[2\vert X;\lambda]!
\otimes [m\vert X;\lambda]!\big)\circ \big[{n+m\atop (n-2,2,m)}\vert X;
\lambda\big]$, and dually by Proposition \ref{combi} and
eq.\ (\ref{combi-eq4}). But in general the ideal $\ker\Hat A_X$ does not
coincide with $\big({\ker}[2\vert X;\lambda]!\big)$. For example the
tensor algebra of polynomials in one variable $\CC\langle x\rangle$ can be
equipped with the braiding
$\Psi(x\otimes x)=q(x\otimes x)$, $q\in\CC^{\,*}$. This is the braided
line as described in \cite{Koo,Ma4}. The braided integers $[n\vert X;q]$
are just the $q$-numbers $[n\vert X;q]=[n]_q:=1+q+\dots +q^{n-1}$
and the braided factorial $[n\vert X;q]!$ equals $[n]_q!$ which is the
usual $q$-factorial. For $q$ a primitive root of unity of order
$n>2$ one obtains $({\ker}[2\vert X;q]!)=0$ which does not
coincide with the ideal $\ker\Hat A_X=(x^{\otimes n})$.}
\end{remark}
\abs
\abs
\section{Differential Calculi and Exterior Differential Hopf Algebras}
\label{sec-diff-calc}
In Chapter \ref{sec-diff-calc} we define differential calculi in
$\otimes$-factor abelian
categories following the ideas of \cite{Mal,Man,Wor}. Special
interest is given to the generalization of bicovariant differential
calculi in $\otimes$-factor braided abelian categories. We prove
the existence of a higher order differential Hopf algebra calculus which
extends a given first order bicovariant differential calculus in a natural
way. This is the braided exterior Hopf algebra of differential forms over
the first order bicovariant differential calculus. The corresponding
results of \cite{Wor} are therefore generalized to $\otimes$-factor braided
abelian categories.
\abs
\subsection*{Differential Calculi}

\begin{definition}\label{diff-calc}
A complex $(\hat X,\hat \dif)$ in $\C^{\I_c}$ is called a
differential calculus if $(\hat X,\hat \dif)$ is an algebra in
$\C^{\I_c}$ and $X_0\langle \dif_n\rangle = X_{n+1}$ for all $n\in\I$.
The category of differential calculi as a full subcategory of
${{\rm Alg}\text{-}\C^{\I_c}}$ will be denoted by
${{\rm Diff}\text{-}\,\C^{\I_c}}$.
\end{definition}
\abs
Sometimes we adopt the notation of \cite{Wor} and we speak of first order
differential calculi if $\I_c={\mathbf 2}_c$ and of higher order differential
calculi if $\I_c=\uNN_c$. From Definition \ref{diff-calc} we derive
some implications in analogy to the results of \cite{Wor} which
we will collect in the next propositions.
\abs
\begin{proposition}\label{1-2-1ZD34}
The following statements for an algebra $(\hat X,\hat \dif)$ in
$\C^{\I_c}$ are equivalent.
\begin{enumerate}
\item $(\hat X,\hat\dif)$ is a differential calculus.
\item For all $n\in \I$ it holds $X_0\langle \dif_n\rangle =X_{n+1}$.
\item For all $n\in \I$ it holds $\langle \dif_n\rangle X_0 =X_{n+1}$.
\item For all $n\in \I$ it holds $X_0\langle \dif_n\rangle X_0 =X_{n+1}$.
\item For all $n\in \I$ it holds $X_0\langle\m^{(n-1)}_{1,\ldots,1}\circ
      (\dif\otimes\ldots\otimes\dif)\rangle=X_{n+1}$.
\end{enumerate}
The differential of a differential calculus $(\hat X,\hat \dif)$ is
determined uniquely by its 0th component $d_0$. 
\end{proposition}

\begin{proof}
The proof of the equivalences is rather straightforward. It can be
performed as in the symmetric category of vector spaces \cite{Wor}
because no braidings are involved. To give an idea of the techniques used
in the proof we will demonstrate the equivalence of (2) and (3).
We only have to prove that (2) implies (3). Because of symmetry reasons
one can argue similarly that (2) follows from (3).

Suppose that $\m_{0,l}\circ(\id_0\otimes\dif_{l-1})$ is epimorphic. Then
it follows by induction, using the identity
\begin{equation}
\begin{gathered}
\m_{k,l+m+1}\circ\big(\id_k\otimes\dif_{l+m}\circ\m_{m,l}\circ(\id_m\otimes
\dif_{l-1})\big)\\
=\m_{k+m+1,l}\circ(\id_{k+m+1}\circ\dif_{l-1})\circ
\big(\m_{k,m+1}\circ(\id_k\otimes\dif_m)\otimes\id_{l-1}\big)
\end{gathered}
\end{equation}
that $\m_{k,l}\circ(\id_k\otimes\dif_{l-1})$ is epimorphic.
$(\hat X,\hat\dif)$ is an algebra in $\C^{\I_c}$ which implies the
relation
\begin{equation}\label{diff-mult}
\dif_n\circ\m_{k,l}=\m_{k+1,l}\circ(\dif_k\otimes\id_l)+
(-1)^k \m_{k,l+1}\circ(\id_k\otimes\dif_l)\,.
\end{equation}
Suppose that $f\circ\m_{k+1,0}\circ(\dif_k\otimes\id_0)=0$
for some morphism $f$. Then $f\circ\dif_n\circ\m_{k,0}\circ(\id_k\otimes
\eta_0)=(-1)^k f_{k+1}\circ(\id_k\otimes\dif_0\circ\eta_0)=0$ because
(\ref{diff-mult}) and $\dif_0\circ\eta_0=0$ hold. Therefore
$f\circ\m_{k,1}\circ
(\id_k\otimes\dif_0)=0$ which yields $f=0$. Hence $\m_{k+1,0}\circ
(\dif_k\otimes\id_0)$ is epimorphic.\end{proof}
\abs
\begin{proposition}\label{34-D34}
Let $(\hat X,\hat \dif_X)$ be a differential calculus in $\C^{\I_c}$
and $(\hat Y,\hat \dif_Y)$ be an algebra in $\C^{\I_c}$.
Then any algebra morphism $\hat f:(\hat X,\hat \dif_X)\to
(\hat Y,\hat \dif_Y)$ is uniquely determined by its $0$th component $f_0$.
\end{proposition}

\begin{proof}
From $f_n\circ\m^X_{0,n}\circ(\id_0\otimes\dif_{X,n-1})=
\m^Y_{0,n}\circ(\id_0\otimes\dif_{Y,n-1})\circ(f_0\otimes f_{n-1})$
the proof can be concluded since
$\m^X_{0,n}\circ(\id_0\otimes\dif_{X,n-1})$ is epimorphic.\end{proof}
\abs
If we impose further conditions on $\C$ we can prove the converse as will
be seen in the following proposition.
\abs
\begin{proposition}\label{diff-calc-ext2}
Suppose that $\C$ is $\otimes$-exact braided abelian.
Let $(\hat X,\hat\dif_X)$ be an algebra in $\C^{\I_c}$, and let
the functor $\hat X\otimes\id_{\C^\I}$ be right exact.
Suppose that any algebra morphism
$\hat f:(\hat X,\hat\dif_X)\to(\hat Y,\hat\dif_Y)$ in $\C^{\I_c}$
is uniquely determined by its $0$th component $f_0$.
Then $(\hat X,\hat\dif_X)$ is a differential calculus in $\C^{\I_c}$.
\end{proposition}

\begin{proof}
Define the object $\hat Y$ in ${\C}^\I$ by
$Y_{n+1}=(\hat X\otimes\hat X\otimes\hat X)_n$ for all
$n\in\I$. For $m<0$ we set $X_m=0$ and $\dif_{X,m}=0$
in the following. We consider the morphism $\hat M_\dif:=\hat\m_{\hat X}
\circ(\hat\id_{\hat X}\otimes\hat\m_{\hat X})\circ(\hat\id_{\hat X}\otimes
\hat\dif_X\otimes\hat\id_{\hat X}):\hat Y\to\hat X$.
It is the morphism which multiplies $\hat X$ from the left and from the
right to $\hat\dif_X(\hat X)$. It holds $M_{\dif,0}=0$.
Because of Proposition \ref{2-3ZD30} (in ${\C}^\I$) we know that
$\hat X\langle\hat\dif_X\rangle\hat X=\im\hat M_\dif$ is an ideal in
$\hat X$. With the help of Proposition \ref{17-32} (in ${\C}^\I$)
one observes that the morphism $\coker\hat M_\dif:{\hat X\to\hat P}$ is
an algebra morphism. In particular $X_0\cong P_0$. For $n>0$ the equation
$\coker(\hat M_\dif)_{n+1}\circ\dif_{X,n}\circ\im (\hat M_\dif)_n=0$
can be derived by use of the relation
\begin{equation}
\begin{split}
&\dif_{X,n}\circ\m^X_{k,l+m+1}\circ\big(\id_{X_k}\otimes
\m^X_{l+1,m}\circ(\dif_{X,l}\otimes\id_{X_m})\big)\\
&=\big((\hat M_\dif)_{k+1,l,m}\circ(\dif_{X,k}\otimes
\id_{X_l\otimes X_m})+(-1)^{k+l+1}(\hat M_\dif)_{k,l,m+1}\circ
(\id_{X_k\otimes X_l}\otimes\dif_{X,m})\big)
\end{split}
\end{equation}
where $k+l+m=n-1$. For $n=0$ it holds
$\coker(\hat M_\dif)_1\circ\dif_{X,0}=
\coker(\hat M_\dif)_1\circ(\hat M_\dif)_1\circ(\eta_X\otimes\id_{X_0}
\otimes\eta_X)=0$. Hence there are unique morphisms
$\dif_{P,n}:P_n\to P_{n+1}$ such that
\begin{equation}\label{coker-dif}
\dif_{P,n}\circ\coker(\hat M_\dif)_n=\coker(\hat M_\dif)_{n+1}
\circ\dif_{X,n}\,,
\end{equation}
moreover $\dif_{P,0}=0$. From (\ref{coker-dif}) it follows that
$\coker(\hat M_\dif):(\hat X,\hat\dif_X)\to(\hat P,\hat\dif_P)$ is an
algebra morphism in $\C^{\I_c}$. Since $\coker(M_\dif)_0 =\id$
we conclude that $\coker(\hat M_\dif)_n=0$ for all $n>0$,
which indeed yields an algebra morphism in $\C^{\I_c}$ and is therefore
unique by assumption.\end{proof}
\abs
In the case of first order differential calculi the result of Proposition
\ref{diff-calc-ext2} can be derived under the weaker condition that
$X_0\otimes\id_\C$ is right exact. Before we will give the definition of
left, right and bi-covariant differential calculi in $\otimes$-factor braided
abelian categories, we state a corollary of the previous results
concerning the extendibility of algebra morphisms to differential algebra
morphisms.

\begin{corollary}\label{proof-2-D55}
Suppose that $(\hat X,\hat\dif_X)$ is a differential calculus and
$(\hat Y,\hat\dif_Y)$ is an algebra in $\C^{\I_c}$. If
$\hat f:\hat X\to\hat Y$ is an algebra morphism in ${\C}^\I$ such that
$f_1\circ\dif_{X,0}=\dif_{Y,0}\circ f_0$ then $\hat f: (\hat X,\hat\dif_X)
\to (\hat Y,\hat\dif_Y)$ is an algebra morphism in $\C^{\I_c}$.\endproof
\end{corollary}
\abs
We will introduce the maximal differential calculus of a given
differential algebra.  It will be used in the last section of Chapter
\ref{sec-diff-calc} for the construction of the braided exterior Hopf algebra
of differential forms. Guided by Proposition \ref{1-2-1ZD34} we give the
following definition.
\abs
\begin{definition}
Let $(\hat A,\hat\m,\hat\eta,\hat\dif)$ be an algebra in $\C^{\I_c}$.
The maximal
differential calculus $\hat A^{\rm Diff}$ is the graded sub-object
of $\hat A$ given by the inclusion
$\hat\i:\hat A^{\rm Diff}\hookrightarrow\hat A$
where\/ $\i_0:=\id_{A_0}$,\/
$\i_1:=A_0\langle\dif_0\rangle$ and\/
$\i_n:=\im(\m^{(n-1)}_{1\ldots,1}\circ \i_1^{\otimes,n})$ for $n>1$.
\end{definition}
\abs
The nomenclature is justified by the following proposition which can
be proven easily with the help of the previous results and techniques
already provided in the article.
\abs
\begin{proposition}\label{max-diff}
On $\hat A^{\rm Diff}$ a differential $\hat\dif^{\rm Diff}$ can be defined
such that $\hat\i:(\hat A^{\rm Diff},\hat\dif^{\rm Diff})\hookrightarrow
(\hat A,\hat\dif)$ is an algebra morphism in $\C^{\I_c}$ and
$(\hat A^{\rm Diff},\hat\dif^{\rm Diff})$ is a differential calculus.
If $(\hat A',\hat\dif')$ is a differential calculus and
$\hat \j: (\hat A',\hat\dif')\hookrightarrow (\hat A ,\hat\dif)$ is an
algebra embedding in $\C^{\I_c}$ then there exists a unique
embedding\/ $\hat\i':(\hat A',\hat\dif')\hookrightarrow
(\hat A^{\rm Diff},\hat\dif^{\rm Diff})$ such that\/
$\hat\j=\hat\i\circ\hat\i'$. If in addition $A'_0=A_0$
then $(\hat A',\hat\dif')\cong
(\hat A^{\rm Diff},\hat\dif^{\rm Diff})$.\endproof
\end{proposition}
\abs
Now let us mimick Woronowicz's approach and define $\hat B$-covariant
differential calculi for graded bialgebras $\hat B$ in $\C^\I$.
\abs
\begin{definition}\label{35-D35}
Let $\hat B$ be a bialgebra in the category ${\C}^\I$ and
$(\hat X,\hat\dif)$ be an object in $\C^{\I_c}$. Suppose that
$\hat B\otimes\id_{\C^\I}$ is left exact.
Then $(\hat X,\hat\dif)$ is called $\hat B$-left covariant,
$\hat{B}$-right covariant  or
$\hat{B}$-bicovariant differential
calculus if it is a differential calculus in the category of
$\hat B$-comodules ${}^{\hat{B}}(\C^{\I_c})$,
$(\C^{\I_c})^{\hat{B}}$ and
${}^{\hat{B}}(\C^{\I_c})^{\hat{B}}$ respectively.
The corresponding categories are denoted by
${\rm Diff}\text{-}{}^{\hat{B}}(\C^{\I_c})$,
${\rm Diff}\text{-}(\C^{\I_c})^{\hat{B}}$ and
${\rm Diff}\text{-}{}^{\hat{B}}(\C^{\I_c})^{\hat{B}}$
respectively.
\end{definition}
\abs
In the remainder of this chapter we are especially interested in
$H$-bicovariant differential calculi of a bi- or Hopf algebra $H$ in
$\C$. We call them bicovariant differential calculi {\it over} $H$ if their
0th component coincides with $H$. We also consider differential calculi which
are bi- or Hopf algebras in $\C^{\I_c}$. They are called differential
bi- or Hopf algebra calculi henceforth.
\abs
\begin{proposition}\label{1-2-1ZD35}
Given a differential bialgebra calculus $(\hat X,\hat\dif)$
in $\C^{\I_c}$. Suppose that $X_0\otimes\id_\C$ is left exact.
Then $(\hat X,\hat\dif)$ is a bicovariant
differential calculus over $X_0$. Conversely, if\/ $\I={\mathbf 2}$ and
if $(X_0,X_1,\dif)$ is a bicovariant differential calculus over $X_0$
then it is a differential bialgebra calculus. If $X_0$ is a Hopf
algebra then $(X_0,X_1,\dif)$ is a differential Hopf algebra calculus.
\end{proposition}

\begin{proof}
We use the results of Proposition \ref{5-6-4ZD11}, especially
(\ref{cat-grade-mod1}). From Proposition \ref{12-8IIIZD34} one derives
that $X_0$ is a bialgebra and $X_n$ are $X_0$-Hopf bimodules for all
$n\in\I$. This means that $(\hat X,\hat\dif)$ is an algebra which is at
the same time $X_0$-bicomodule in ${\C}^\I$. Proposition
\ref{21-22-37--6-1ZD11} ensures that the image condition of the
differential calculus $(\hat X,\hat\dif)$ in ${\C}^\I$ can be taken
over to ${}^{X_0}(\C^{\I_c})^{X_0}$. Hence $(\hat X,\hat\dif)$ is an
$X_0$-bicovariant differential calculus. Conversely, suppose that
$(X_0,X_1,\dif)$ is a bicovariant differential calculus over $X_0$. Then
$(X_0,X_1,\dif)$ is a differential calculus and an $X_0$-bicomodule
algeba in ${\C}^{{\mathbf 2}_c}$. In particular $X_0$ is a bialgebra and
$X_1$ is an $X_0$-Hopf bimodule. As a consequence of Remark
\ref{mod-alg-equiv}, $(X_0,X_1,\dif)$ is a differential bialgebra calculus
in ${\C}^{{\mathbf 2}_c}$.\end{proof}
\abs
Hence Definition \ref{35-D35} generalizes the notations of
\cite{Wor} because we know from Proposition \ref{1-2-1ZD35} that a
differential bialgebra calculus $(\hat X,\hat\dif)$ in $\C^{\I_c}$
is a braided bicovariant differential calculus over $X_0$. In the case of
$\I={\mathbf 2}$ this is just a braided version of the definition
of bicovariant differential calculi over the bialgebra $X_0$ given in
\cite{Wor}.
\abs
\begin{proposition}
\label{Bialg-ext}
Suppose that $\hat A$ is a bialgebra in $\C^\I$ and
$(\hat A,\hat\dif)$ is an algebra in $\C^{\I_c}$ which is generated
by its 0th and 1st component. 
If the identity
\begin{equation}
\label{comul-diff}
\Delta_{n+1}\circ\dif_n=\dif^{\otimes}_n\circ\Delta_n\;:\;
A_n\to(\hat A\otimes\hat A)_n
\end{equation}
holds for $n=0,1$ then it is true for any $n\in\I$ and therefore
$\hat A$ is bialgebra in ${\cal C}^{\I_c}$.
If  in addition $(A_0,A_1,\dif_0)$ is a bialgebra in
${\cal C}^{{\mathbf 2}_c}$ and $\hat A^{\rm Diff}\hookrightarrow\hat A$
is a sub-bialgebra in $\C^\I$ then $(\subdif{\hat A},\subdif{\hat\dif})$ is
bialgebra in $\C^{\I_c}$. In particular it is a bicovariant differential
calculus over $A_0$ if $A_0\otimes\id_\C$ is left exact.
\end{proposition}

\begin{proof}
To prove the first part we only have to consider $\I=\uNN$. Suppose that
(\ref{comul-diff}) holds
for $n\ge 1$. Then
\begin{equation*}
\begin{split}
(\Delta_{n+2}\circ\dif_{n+1})_{k,l}\circ\m_{1,n}
&=(\Delta_{n+2}\circ\dif_{n+1}\circ\m_{1,n})_{k,l;1,n}\\
&=(\m_{n+2}^\otimes\circ(\Delta\otimes\Delta)_{n+2}\circ
 \dif^\otimes_{n+1})_{k,l;1,n}\\
&=(\m_{n+2}^\otimes\circ\dif^{\otimes^2}_{n+1}\circ
  (\Delta\otimes\Delta)_{n+1})_{k,l;1,n}\\
&=(\dif^\otimes_{n+1}\circ\Delta_{n+1}\circ\m_{n+1})_{k,l;1,n}
 =(\dif^\otimes_{n+1}\circ\Delta_{n+1})_{k,l}\circ\m_{1,n}\,.
\end{split}
\end{equation*}
In the third equation use has been made of the induction hypothesis
since $1,n < n+1$. Therefore (\ref{comul-diff}) holds for $n+1$, because
$\m_{1,n}$ is supposed to be an epimorphism.

For the proof of the second part we note that
\begin{equation*}
\begin{split}
\Delta_2\circ\dif_1\circ\m_{0,1}\circ(\id_A\otimes\dif_0)
&=\Delta_2\circ\m_{1,1}\circ(\dif_0\otimes\dif_0)\\
&=\m^{\otimes}_{1,1}\circ
        (\Delta_1\circ\dif_0\otimes\Delta_1\circ\dif_0)\\
&=\m^{\otimes}_{1,1}\circ
     (\dif^{\otimes}_0\otimes\dif^\otimes_0)\circ(\Delta_0\otimes\Delta_0)\\
&=\dif^\otimes\circ\m^\otimes_{0,1}\circ
      (\Delta_0\otimes\Delta_1\circ\dif_0)\\
&=\dif^\otimes\circ\Delta_1\circ\m_{0,1}\circ(\id_A\otimes\dif_0)\,,
\end{split}
\end{equation*}
from which follows $\subdif{\Delta_2}\,\circ\,\subdif{\dif_1}=
\subdif{\dif_1}{}^\otimes\,\circ\,\subdif{\Delta_1}$. Since
$(A_0,A_1,\dif_0,\Delta_0,\Delta_1)$ is a bialgebra in
${\cal C}^{{\mathbf 2}_c}$
the analogous relation will be decuced for the 0th component.
Then we can apply the first part of the proposition to finish the
proof.\end{proof}
\abs
The following proposition is important for the subsequent investigation as
well as for the construction of the braided exterior Hopf algebra of
differential forms in the next section.
\abs
\begin{proposition}\label{22-D22} Let $\I={\mathbf 2}$ or $\I=\uNN$ and
$\C$ be a $\otimes$-factor braided monoidal
category. Suppose that $\hat A$ is an algebra in $\C^\I$ and
$x:\E\to A_1$ is a morphism in $\cal C$ such that
$\m_{2,n}\circ\big(\m_{1,1}\circ(x\otimes x)\otimes\id_n)=
\m_{n,2}\circ\big(\id_n\otimes\m_{1,1}\circ(x\otimes x))$ for all $n\in\I$
if $\I=\uNN$. Then the morphisms 
\begin{equation}
[x,\cdot]_n:=
\mu_{1,n}\circ(x\otimes\id_{A_n})-
(-1)^n\mu_{n,1}\circ(\id_{A_n}\otimes x):A_n\to A_{n+1}
\end{equation}
for any $n,n+1\in\I$, turn $(\hat A,[x,\cdot]^\wedge)$ into an algebra in
${\cal C}^{\I_c}$.

If $(\hat A,\hat\Delta,\hat\varepsilon)$ is a bialgebra in
${\cal C}^\I$, $x$ is bi-invariant, i.\ e.\ $\Delta_{0,1}\circ x
=\eta\otimes x$ and $\Delta_{1,0}\circ x=x\otimes\eta$, and if
the maximal differential calculus $\subdif{\hat A}$ of
$(\hat A,[x,\cdot]^\wedge)$ is a sub-bialgebra of $\hat A$ in $\C^\I$, then
$(\subdif{\hat A},\subdif{([x,\cdot]^\wedge)})$ is a bialgebra in
${\cal C}^{\I_c}$.
\end{proposition}

\begin{proof}
The following identities prove the differential algebra property for
$(\hat A,[x,\cdot]^\wedge)$.
\begin{equation*}
\begin{split}
[x,\cdot]_{m+n}\circ\m_{m,n}
&=\m_{1,m+n}\circ(x\otimes\m_{m,n})+
  (-1)^{m+n+1}\m_{m+n,1}\circ(\m_{m,n}\otimes x)\\
&=\m_{1,m,n}\circ(x\otimes\id_{A_m\otimes A_n})+
  (-1)^{m+n+1}\m_{m,n,1}\circ(\id_{A_m\otimes A_n}\otimes x)\\
&=\m_{m+1,n}\circ\Big(\big(
   \m_{1,m}\circ(x\otimes\id_{A_m})
   -(-1)^m\m_{m,1}\circ(\id_{A_m}\otimes x)\big)\otimes\id_{A_n}\Big)\\
&\quad +(-1)^m\m_{m,n+1}\circ\Big(\id_{A_m}\otimes\big(
   \m_{1,n}\circ(x\otimes\id_{A_n})
   -(-1)^n\m_{n,1}\circ(\id_{A_n}\otimes x)\big)\Big)\\
&=\m_{m+1,n}\circ([x,\cdot]_m\otimes\id_{A_n})
 +(-1)^m\m_{m,n+1}\circ(\id_{A_m}\otimes[x,\cdot]_n)
\end{split}
\end{equation*}
and
\begin{equation*}
\begin{split}
[x,\cdot]_{n+1}\circ[x,\cdot]_n &=
 \m_{1,1,n}\circ(x\otimes x\otimes\id_{A_n})+
 (-1)^{n+1}\m_{1,n,1}\circ(x\otimes \id_{A_n}\otimes x)+\\
 &\quad +(-1)^{n+2}\m_{1,n,1}\circ(x\otimes \id_{A_n}\otimes x)+
 -\m_{n,1,1}\circ(\id_{A_n}\otimes x\otimes x)=0\,.
\end{split}
\end{equation*}
The assumption that $\hat A$ is a bialgebra in $\C^\I$ and the
bi-invariance property of $x$ according to the second part of the
proposition verify directly that $\bigl((A_0,A_1),[x,\cdot]_0,\Delta_0,
\Delta_1\bigr)$ is a bialgebra in $\C^{{\mathbf 2}_c}$. The application of
Proposition \ref{Bialg-ext} then completes the proof.\end{proof}
\abs
In the remainder of this section we will look in more detail to first order
differential calculi providing the corresponding facts which are
necessary for the understanding of the final section on exterior differential
Hopf algebras. We will begin with the study of the category of derivations
$\Der(A)$ over an algebra $A$ in $\C$.
\abs
\begin{definition}\label{deriv}
Let $A$ be an algebra in $\C$. Then the category $\Der (A)$ of
derivations over $A$ consists of objects $(A,X,\dif_X)$ which are
algebras in $\C^{\I_c}$, and of algebra morphisms in $\C^{\I_c}$
of the form $(\id_A,f):(A,X,\dif_X)\to (A,Y,\dif_Y)$.
\end{definition}
\abs
A morphism in $\Der (A)$ can therefore equivalently be seen as an
$A$-bimodule morphism $f:X\to Y$ such that $f\circ\dif_X=\dif_Y$.
As a side-step the following proposition and corollary will be derived
using derivations. The statements are braided analogues of the
corresponding results of \cite{Wor}.
\abs
\begin{proposition}\label{26ab-D27}
Given an algebra $(A,\m,\eta)$ in $\C$. We consider $A_{(2)}:=\Ker(\m)$
and the unique morphism $D_A:A\to A_{(2)}$ such that
$\ker(\m)\circ D_A= \eta\otimes\id_A - \id_A\otimes\eta$.
Then $(A,A_{(2)},D_A)$ is a first order differential calculus. It is an
initial object in $\Der(A)$. If $B$ is a bialgebra in $\C$ with left
exact functor $B\otimes\id_\C$ then $(B,B_{(2)},D_B)$ is a bicovariant
differential calculus over $B$.
\end{proposition}

\begin{proof}
One can easily verify within the category $\Der(A)$ that
\begin{equation*}
(A,A_{(2)},D_A)=\Ker\bigl\{
       (A,A\otimes A,\eta\otimes\id_A - \id_A\otimes\eta)
       @>(\id_A,\m)>>
       (A,A,0) \bigr\}
       \qquad\in\Ob(\Der(A))\,.
\end{equation*}
The multiplication $\m$ is an epimorphism and therefore
$\m=\coker(\id_{A\otimes A}-\m\otimes\eta)$ because
$f\circ(\id_{A\otimes A}-\m\otimes\eta)=0$ implies that $f$ factorizes over
$\m$ in the form $f=f\circ(\id\otimes\eta)\circ\m$. So
$A\langle D_A\rangle=\Im(\id_{A\otimes A}-\m\otimes\eta)
=\Ker\m=A_{(2)}$ and $(A,A_{(2)},D_A)$ is a differential calculus.
By Proposition \ref{34-D34} there is at most one morphism in $\Der(A)$
from the differential calculus $(A,A_{(2)},D_A)$ to any other object
$(A,X,\dif_X)$. We have
\begin{equation*}
\big(\mu^X_\ell\circ(\id\otimes\dif_X)+\mu^X_r\circ(\dif_X\otimes\id)\big)
\circ\ker\m =\dif_X\circ\m\circ\ker\m=0\,.
\end{equation*}
Denoting $\pi:=\mu_\ell^X\circ(\id\otimes\dif_X)\circ\ker\m=
-\mu_r^X\circ(\dif_X\otimes\id)\circ\ker\m$ one can directly verify that
$(\id_A,\pi)\,:\,(A,A_{(2)},D_A)\to(A,X,\dif_X)$ is a morphism in $\Der(A)$.
\end{proof}
\abs
In a similar manner as in \cite{Wor} we can describe bicovariant first
order differential calculi over a Hopf algebra $H$ as certain
sub-bimodules of $H_{(2)}$. The corresponding result can be derived as a
corollary of Proposition \ref{26ab-D27} using Lemma \ref{univ-diff-calc}.
\abs
\begin{corollary}\label{class-bicov}
The bicovariant first order differential calculi over a flat Hopf algebra
$H$ with isomorphic antipode in $\C$ are in one-to-one correspondence
with either
\begin{enumerate}
\item Hopf sub-bimodules of $H_{(2)}$,
\item crossed sub-bimodules of $\ker\varepsilon$
\end{enumerate}
where the counit $\varepsilon$ is regarded as crossed module morphism
$\varepsilon:H^{\rm ad}\to\E$.
\end{corollary}

\begin{proof}
Let $A$ be any algebra in $\cal C$ such that $A\otimes\id_\C$ is
right exact. Let $(A,X,\dif_X)$ be a differential calculus and
$(\id_A,f)\,:\,(A,X,\dif_X)\to(A,A_{(2)},D_A)$ be a morphism in $\Der(A)$.
Then $f$ is an epimorphism because $f\circ \dif_X=D_A$ and
$(A,A_{(2)},D_A)$ is a differential calculus.
Conversely, if $f:A_{(2)}\to X$ is a bimodule epimorphism
then $(A,X,f\circ D_A)$ is a differential calculus.
So first order differential calculi over $A$ are in one-to-one
correspondence with sub-bimodules $\Ker f$ of $A_{(2)}$.
One can specify this result to the Hopf algebra $H$ which is an algebra in
${}^H\C^H$ through its comultiplication. Then the first statement of the
corollary is derived.

Lemma \ref{univ-diff-calc} and the equivalence
$\hhchh\simeq\DY({\cal C})^H_H$ from Theorem \ref{yd-hopfbi} imply the
correspondence between $H_{(2)}\in\Ob\bigl(\hhchh\bigr)$ and
$\Ker\varepsilon\in\Ob\bigl(\DY({\cal C})^H_H\bigr)$,
and therefore between Hopf sub-bimodules of $H_{(2)}$ and
crossed sub-bimodules of $\Ker\varepsilon$.\end{proof}
\abs
For our further investigations it is convenient to work with
the notion of a comma category in the simplest case. Its definition is
outlined in the Appendix.
In Proposition \ref{27-D28} the comma category $(\E\downarrow {\C})$
will be used to formulate the braided
Woronowicz construction of a Hopf bimodule with bi-invariant
element out of a first order bicovariant differential calculus.
Before, we derive some results relating bimodules in $(\E\downarrow {\C})$
and derivations.
\abs
\begin{lemma}\label{comma-alg-mod}
Let $(A,\m,\eta)$ be an algebra in $\C$ then
$\big((A,\eta),\m,\eta\big)$ is an algebra in $(\E\downarrow {\C})$.
An analogous result holds for bi- and Hopf algebras.
A module $(X,\mu)$ over the algebra $A$ in $\C$ for which
$(X,x)$ is an object in $(\E\downarrow {\C})$ is also an
$(A,\eta)$-module $((X,x),\mu)$ in $(\E\downarrow {\C})$.
\end{lemma}

\begin{proof}
Of course, trivially $\m\circ(\eta\otimes\eta)=\eta$, $\eta\circ\id_\E=
\eta$ and similarly for the comultiplication, counit and antipode
of bi- and Hopf algebras. Therefore $\big((A,\eta),\m,\eta\big)$ is an
algebra in $(\E\downarrow {\C})$. Analogously the equation
$\mu_l\circ(\eta\otimes x)=x$ shows that $\big((X,x),\mu_l\big)$ is
an $(A,\eta)$-left module in $(\E\downarrow {\C})$.\end{proof}
\abs
Restricting to $\I={\mathbf 2}$ one finds a similar result as in
\cite{Wor} which is converse to Proposition \ref{22-D22}.
\abs
\begin{proposition}\label{23-D22}
Given a derivation $(X,\dif)$ over an algebra $A$ in $\C$. Then an
$(A,\eta)$-bimodule $A\oplus_\dif X:=\left(
 \left(A\oplus X,\binom{\eta}{0}\right),\mu_{\dif,l},\mu_{\dif,r}\right)$
in $(\E\downarrow {\C})$ can be constructed where the modul actions are
given through
\begin{equation}\label{comma-mod-act}
\mu_{\dif,l}=
\left(\begin{matrix}
       \m &  0\\
        0  & \mu_l
      \end{matrix}\right)
  \quad {\it and}\quad
\mu_{\dif,r}=
\left(\begin{matrix}
   \m &  0\\
   \mu_l\circ(\id_A\otimes\dif)  & \mu_r
      \end{matrix}\right)\,.
\end{equation}
Therefore in particular
$\left(A,A\oplus_\dif X,[\binom{\eta}{0},\cdot]\right)$ is an algebra in
$\C^{{\mathbf 2}_c}$ according to Proposition \ref{22-D22}.
Moreover $\left[\binom{\eta}{0},\cdot\right]={0\choose\dif}$.
\end{proposition}

\begin{proof}
For as the proof of the proposition does not require difficult techniques
we exemplarily outline that $\mu_{\dif,r}$ is a right module action. It holds
\begin{equation}\label{comma-mod-unit}
 \mu_{\dif,r}\circ(\id_{A\oplus_\dif X}\otimes\eta)
 =\left(
  \begin{matrix}
  \m                           & 0    \\
  \mu_l\circ(\id_A\otimes\dif) & \mu_r
  \end{matrix}
  \right)\circ
  \left(
  \begin{matrix}
  \id_A\otimes\eta & 0               \\
  0                & \id_X\otimes\eta
  \end{matrix}
  \right)
 = \id_{A\oplus_\dif X}
\end{equation}
and
\begin{equation}\label{comma-mod-mult}
\begin{split}
\mu_{\dif,r}\circ(\mu_{\dif,r}\otimes\id_A)
&=\left(
  \begin{matrix}
  \m                           & 0    \\
  \mu_l\circ(\id_A\otimes\dif) & \mu_r
  \end{matrix}
  \right)\circ
  \left(
  \begin{matrix}
  \m\otimes\id_A                           & 0    \\
  \mu_l\circ(\id_A\otimes\dif)\otimes\id_A & \mu_r\otimes\id_A
  \end{matrix}
  \right)\\
&=\left(
  \begin{matrix}
  \m\circ(\id_A\otimes\m) & 0\\
  \mu_l\circ(\id_A\otimes\dif\circ\m) & \mu_r\circ(\id_X\otimes\m)
  \end{matrix}
  \right)\\
&=\mu_{\dif,r}\circ(\id_{A\oplus_\dif X}\otimes\m)\,.
\end{split}
\end{equation}  
The equations (\ref{comma-mod-unit}) and (\ref{comma-mod-mult}) prove the
statement. The relation 
$\mu_{\dif,r}\circ\left(\binom{\eta}{0}\otimes\eta\right)=\binom{\eta}{0}$,
which is a consequence of (\ref{comma-mod-unit}), implies that
$\mu_{\dif,r}$ is a morphism in $(\E\downarrow {\C})$.\end{proof}
\abs
The results of Propositions \ref{22-D22} and \ref{23-D22}
are used to formulate the following theorem which shows that the categories
${}_{(A,\eta)}(\E\downarrow{\C})_{(A,\eta)}$ and $\Der(A)$ are
mutually adjoint.
\abs
\begin{theorem}\label{24-D24}
Let $A$ be an algebra in $\cal C$. Then the functors
\begin{equation}\label{adjoint1}
\bfig
\puttwohmorphisms(0,0)[\Der(A)`{}_{(A,\eta)}(\E\downarrow\C)_{(A,\eta)};%
 \scriptstyle{A\oplus_\dif (\_\,)}`\scriptstyle{{[\_\,]}{}_A}]{1000}1{-1}
\efig
\end{equation}
are defined on the objects by $A\oplus_\dif(X,\dif)= A\oplus_\dif X$ and
$[(X,x)]{}_A= (X,[ x,\cdot ])$ and on morphisms by
$A\oplus_\dif(f)=\id_A\oplus f$ and $[ g ]{}_A=g$ respectively.
The functor $A\oplus_\dif(\_\,)$ is left adjoint of the functor $[\_\,]_A$
through the (inverse of the) natural bijection
$\Gamma_{(X,\dif),(Y,y)}:
{\Hom}_{\Der(A)}\big((X,\dif),(Y,[y,\cdot])\big)\to
\Hom_{{}_{(A,\eta)}(\E\downarrow{\cal C})_{(A,\eta)}}\big(A\oplus_\dif X,
(Y,y)\big)$, $f\mapsto \big(\mu_l^Y\circ(\id_A\otimes y), f\big)$.
\end{theorem}

\begin{proof}
Straightforward calculations like in Proposition \ref{23-D22} show
that $\big(\mu_l^Y\circ(\id_A\otimes y), f\big)$ is an $A$-bimodule morphism.
Moreover, for any
$(g,f)\in\Hom_{{}_{(A,\eta)}(\E\downarrow{\cal C})_{(A,\eta)}}
\big(A\oplus_\dif X,(Y,y)\big)$ the 
condition $(g,f)\circ\binom{\eta}{0}=y$ is
equivalent to $g=\mu_l^Y\circ(\id_A\otimes y)$.\end{proof}
\abs
\begin{lemma}\label{26c-D27}
Suppose that $B$ is a bialgebra in $\C$ with left exact functor
$B\otimes\id_\C$. The full subcategory of $\Der(B)$ generated
by left, right or bicovariant first order differential calculi
over $B$ is a subcategory of the left, right or bicovariant
first order differential calculi over $B$ respectively.
We denote it by $\LCD(B)$, $\RCD(B)$ and $\BiCD(B)$ respectively.
\end{lemma}

\begin{proof}
It only has to be proved that the morphisms are (bi-)comodule
morphisms. We demonstrate the case of left comodules. Suppose that
$f:(X,\dif_X)\to(Y,\dif_Y)$ is a morphism in $\LCD(B)$. Then the relations
in Figure 3 are fulfilled.
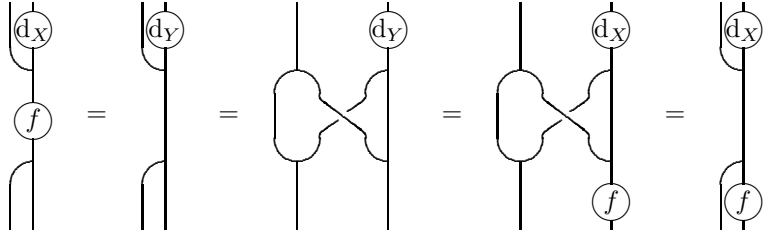
\begin{figure}
\begin{center}
\unitlength=0.4ex  
\begin{picture}(165,50)

\put(0,0){\figeins}
\put(19,25){\makebox(0,0){$=$}}
\put(29,0){\figzwei}
\put(48,25){\makebox(0,0){$=$}}
\put(63,0){\figdrei}
\put(97,25){\makebox(0,0){$=$}}
\put(112,0){\figvier}
\put(146,25){\makebox(0,0){$=$}}
\put(156,0){\figfuenf}
\end{picture}
\end{center}
\caption{Proof of the comodule morphism property.}
\end{figure}  
In the first equation of Figure 3 we used that $f$ is a
$B$-bimodule morphism and the identity $f\circ\dif_X=\dif_Y$. The
second equality is obtained by applying the obvious Hopf module properties
of a left covariant differential calculus. With the same
techniques as exploited in the first and second equations we derive the
equalities three and four. Since by assumption
$\mu_l^X\circ(\id_B\otimes\dif_X)$ is an epimorphism we conclude from
Figure 3 that $f$ is a $B$-left comodule morphism.\end{proof}
\abs
Finally Lemma \ref{26c-D27} yields the proposition which is a formal
Woronowicz construction of Hopf bimodules with invariant elements.
\abs
\begin{proposition}\label{27-D28}
Let $H$ be a flat Hopf algebra with isomorphic antipode in $\C$.
Then the following functors can be defined canonically.
\begin{equation}\label{diff-comma}
\begin{alignedat}{3}
H\oplus_\dif (\_\,) &: \LCD(H)&\ &\longrightarrow &\ 
{}^{(H,\eta)}_{(H,\eta)}(\E&\downarrow{\C})_{(H,\eta)}\,,\\
H\oplus_\dif (\_\,) &: \RCD(H)&\ &\longrightarrow &\ 
{}_{(H,\eta)}(\E&\downarrow{\C})_{(H,\eta)}^{(H,\eta)}\,,\\
H\oplus_\dif (\_\,) &: \BiCD(H)&\ &\longrightarrow &\ 
{}^{(H,\eta)}_{(H,\eta)}(\E&\downarrow{\C})_{(H,\eta)}^{(H,\eta)}
\end{alignedat}
\end{equation}
where the coactions on the objects $H\oplus_\dif(X,\dif_X)$ are 
given by
\begin{math}
\nu_{\dif,l}=\left(
\begin{smallmatrix}\Delta &0\\
                     0    &\nu_l
\end{smallmatrix}\right)
\end{math} 
and
\begin{math}
\nu_{\dif,r}=\left(
\begin{smallmatrix}\Delta &0\\
                     0    &\nu_r
\end{smallmatrix}\right).
\end{math}
\end{proposition}

\begin{proof}
With the help of Proposition \ref{22-D22}, Theorem \ref{24-D24} and
Lemma \ref{26c-D27} the proof of the proposition follows
straightforwardly.\end{proof}
\abs
\begin{remark}{\rm In particular for a bicovariant differential calculus
$(X,\dif)$ over $H$ this means that $\binom{\eta}{0}:\E\to H\oplus_\dif X$
is bi-invariant. I.\ e.\
$\nu_l\circ\binom{\eta}{0}=\eta\otimes\binom{\eta}{0}$
and $\nu_r\circ\binom{\eta}{0}=\binom{\eta}{0}\otimes\eta$.}
\end{remark}
\abs
\begin{remark}{\rm We have been working
with Hopf algebras which have an isomorphic antipode
in $\C$. This condition could be omitted at several places.
Then instead one would obtain {\it pseudo-braided} categories of Hopf
bimodules. Pseudo-braided categories are pre-braided categories
where the pre-braiding is an isomorphism if one of the tensor factors
is the unit object. Many results obtained in the article could be
transferred to this case without difficulties.}
\end{remark}
\abs
\subsection*{Exterior Differential Hopf Algebras}

In this section we construct the braided exterior Hopf algebra of
differential forms of a given first order bicovariant differential
calculus. This object is a differential Hopf algebra calculus and
extends the initial first order differential calculus uniquely.

In Lemma \ref{7-D40} we build the differential algebra
$\big((H\oplus_\dif X)^{\wedge_H},\hat\Dif\big)$ out of the bicovariant
differential calculus $(X,\dif)$ over $H$ according to Propositions
\ref{27-D28}, \ref{22-D22} and Definition \ref{ext-alg}. Then in Theorem
\ref{8-10-D40} the exterior algebra $(X^{\wedge_H},\dif^\wedge)$ is
derived as the maximal calculus of
$\big((H\oplus_\dif X)^{\wedge_H},\hat\Dif\big)$ .
\abs
\begin{lemma}\label{7-D40}
Suppose that $(X,\dif)$ is a first order bicovariant differential calculus
over the flat Hopf algebra $H$ with isomorphic antipode in the category $\C$.
Then the graded Hopf algebra $(H\oplus_\dif X)^{\wedge_H}$ is a differential
algebra in ${\C}^{\uNN_c}$ through the differential
$\hat\Dif:=[\binom{\eta}{0},\cdot]=\binom{0}{\dif}$ which is derived from
$\binom{\eta}{0}:\E\to H\oplus X$ according to Proposition \ref{22-D22}.
\end{lemma}

\begin{proof}
We have to prove that 
$\m_{2,n}\circ\big(\m_{1,1}\circ(x\otimes x)\otimes\id_n)=
\m_{n,2}\circ\big(\id_n\otimes\m_{1,1}\circ(x\otimes x))$ for all $n\in\I$,
where $\m:=\m^{(H\oplus_\dif X)^{\wedge_H}}$ and $x:=\binom{\eta}{0}$.
Then all the conditions of the first part of Proposition \ref{22-D22}
are satisfied.

We are using Theorem \ref{bialg-proj} and the relation corresponding to
(\ref{mult-ext-hopf}) involving the multiplication
$\m^{(H\oplus_\dif X)^{\wedge_H}}$. Then we derive
\begin{equation}
\m^{(H\oplus_\dif X)^{\wedge_H}}_{1,1}\circ
\left(\binom{\eta}{0}\otimes\binom{\eta}{0}\right)=
\coim(\hat A)_2\circ \left(\binom{\eta}{0}\otimes{}_{\scriptscriptstyle
H\oplus_\dif X}\p\circ\binom{\eta}{0}\right)=0
\end{equation}
because
\begin{equation}
\begin{gathered}
{}^{\scripthhchh}\Psi_{H\oplus_\dif X,H\oplus_\dif X}\circ
  \left(\binom{\eta}{0}\otimes{}_{\scriptscriptstyle H\oplus_\dif X}\p
  \circ\binom{\eta}{0}\right)
 =\left(\binom{\eta}{0}\otimes{}_{\scriptscriptstyle H\oplus_\dif X}\p
  \circ\binom{\eta}{0}\right)\,,\\
(\hat A_{H\oplus_\dif X})_2 = \id_{{\cal T}_{\C}^H(H\oplus_\dif X)_2}-
 {}^{\scripthhchh}\Psi_{H\oplus_\dif X,H\oplus_\dif X}\,.
\end{gathered}
\end{equation}
This concludes the proof.\end{proof}
\abs
\begin{remark}
{\rm It can be shown that the object
$\big((H\oplus_\dif X)^{\wedge_H},\hat D\big)$ is even a differential
Hopf algebra. However, this fact will not be used in the sequel.
Therefore we only mention it here in the remark.}
\end{remark}
\abs
We exploit Lemma \ref{7-D40} to prove the central theorem
in Chapter \ref{sec-diff-calc}.
\abs
\begin{theorem}\label{8-10-D40}
Let $H$ be a flat Hopf algebra with isomorphic antipode in $\C$,
and $(X,\dif)$ be a first order bicovariant differential calculus
over $H$. Then there exists a differential $\dif^\wedge$ such that
$(X^{\wedge_H},\dif^\wedge)$ is a differential Hopf algebra calculus in
$\C^{\uNN_{c}}$ which extends the first order calculus $(X,\dif)$,
i.\ e.\ $\dif^\wedge_0=\dif$.

We call $(X^{\wedge_H},\dif^\wedge)$ the braided exterior tensor Hopf
algebra of differential forms over $(X,\dif)$.
Explicitely, $(X^{\wedge_H},\dif^\wedge)=
\subdif{\big((H\oplus_\dif X)^{\wedge_H},\hat D\big)}$.
\end{theorem}

\begin{proof}
Since $\binom{0}{\id_X}:X\to H\oplus_\dif X$ is a Hopf bimodule monomorphism,
Proposition \ref{5-D11XVII} yields a monomorphic Hopf algebra morphism
$\hat\j:X^{\wedge_H}\to(H\oplus_\dif X)^{\wedge_H}$.
According to Proposition \ref{max-diff} the 0th and the 1st component of the
embedding of the maximal differential calculus into
$((H\oplus_\dif X)^{\wedge_H},\hat \Dif)$ are given by $\i_0=\id_H$
and $\i_1=\im (\m^{(H\oplus_\dif X)^{\wedge_H}}\circ(\id_H\otimes \Dif_0))
=\binom{0}{\im (\mu_l\circ(\id_H\otimes\dif))}=\binom{0}{\id_X}\cong\j_1$
where it has been used that the left action of $H$ on $H\otimes_\dif X$
is diagonal and that $(X,\dif)$ is a differential calculus.
Because $\hat \i$ and $\hat\j$ are algebra morphisms, and both
$X^{\wedge_H}$ and $\subdif{(H\otimes_\dif X)^{\wedge_H}{}}$ are generated
multiplicatively by their 0th and 1st component, it follows that
$\hat\i=\hat\j$. Hence $(X^{\wedge_H},\dif^{\wedge})=
\subdif{\big((H\oplus_\dif X)^{\wedge_H},\hat \Dif\big)}$ and
$\dif^\wedge_0=\dif$. Since $\binom{\eta}{0}$ is bi-invariant by
Proposition \ref{27-D28} we can apply the Propositions \ref{22-D22} and
\ref{7-8-6ZD11} to finish the proof of the theorem.\end{proof}
\abs
Using Theorem \ref{23-8XIIID34} and Corollary
\ref{proof-2-D55} the following corollary can be deduced.

\begin{corollary}\label{2-D55}
Suppose that $(\hat Y,\hat\dif_Y)$ is a differential Hopf algebra calculus.
If $(g_0,g_1):(Y_0,Y_1,\dif_{Y,0})\to (H,X,\dif)$ is a bialgebra morphism in
$\C^{{\mathbf 2}_c}$ then there exists a unique bialgebra morphism
$g^{\wedge_H}:(\hat Y,\hat\dif_Y)\to (X^{\wedge_H},\dif^{\wedge})$ in
$\C^{\uNN_c}$ which extends $(g_0,g_1)$.\endproof
\end{corollary}
\abs
\begin{remark} {\rm In \cite{q-alg/9609005,q-alg/9702020} the 
intrinsic Hopf algebra structure of the higher order Woronowicz 
differential calculus over a first order bicovariant differential calculus
is shown to generate very naturally a bicovariant algebra
of four basic objects within a differential calculus on quantum groups
as the Heisenberg double of two mutually dual
graded Hopf algebras. Namely coordinate functions, differential 1-forms, 
Lie derivatives, and inner derivations.

This construction works completely analogous in the braided setting.
The definition of the Heisenberg double can be directly generalized to the
braided case if we impose some further conditions on the braided category $\C$
like, for instance, rigidity (see \cite{BKLT:int} for such a
possibility). A more detailed description of these additional properties 
of $\C$ is required for the rigorous mathematical treatment -- this will
be published elsewhere. Then the generalization is based on the 
$(((H^\vee)_{\rm op})^{\rm op})$-left module algebra structure on $H$.
And the covariance properties of the corresponding cross product algebra
$H\# (((H^\vee)_{\rm op})^{\rm op})$ are checked literally.

For dually paired Hopf algebras $\big((A,X),(A^*,X^*)\big)$ in 
$\C^{\mathbf 2}$ the external Hopf algebras $X^\wedge$, $X^{*\wedge}$ in 
$\C^\uNN$ are dually paired in a natural way. And the formal analogue of the 
Heisenberg double $X^\wedge \# X^{*\wedge}$ can be constructed.
Let moreover $(A,X,\dif)$ be a bicovariant differential calulus.
Note that in this case the ``Cartan identity'' is a special case of
the differential coalgbera property
$\Delta_{1,n}\circ\dif=
(\dif\otimes\id)\circ\Delta_{0,n}-(\id\otimes\dif)\circ\Delta_{1,n-1}$.
We simply have to switch the left ``coaction" $\Delta_{1,\_}$ of $X$ into
a right action of the dual $X^*$.}
\end{remark}
\abs
The results of our article show that it is possible to generalize 
Woronowicz's and others' results on (higher order bicovariant) differential 
calculi to a very general class of braided abelian categories.
The essential tools for the succesful approach are the braided
Hopf bimodules and crossed modules. In rigid braided monoidal categories
an appropriate definition of a (non-degenerate) pairing of objects
should basically lead to versions of Woronowicz's
``quantum Lie algebras" in braided categories and to a description of
differential calculi in terms of ``braided linear functionals". The
resulting ``quantum Lie algebras" may fit into the framework of Majid's
braided Lie algebras \cite{Ma7}. Once a suitable braided category is fixed
the ``functional" description of bicovariant differential calculi could
provide another tool for attacking problems like their classification
in a way similar as in \cite{BGMST,Ma8,Ros,SS}. This will be our direction
for further investigations on this subject.
\abs
\abs
\appendix
\section{Appendix}\label{appendix}
In the appendix we set out some general categorical results and
derivations of previous outcomes which will be needed in the course of
the present work or are related with it.
\abs
\begin{proposition}\label{21-22-37--6-1ZD11}
Let $\C$ be a braided abelian category with bi-additive
tensor product. Let $A$ be an algebra in $\C$ and $A\otimes\id_{\C}$
be right exact. Then the categories of modules ${}_A{\C}$, ${\C}_A$
and the category of bimodules ${}_A{\C}_A$ are abelian.
The dual result holds for the comodules of a coalgebra
$C$ if the functor $C\otimes\id_{\C}$ is supposed to be left exact.
If $B$ is a flat bialgebra in $\C$, then also
the several categories of Hopf (bi-)modules
${}^B_B{\C}$, ${\C}_B^B$, ${}_B^B{\C}_B$,
${}_B{\C}_B^B$,$\ldots$,${}_B^B{\C}_B^B$ and $\DY({\C})^B_B$
are abelian.
Let $H$ be a flat Hopf algebra with invertible antipode in $\C$ then
the categories $\hhchh$ and $\DY({\C})^H_H$ are braided abelian with
bi-additive tensor product. In this case in particular the equivalences
$H\ltimes(\_\,):\DY({\C})^H_H\leftrightarrows\hhchh:{}_H(\_\,)$ of Theorem
\ref{yd-hopfbi} and the several forgetful functors, especially
$\hhchh\longrightarrow {\C}$, $\DY({\C})^H_H\longrightarrow{\C}$
are exact.
\end{proposition}

\begin{proof}
We have $\id_A\otimes\coker f =\coker(\id_A\otimes f)$ when $A\otimes\id_\C$
is right exact. Then for any $A$-left module morphism $f:X\to Y$ the
cokernel is an $A$-module morphism because there is a unique morphism
$\mu_l'$ such that $\mu_l'\circ(\id\otimes\coker f)=\mu_l'\circ
\coker(\id\otimes f)=\coker f\circ\mu_l^Y$. The module properties can be
verified without problems. The kernel of an $A$-module morphism is an
$A$-module morphism even if $A\otimes\id_\C$ is not supposed to be right
exact. Then one proves easily that the
$H$-left, $H$-right and $H$-bimodule categories ${}_H{\C}$,
${\C}_H$ and ${}_H{\C}_H$ respectively are abelian if $\C$
is abelian. If $B$ is a flat bialgebra in $\C$ it follows immediately
that ${}^B_B{\C}$, ${\C}_B^B$, ${}_B^B{\C}_B$,
${}_B{\C}_B^B$,$\ldots$,${}_B^B{\C}_B^B$ and $\DY({\C})^B_B$
are abelian. Suppose that $H$ is a flat Hopf algebra, then according to
Theorem \ref{Hopf-br} the tensor product of the braided abelian category
$\hhchh$ is bi-additive. Similarly it can be proven that
$\DY({\C})^H_H$ is braided abelian with bi-additive tensor product.
Obviously the forgetful functors $\hhchh\to {\C}$ and
$\DY({\C})^H_H\longrightarrow{\C}$ are exact. The flatness of $H$
implies the exactness of the functor $H\ltimes(\_\,)$. To prove that
${}_H(\_\,)$ is exact we consider an $H$-Hopf bimodule morphism $g$.
It holds $g\circ\i\circ\ker (\p\circ g\circ\i)
=0$ and $\p\circ g\circ\i\circ(\p\circ\ker g\circ\i)=0$
where $\i$ and $\p$ are the universal morphisms coming
from the construction of ${}_H(\_\,)$. Hence there are
unique morphisms $k:\Ker({}_Hg)\to\Ker g$ and
$l:{}_H(\Ker g)\to\Ker({}_Hg)$ such that 
$\ker({}_Hg)=\p\circ g\circ k$ and ${}_H(\ker g)=
\ker({}_Hg)\circ l$. Therefore $\p\circ k\circ\l=\id$ and
$l\circ\p\circ k=\id$ which yields $\ker({}_Hg)={}_H(\ker g)$ up to
isomorphism. The dual result for the cokernels will be derived in a
similar way.\end{proof}
\abs
\begin{corollary}\label{4-D11}
Let $\C$ be a braided abelian category with bi-additive tensor product.
Assume that $H$ is a flat Hopf algebra with invertible antipode in
$\C$. If $X$ is an $H$-Hopf bimodule and $X\otimes\id_{\C}$ is
left or right exact in $\C$, then $X\otimes_H\id_{\C}$ is
left or right exact respectively. Conversely if $X\otimes_H\id_{\C}$ is
left/right exact then ${}_HX\otimes\id_{\C}$ is left/right exact in
$\C$.
\end{corollary}

\begin{proof}
Form the proof of Proposition \ref{21-22-37--6-1ZD11} we know that
${}_H(\ker f)= {}_Y\p\circ \ker f\circ{}_X\i$ for every $H$-Hopf bimodule
morphism $f:X\to Y$. It follows that $X\otimes_H\id_\C$ is left/right
exact when $X\otimes\id_\C$ is left/right exact. Suppose that
$\id_\C\otimes_H X$ is left/right exact. Then Theorem
\ref{Hopf-br} implies that $\id_\C\otimes_H X$ is left/right
exact.\end{proof}
\abs
\begin{proposition}\label{4-D11XV--1-D39}
If $\C$ is $\otimes$-factor braided abelian and $H$ is a flat Hopf algebra
in $\C$ with isomorphic antipode, then $\hhchh$ and
$\DY({\C})^H_H$ are $\otimes$-factor braided abelian. The categories
$\hhchh$ and $\DY({\C})^H_H$ are $\otimes$-exact braided abelian
if $\C$ is $\otimes$-exact braided abelian.
\end{proposition}

\begin{proof}
For the proof of the $\otimes$-factor abelian property we remark that
the projection $\Pi:=\p\circ\i$ commutes with Hopf bimodule morphisms.
Now suppose that $\C$ is $\otimes$-exact abelian.
Let $f_1$ and $f_2$ be epimorphisms and $g$ and $h$ be morphisms
in $\hhchh$ such that $g\circ(\id\otimes_H f_2)=
h\circ(f_1\otimes_H\id)$. Using Theorem \ref{Hopf-br} we obtain
$g\circ(\id\otimes\p)\circ(\id\otimes f_2)=h\circ(\id\otimes\p')\circ
(f_1\otimes\id)$ where $\p$ and $\p'$ are the corresponding universal
morphisms of coinvariants. 
Since $\C$ is $\otimes$-exact abelian there is a unique
morphism $t$ in $\C$ obeying the relation
\begin{equation}
\begin{split}
g\circ(\id\otimes\p) &= t\circ(f_1\otimes\id)\,,\\
h\circ(\id\otimes\p')&= t\circ(\id\otimes f_2)\,.
\end{split}
\end{equation}
$f_2$ is an $H$-Hopf bimodule morphism. Thus there is the unique morphism
$t\circ(\id\otimes\i)$ which fulfills the equations
\begin{equation}
\begin{split}
g &= t\circ(\id\otimes\i)\circ(f_1\otimes_H\id)\,,\\
h &= t\circ(\id\otimes\i)\circ(\id\otimes_H f_2)\,.
\end{split}
\end{equation}
The morphisms $f_1\otimes_H\id$
and $\id\otimes_H f_2$ are epimorphic in $\C$. Therefore it follows
that $t\circ(\id\otimes\i)$ is a morphism in $\hhchh$. In an analogous
manner the pull-back condition of Definition \ref{tensor-exact} can be
proved for $\hhchh$. Similarly one verifies that $\DY({\C})^H_H$
is $\otimes$-exact abelian.\end{proof}
\abs
In the next proposition we derive a result which is useful for the study
of graded (co-)modules, Hopf bimodules and crossed modules.
\abs
\begin{proposition}\label{5-6-4ZD11}
Let $B$ be a bialgebra in the $\otimes$-factor braided abelian category $\C$
and suppose that the functor $B\otimes\id_\C$ is left exact. Then there are
canonical braided monoidal isomorphisms of comodule categories according to
\begin{equation}\label{cat-grade-mod1}
{}^B({{\C}^\I})\cong ({}^B{\C})^\I\quad{\it and}\quad
({}^B{\C^{\I_c}})\cong ({}^B{\C})^{\I_c}\,.
\end{equation}
Analogous isomorphisms are obtained for the categories of
right $B$-comodules and $B$-bicomodules. If $B$ is right exact the
corresponding results are obtained for the module categories.
Suppose that $H$ is a flat Hopf algebra with invertible antipode in $\C$,
then there are exact braided monoidal
isomorphisms of $\otimes$-factor braided abelian categories
\begin{equation}\label{cat-grade-mod2}
\begin{aligned}
{}^H_H({{\C}^\I})^H_H&\cong ({\hhchh})^\I\,,\\
\DY({\C}^\I)^H_H&\cong(\DY({\C})^H_H)^I\,,
\end{aligned}
\quad and \quad
\begin{aligned}
{}^H_H({\C^{\I_c}})^H_H &\cong ({\hhchh})^{\I_c}\,,\\
\DY(\C^{\I_c})^H_H&\cong(\DY({\C})^H_H)^{\I_c}\,.
\end{aligned}
\end{equation}
\end{proposition}

\begin{proof} We use (\ref{1-1b-1ZD11}) to imbed the catgory
$\C$ into ${\C}^\I$ or $\C^{\I_c}$. Then we consider the
functor $T:({}^B{\C})^{\uNN}\longrightarrow {}^B({\C}^{\uNN})$,
$(X_n,\nu_n)_n\mapsto ((X_n)_n,(\nu_n)_n)$, $\hat f\mapsto \hat f$. We
complete the proof using Propositions \ref{21-22-37--6-1ZD11},
\ref{4-D11XV--1-D39} and Corollary \ref{4-D11}. Similarly all other cases of
{(co-)}modules, Hopf bimodules and crossed modules can be treated.\end{proof}
\abs
A derivation of \cite{BD} is the following lemma.
\abs
\begin{lemma}\label{univ-diff-calc}
Let $H$ be a Hopf algebra with invertible antipode in the braided category
$\C$ with split idempotents,
and let $H^{\rm ad}$ be the $H$-right crossed module defined in Example
\ref{r-adj-cross}. Then the diagram of $H$-Hopf bimodule morphisms
\begin{equation}\label{hopf-bi-diag}
\bfig
\putsquare<1`1`1`1;1500`500>(0,-250)[H\otimes H`H\ltimes
 H^{\rm ad}`H`H\ltimes\E ;
 (\m\otimes\id_H)\circ(\id_H\otimes\Delta)`\m`H\ltimes(\varepsilon)`]
\efig
\end{equation}
is commutative. The horizontal morphisms are isomorphisms.
\end{lemma}

\begin{proof} Straightforward. We only note that the inverse of
$(\m\otimes\id_H)\circ(\id_H\otimes\Delta)$ is given by
$(\m\otimes\id_H)\circ(\id_H\otimes S\otimes\id_H)\circ
(\id_H\otimes\Delta)$.\end{proof}
\abs
We recall the definition and results on the comma category
$(\E\downarrow{\C})$ \cite{Mac}. It is the category which is formed
by objects $(X,x)$ where $X$ is an object in $\C$ and
$x:\E\to X$ is a morphism in $\C$. The morphisms 
$f:(X,x)\to(Y,y)$ in $(\E\downarrow{\C})$ are morphisms $f:X\to Y$
in $\C$ such that $f\circ x=y$. Without problems the next lemma can
be proved.
\abs
\begin{lemma}\label{comma-braid}
Suppose that $\C$ is braided with split idempotents. Then the category
$(\E\downarrow {\C})$ is braided with split
idempotents. The unit object is given by $(\E,\id_\E)$, the tensor
product of two objects equals
$(X,x)\otimes (Y,y)=(X\otimes Y,x\otimes y)$. The tensor product of two
morphisms $f$ and $g$ is given by the tensor product $f\otimes g$ in
$\C$. The braiding in $(\E\downarrow {\C})$ is defined through
$\Psi^{(\E\downarrow {\C})}_{(X,x),(Y,y)}=\Psi_{X,Y}$.\endproof
\end{lemma}
\abs
\abs

\Abs
\ouraddress

\begin{thebibliography}{MMMM}
\small
\bibitem[Art]{Art} E.\ Artin, {\it Theory of braids}, Ann.\ of Math.\ 
 {\bf 48}, 101 (1947).
 \nl
 {------}, {\it Braids and permutations},  Ann.\ of Math.\ {\bf 48},
 643 (1947).
\bibitem[BD1]{BD} Yu.\ {Bespalov} and B.\ {Drabant}, {\it Hopf
 (Bi-)Modules and Crossed Modules in Braided Monoidal Categories},
 to appear in J.\ Pure Appl.\ Algebra.
\bibitem[BD2]{BD2} Yu.\ {Bespalov} and B.\ {Drabant}, {\it Bicovariant
 Differential Calculi and Cross Products on Braided Hopf Algebras}, in:
 ``Quantum Groups and Quantum Spaces'', Banach Center Publications {\bf 40},
 Institute of Math., Polish Acad. Sci., eds.\ R.\ Budzynski, W.\ Pusz and
 S.\ Zakrzewski. Warsaw (1997).           
\bibitem[Bes]{Bes} Yu.\ {Bespalov}, {\it Crossed Modules and Quantum
 Groups in Braided Categories}, Appl.\ Categorical Structures {\bf 5}, No.\ 2
 (1997).
\bibitem[BGMST]{BGMST} F.\ Bonechi, R.\ Giachetti, R.\ Maciocco, E.\
 Sorace and M.\ Tarlini, {\it Quantum Double and Differential Calculi},
 preprint q-alg/9507004 (1995).
\bibitem[BKLT]{BKLT:int}
 Yu.\ Bespalov, T.\ Kerler, V.\ V.\ Lyubashenko and V.\ Turaev,
 {\it Integals for braided Hopf algebras}, in preparation.
\bibitem[BM]{BM} T.\ {Brzezinski}, {\it Crossed Products by a Coalgebra},
 preprint DAMTP/96-28 (1996).
 \nl
 S.\ {Majid}, {\it Diagrammatics of Braided Group Gauge Theory},
 preprint DAMTP/96-31 (1996).
\bibitem[Brz]{Brz} T.\ {Brzezinski}, {\it Remarks on Bicovariant
 Differential Calculi and Exterior Hopf Algebras}, Lett.\ Math.\ Phys.\
 {\bf 27}, 287 (1993).
\bibitem[CF]{CF} L.\ Crane and I.\ Frenkel, {\it Four
 dimensional topological quantum field theory, Hopf categories, and the
 canonical bases}, J.\ Math.\  Phys.\ {\bf 35}, 5136 (1994).
\bibitem[Dra]{Dr2} B.\ {Drabant}, {\it Differential Hopf Algebra Structure
 of the Quantum Standard Complex}, to appear in J.\ Math.\ Phys.
\bibitem[Dri]{Dri} V.\ G.\ {Drinfeld}, {\it Quantum groups}, Proceedings
 of the ICM, A.\ Gleason (ed.), Rhode Island, 798 (1987).
\bibitem[FRT]{FRT} L.\ D. Faddeev, N.\ Yu.\ Reshetikhin and L.\ A.\
 Takhtajan, {\it Quantization of Lie Groups and Lie Algebras}, Leningrad
 Math.\ J. {\bf 1}, 193 (1990).
\bibitem[FY]{FY} P.\ {Freyd} and D.\ {Yetter}, {\it Braided compact closed
 categories with applications to low dimensional topology}, Adv.\ Math.\
 {\bf 77}, 156 (1989).
\bibitem[Hus]{Hus} D.\ {Husemoller}, {\it Cyclic Homology}, Tata
 Institute Lecture Notes on Mathematics {\bf 83}, Springer (1991).
\bibitem[IV]{IV} A.\ P.\ {Isaev} and A.\ A.\ {Vladimirov},
 {\it $GL_q(n)$-Covariant Braided Differential Bialgebras}, Lett.\ Math.\
 Phys.\ {\bf 33}, 297 (1995).
\bibitem[JS]{JS} A.\ {Joyal} and R.\ {Street}, {\it Braided monoidal
 categories}, {Mathematics Reports 86008}, Macquarie University (1986).
\bibitem[Koo]{Koo} T.\ H.\ Koornwinder, {\it Orthogonal Polynomials in
 Connection with Quantum Groups}, NATO ASI Series {\bf C 294}, P.\ Nevai
 (ed.), Kluver Acad.\ Publishers (1990).
\bibitem[Lyu]{Lyu} V.\ V.\ Lyubashenko, {\it Tangles and Hopf algebras in
 braided categories}, J.\ Pure Appl.\ Algebra {\bf 98}, 245 (1995).
\bibitem[Ma1]{Ma1} S.\ {Majid}, {\it Algebras and Hopf algebras in
 braided categories}, Advances in Hopf Algebras, Lecture Notes in Pure
 and Appl.\ Math.\ {\bf 158}, 55, Dekker (1994).
\bibitem[Ma2]{Ma2} S.\ {Majid}, {\it Cross Products by Braided Groups and
 Bosonization}, J.\ Algebra {\bf 163}, 165 (1994).
\bibitem[Ma3]{Ma3} S.\ {Majid}, {\it Transmutation theory and rank for
 quantum braided groups}, Math.\ Proc.\ Camb.\ Phil.\ Soc.\ {\bf 113}, 45
 (1993).
\bibitem[Ma4]{Ma4} S.\ {Majid}, {\it Free Braided Differential Operators,
 Braided Binomial Theorem, and the Braided Exponential Map}, J.\ Math.\
 Phys.\ {\bf 34}, 4843 (1993).
 \nl
 ------, {\it $q$-epsilon tensor for quantum and braided matrices},
 J.\ Math.\ Phys.\ {\bf 36}, 1991 (1995).
\bibitem[Ma5]{Ma5} S.\ {Majid}, {\it Braided Momentum in the $q$-Poincar\'e
 Group}, J.\ Math.\ Phys. {\bf 34}, 2045 (1993).
\bibitem[Ma6]{Ma6} S.\ {Majid}, {\it Braided Groups}, J.\ Pure Appl.\
 Algebra {\bf 86}, 187 (1993).
\bibitem[Ma7]{Ma7} S.\ Majid, {\it Quantum and braided Lie algebras},
 J.\ Geom.\ Phys.\ {\bf 13}, 307 (1994).
\bibitem[Ma8]{Ma8} S.\ Majid, {\it Classification of Bicovariant
 Differential Calculi}, preprint q-alg/9608016 (1996).
\bibitem[Ma9]{Ma9} S.\ Majid, {\it Braided matrix structure of the
 Sklyanin algeba and of the quantum Lorentz group}, Commun.\ Math.\ Phys.\
 {\bf 156}, 607 (1993).
\bibitem[Mac]{Mac} S.\ {Mac Lane}, {\it Categories. For the Working
 Mathematician}, GTM {\bf 5}, Springer (1972).
 \nl
 ------, {\it Homology}, Grundlehren der mathematischen Wissenschaften
 {\bf 114}, Springer (1963).
\bibitem[Mal]{Mal} G.\ {Maltsiniotis}, {\it Le langage des espaces et des
 groupes quantiques}, Commun.\ Math.\ Phys.\ {\bf 151}, 275 (1993).
\bibitem[Man]{Man} Yu.\ I.\ Manin, {\it Notes on quantum groups and quantum
 de Rham complexes}, J.\ Theor.\ and Math.\ Phys.\ {\bf 92}, 425, Russian
 Academy of Sciences (1992).
\bibitem[MP]{MP} S.\ Mac Lane and R.\ Par\'e, {\it Coherence for Bicategories
 and Indexed Categories}, J.\ Pure Appl.\ Algebra {\bf 37}, 59 (1985).
\bibitem[Rad]{Rad} D.\ E.\ {Radford}, {\it The Structure of Hopf Algebras
 with a Projection}, J.\ Algebra {\bf 92}, 322 (1985).
\bibitem[Ros]{Ros} M.\ {Rosso}, {\it Alg\`ebres Enveloppantes
 Quantifi\'ees, Groupes Quantiques Compacts de Matrices et Calcul
 Differentiel Non Commutatif}, Duke Math.\ J.\ {\bf 61}, 11 (1990).
\bibitem[RT]{RT} D.\ E.\ {Radford} and J.\ {Towber},
 {\it Yetter-Drinfel'd categories associated to an arbitrary bialgebra},
 J.\ Pure Appl.\ Algebra {\bf 87}, 259 (1993).
\bibitem[RV]{q-alg/9702020}
 O.\ V.\ Radko and A.\ A.\ Vladimirov, {\it On the algebraic structure of
 differential calculus on quantum groups}, JINR preprint E2-97-45 (1997). 
 q-alg/9702020.
\bibitem[SS]{SS} K.\ Schm\"udgen and A.\ Sch\"uler, {\it Classification
 of Bicovariant Differential Calculi on Quantum Groups}, Commun.\ Math.\
 Phys.\ {\bf 170}, 315 (1995).
\bibitem[SZ]{SZ} M.\ Schlieker and B.\ Zumino, {Braided Hopf Algebras
 and Differential Calculus}, Lett.\ Math.\ Phys.\ {\bf 33}, 33 (1995).
\bibitem[Vla]{q-alg/9609005}
 A.\ A.\ Vladimirov, {\it Implications of the Hopf algebra properties
 of noncommutative differential calculi}, Czech. J. Phys. {\bf 47}, 131
 (1997). Proceedings of the 5th Colloquium on "Quantum Groups and
  Integrable Systems", (Prague, 1996).
\bibitem[Wor]{Wor} S.\ L.\ {Woronowicz}, {\it Differential Calculus on
 Compact Matrix Pseu\-do\-groups (Quantum Groups)},
 Commun.\ Math.\ Phys.\ {\bf 122}, 125 (1989).
\bibitem[Yet]{Yet} D.\ {Yetter}, {\it Quantum groups and representations of
 monoidal categories}, Math.\ Proc.\ Camb.\ Phil.\ Soc.\ {\bf 108}, 261
 (1990).
\end{thebibliography}
\end{document}